\documentclass[a4paper,12pt]{report} 

\raggedright                            

\usepackage{eucal}
\usepackage{calc}
\usepackage[all]{xy}
\usepackage[centertags]{amsmath}
\usepackage{latexsym}
\usepackage{amsfonts}
\usepackage{amssymb}
\usepackage{amsthm}
\usepackage[toc,page]{appendix}
\usepackage{subeqnarray}
\usepackage{makeidx}
\usepackage[english,french]{babel}
\usepackage{graphicx}
\usepackage{t1enc}
\usepackage{epsfig}
\usepackage{titlesec}
\usepackage{color}
\usepackage{url}
\usepackage{tocloft}
\makeindex
\newtheorem{remark}{Remarque}[chapter]
\newenvironment{prof}[1][\underline{Preuve}]{\textbf{#1.} }{\ \rule{0.5em}{0.5em}}
\theoremstyle{plain}
\newtheorem{lemma}{Lemme}[chapter]
\newtheorem{proposition}{Proposition}[chapter]
\newtheorem{theorem}{Théorème}[chapter]
\newtheorem{definition}{Définition}[chapter]
\newtheorem{corollary}{Corollaire}[chapter]

\newtheorem{example}{Exemple}[chapter]

\newtheorem{question}{Question}[chapter]
\theoremstyle{plain}

\renewcommand{\cftchappresnum}{\chaptername~} 

\title{\bf Actions des groupes topologiques sur les objects universels}  
\author{Brice Rodrigue MBOMBO\\
D\'epartement de Math\'ematiques\\ Facult\'e des
Sciences, Universit\'e de Yaoundé 1 (Cameroun)}              
\date{\today}                           
\begin{document}

\makeatletter
    \def\thebibliography#1{\chapter*{References\@mkboth
      {REFERENCES}{REFERENCES}}\list
      {[\arabic{enumi}]}{\settowidth\labelwidth{[#1]}\leftmargin\labelwidth
	\advance\leftmargin\labelsep
	\usecounter{enumi}}
	\def\newblock{\hskip .11em plus .33em minus .07em}
	\sloppy\clubpenalty4000\widowpenalty4000
	\sfcode`\.=1000\relax}
    \makeatother
\pagenumbering{roman}
\maketitle                              
\setcounter{page}{2}
\tableofcontents                        


\pagenumbering{arabic}

\chapter*{Dédicaces}\addcontentsline{toc}{chapter}{Dédicaces}
Je dédie ce travail à :
\begin{enumerate}
  \item [$\bullet$] A mes défunts grands-parents et à mon défunt papa MBOMBO Emmanuel. Que ceci soit un début d'accomplissement de
vos rêves.
  \item [$\bullet$] A Ma mère DONGMO Anne pour son soutien constant.
\end{enumerate}

\chapter*{Remerciements}\addcontentsline{toc}{chapter}{Remerciements}
En tout premier lieu, j'adresse mes très sincères remerciements aux Professeurs {\bf{Vladimir Pestov}} de l'université d'Ottawa (Canada) et {\bf{François Wamon}} de l'université de Yaoundé $1$ qui ont accepté de codiriger ma thèse.
Leur générosité, leur rigueur et leur précieux conseils ont été d'une importance primordiale pour l'accomplissement de ce travail.\\
Je suis très reconnaissant aux Professeurs {\bf{Marcel Dossa}} de l'université de Yaoundé $1$ et {\bf{Valentin Ferenczi}} de l'université de Sao Paulo (Brésil) d'avoir accepté d'écrire un rapport sur mon travail.

Je remercie également:
 \begin{enumerate}
  \item [$\bullet$] Le {\bf{MAÉCI}} (Ministère des Affaires Étrangères et du Commerce International) Canada qui par l'entremise du Bureau canadien de l'éducation internationale ({\bf{BCEI}}) m'a accordé une subvention de recherche pour une visite de recherches à l'université d'Ottawa.
   \item [$\bullet$] Le Natural Sciences and Engineering Research Council of Canada ({\bf{NSERC}}).
   \item [$\bullet$] Le {\bf{Département de Mathématiques et Statistiques de l'université d'Ottawa}} en général et le groupe de recherche en analyse en particulier pour l'hospitalité au cours de mes visites.
 \end{enumerate}
Mes très sincères remerciements vont également à tous les enseignants du département de Mathématiques de l'Université de Yaoundé I qui ont contribué à ma formation.

Je remercie très sincèrement mes amis du Département de Mathématiques et Informatique de la Faculté des Sciences de l'Université de Douala dont les nombreux conseils et discussions auront contribué de façon significative à la réalisation du présent travail.

Toutes les personnes qui me sont proches ont contribué chacun en sa manière propre à la réalisation de ce travail. Je pense à toute la grande famille {\bf{NJI TAMANJE}} pour sa patience et surtout pour toutes les privations consenties pendant les nombreuses années de préparation de cette thèse. Je pense très particulièrement à mon très grand ami {\bf{Bruno S. Djieutcheu}} et son épouse pour leurs nombreux conseils.\\









Enfin, je remercie {\bf{Dieu}} le père tout puissant qui nous donne le souffle de vie et sans qui aucune chose n'est possible.
\eject
\makeatletter
    \def\thebibliography#1{\chapter*{Références\@mkboth
      {REFERENCES}{REFERENCES}}\list
      {[\arabic{enumi}]}{\settowidth\labelwidth{[#1]}\leftmargin\labelwidth
	\advance\leftmargin\labelsep
	\usecounter{enumi}}
	\def\newblock{\hskip .11em plus .33em minus .07em}
	\sloppy\clubpenalty4000\widowpenalty4000
	\sfcode`\.=1000\relax}
    \makeatother                       

\tableofcontents                        


\chapter*{Résumé}
\addcontentsline{toc}{chapter}{Résumé}
Dans cette thèse, nous étudions l'existence des objects universels de deux types différents dans la théorie des groupes topologiques et leurs actions sur les espaces compacts.\\
Dans la première partie, nous contribuons au problème d'existence des espaces test pour la moyennabilité. Nous observons qu'un groupe polonais est moyennable si et seulement si toute action continue de celui-ci sur le cube de Hilbert poss\`ede une mesure de probabilit\'e invariante. Ceci g\'en\'eralise un r\'esultat de Bogatyi et Fedorchuk. Nous d\'emontrons \'egalement que les actions continues sur l'espace de Cantor permettent de tester la moyennabilit\'e, la moyennabilit\'e extr\^eme des groupes polonais non archimédiens, et la moyennabilit\'e \`a l'infini des groupes discrets d\'enombrables. Il en r\'esulte que cette derni\`ere propri\'et\'e peut \'egalement \^etre test\'ee par les actions sur le cube de Hilbert. Ces r\'esultats g\'en\'eralisent un crit\`ere de Giordano et de la Harpe. \\
Dans la deuxième partie, nous contribuons au problème d'existence d'un groupe topologique universel de poids non dénombrable. Nous montrons que le groupe des isométries de l'espace d'Urysohn-Kat\v etov de densité non dénombrable n'est pas universel pour la classe des groupes de poids non-dénombrable. Il ne contient pas le groupe des isométries de l'espace classique polonais universel d'Urysohn. Ce qui semble contradictoire avec le résultat d'Uspenkij ($1990$) sur l'universalité du groupe de isométries de l'espace classique polonais universel pour la classe des groupes polonais.
\subsection*{Mots cl\'{e}s:}
Moyennabilité, Moyennabilité extrême, Moyennabilité topologique, Espace test, Espace métrique d'Urysohn-Kat\v etov, Groupe topologique universel.
\chapter*{Abstract}\addcontentsline{toc}{chapter}{Abstract}
In this thesis, we study the existence of universal objects of two differents types in the theory of topological groups and their actions on compact spaces. \\
In the first part, we contribute to the problem of existence of test spaces for amenability. We observe that a Polish group is amenable if and only if every continuous action on the Hilbert cube admits an invariant probability measure. This generalizes a result of Bogatyi and Fedorchuk. We also show that actions on the Cantor space can be used to detect amenability and extreme amenability of Polish non-archimedean groups as well as amenability at infinity of discrete countable groups. As corollary, the latter property can also be tested by actions on the Hilbert cube. These results generalise a criterion due to Giordano and de la Harpe. \\
In the second part of this thesis we contribute to the problem of existence of an universal topological group of uncountable weight. We show that the group of isometries of a universal Urysohn-Kat\v etov metric space of uncountable density is not a universal group of the corresponding weight. For instance, it does not contain the group of isometries of the classical separable Urysohn metric space. This stands in sharp contrast with Uspenskij's $1990$ result about the group of isometries of the classical separable Urysohn metric space being a universal Polish group.

\section*{Keywords}
Amenability, Extreme amenability, Topological amenability, Test space, Urysohn-Kat\v etov metric space, Universal topological group.
\chapter*{Introduction Générale}
\addcontentsline{toc}{chapter}{Introduction Générale}
\pagenumbering{arabic}\setcounter{page}{1}
Cette thèse traite principalement de deux sujets :
les espaces test pour la {\em moyennabilité}, et le problème d'existence des {\em groupes topologiques universels}. On s'intéresse  dans la première partie aux espaces test pour la moyennabilité, la moyennabilité extrême et la moyennabilité topologique.\\
Les groupes moyennables ont été introduits en 1929 par J. Von Neumann \cite{von} dans son étude du paradoxe de Banach-Tarski \cite{bantark}:\\
Si $\mathcal{B}$ et $\mathcal{B}^{\prime}$ sont deux boules, n'ayant pas nécessairement la même taille, il est possible de
découper $\mathcal{B}$ en un nombre fini de morceaux, et réarranger ces morceaux (par des isométries directes) de
façon à ce qu'ils forment $\mathcal{B}^{\prime}$. Une version fantaisiste serait de dire qu'il est possible de découper un petit pois en morceaux et réarranger ces morceaux pour obtenir une boule de la taille du soleil.\\
Par extension de cette idée, un groupe $G$ admet une décomposition paradoxale s'il
existe des sous-ensembles de $G$ deux à deux disjoints $A_{1},A_{2},...,A_{n},B_{1}B_{2},...,B_{m}$ et
des éléments $g_{1}, g_{2}...,g_{n}, h_{1}, h_{2},...,h_{m}$ de $G$ tels que:
$$\underset{i=1}{\overset{n}{\bigcup}}g_{i}A_{i}=\underset{j=1}{\overset{m}{\bigcup}}h_{j}B_{j}.$$
Le groupe libre à deux générateurs $\mathbb{F}_{2}$ admet par exemple une décomposition paradoxale. Un groupe discret dénombrable $G$ est dit moyennable s'il n'admet pas de décomposition paradoxale.
On peut introduire la moyennabilité par plusieurs définitions équivalentes qui, en
apparence, sont assez éloignées les unes des autres: existence de moyennes invariantes; existence de mesures de probabilités invariantes; propriété du point fixe ; propriétés de Reiter et de
Reiter-Glicksberg ; propriétés de convolution dans $L^{p}$; propriétés combinatoires du type F{\o}lner,
etc...\\
Avant tout, la théorie consiste à prouver l'équivalence de ces diverses définitions. La classe des groupes moyennables contient les groupes finis, les groupes commutatifs, les groupes résolubles, le groupe unitaire $\mathcal{U}(\ell^{2})$ muni de la topologie forte (Gromov et  Milman\cite{grommil})...\\
Un groupe topologique $G$ est dit {\em moyennable} si toute action continue de $G$ sur un espace compact $X$ poss\`ede une mesure de probabilit\'e bor\'elienne invariante.\\
Un espace topologique compact métrisable $K$ est un espace test pour une classe $\mathcal{C}$ de groupes topologiques si
tout groupe $G\in \mathcal{C}$ est moyennable si et seulement si toute action continue de $G$ sur $K$ possède une mesure de probabilité
invariante.\\
En r\'eponse \`a une question de Grigorchuk, Giordano et de la Harpe \cite{gio} ont montr\'e  qu'un groupe discret d\'enombrable $G$ est moyennable si et seulement si toute action continue de $G$ sur l'ensemble de Cantor $D^{\aleph_{0}}$ poss\`ede une mesure de probabilit\'e invariante. On peut dire que l'ensemble de Cantor est un {\em espace test} pour la moyennabilit\'e des groupes discrets d\'enombrables.
Dans le m\^eme sens, Bogatyi et Ferdorchuk \cite{boga} ont r\'epondu \`a une question de \cite{gio} et demontr\'e que le cube de Hilbert $I^{\aleph_{0}}$ est \'egalement un {\em espace test} pour la moyennabilit\'e des groupes discrets d\'enombrables.\\
Dans cette thèse, nous d\'emontrons que le cube de Hilbert reste un espace test pour la moyennabilit\'e de tous les groupes polonais. Nous d\'emontrons \'egalement que le r\'esultat de Giordano et de la Harpe reste vrai pour les groupes polonais non archim\'ediens.\\
Un groupe topologique $G$ est dit {\em extr\^emement moyennable} si toute action continue de $G$ sur un espace compact poss\`ede un point fixe. Un
tel groupe non-trivial est moyennable, mais n'est jamais localement compact (Théorème de Veech \cite{veech}). Des exemples des groupes extrêmement moyennables sont nombreux et ils comprennent le groupe $Aut\,(X,\mu)$ des automorphismes mesurables préservant la mesure $\mu$ d'un espace borelien $(X,\mu)$ muni de la topologie faible (Giordano et Pestov \cite{giopes1}) et le groupe $Aut\,(\mathbb{Q},\leq)$ des bijections de $\mathbb{Q}$ dans lui-même qui préservent l'ordre muni de la topologie de la convergence simple (Pestov \cite{vp2}).\\
Nous d\'emontrons que l'ensemble de Cantor est un espace test pour la moyennabilit\'e extr\^eme des groupes polonais non archim\'ediens.
La question d'existence d'un espace test pour la moyennabilit\'e extr\^eme des groupes polonais reste ouverte, car ni l'espace de Cantor ni le cube de Hilbert ne possèdent cette propriété.\\
Si un groupe discret d\'enombrable $G$ op\`ere par hom\'eomorphismes sur un espace compact $X$, alors l'action de $G$ sur $X$ est {\em moyennable} s'il existe une suite $(b^{n})_{n\in\mathbb{N}}$ d'applications continues de $X$ dans $\mathbb{P}(G)$ telle que: $\underset{n\longrightarrow \infty}{\overset{}{\lim}}\,\,\underset{x\in X}{\overset{}{\sup}}\|gb^{n}_{x}-b^{n}_{gx}\|_{1}=0$ pour tout $g\in G$, o\`u $\mathbb{P}(G)$ d\'esigne l'espace des mesures de probabilit\'e sur $G$, muni de la topologie vague.
Un groupe discret d\'enombrable $G$ est dit {\em moyennable \`a l'infini} \cite{higson},\cite{anan}, ou {\em topologiquement moyennable}, s'il existe un espace compact $X$ et une action par hom\'eomorphismes de $G$ sur $X$ qui est moyennable. Les exemples des groupes moyennables \`a l'infini comprennent les groupes moyennables, les groupes d'automorphismes d'arbres enracin\'es, les groupes hyperboliques (en particulier, les groupes libres). Voir \cite{anan} et \cite{BO}, Ch. 5.\\
Par analogie avec le r\'esultat de Giordano et de la Harpe, nous d\'emontrons qu'un groupe discret d\'enombrable $G$ est moyennable \`a l'infini si et seulement si $G$ poss\`ede une action moyennable sur l'ensemble de Cantor ou sur le cube de Hilbert. Autrement dit, l'ensemble de Cantor et le cube de Hilbert sont des espaces test pour la moyennabilit\'e topologique des groupes discrets d\'enombrables.\\
Le chapitre $1$ regroupe les notions de bases sur les systèmes dynamiques topologique abstraits et des généralités sur la moyennabilité des groupes indispensables pour la compréhension de la suite de la thèse.\\
Au chapitre $2$, on utilise la décomposition du compactifié de Samuel équivariant $S(G)$ d'un groupe polonais $G$ en limite inverse d'un système de $G$-espaces compacts métrisables pour établir le résultat:
\begin{theorem}
Un groupe polonais $G$ est moyennable si et seulement si toute action continue de $G$ sur le cube de Hilbert $I^{\aleph_{0}}$ possède une mesure de probabilité borélienne invariante.
\end{theorem}
Si $G$ est de plus non archimédien, nous observons que nous pouvons décomposer le compactifié de Samuel équivariant $S(G)$ en limite inverse d'un système de $G$-espaces $X_{\alpha}$ où les $X_{\alpha}$ sont tous de Cantor. Cette observation est cruciale pour établir d'une part:
\begin{theorem}Un groupe polonais non archimédien $G$ est moyennable si et seulement si toute action continue de $G$ sur $D^{\aleph_{0}}$ possède une mesure de probabilité borélienne invariante.
\end{theorem}
et d'autre part:
\begin{theorem}
Un groupe polonais non archimédien $G$ est extrêmement moyennable si et seulement si toute action continue de
$G$ sur l'ensemble de Cantor $D^{\aleph_{0}}$ possède un point fixe.
\end{theorem}
Nous donnons une réponse partielle à la question d'existence d'un espace test pour la moyennabilité extrême des groupes polonais. Nous obtenons l'existence d'un espace test compact, séparable, mais non nécessairement métrisable pour la moyennabilité extrême des groupes polonais.\\
Grâce au théorème du point fixe de Schauder, nous observons que toute action continue d'un groupe monothétique ou d'un groupe solénoïde sur le cube de Hilbert $I^{\aleph_{0}}$ admet un point fixe.\\
Le chapitre $3$ est consacré à l'étude de la moyennabilité à l'infini. Nous y démontrons l'équivalence entre plusieurs caractérisations bien connues de la moyennabilité à l'infini. La décomposition du compactifié de Samuel équivariant $S(G)$ pour un groupe $G$ polonais et non archimédien obtenue au chapitre $2$, nous permet d'obtenir une preuve alternative du résultat suivant:
\begin{theorem}
Un groupe discret dénombrable $G$ est moyennable à l'infini si et seulement s'il admet une action moyennable sur l'ensemble de Cantor $D^{\aleph_{0}}$.
\end{theorem}
C'est par ailleurs le principal résultat de \cite{youssef}.
Comme corollaire et en utilisant le théorème de Keller( Toute partie métrisable, convexe, compacte et de dimension infinie d'un espace de Fréchet est homéomorphe au cube de Hilbert $I^{\aleph_{0}}$ \cite{bes}), on obtient:
\begin{theorem}\label{corohilbert}
Un groupe discret dénombrable $G$ est moyennable à l'infini si et seulement s'il admet une action moyennable sur le cube de Hilbert $I^{\aleph_{0}}.$
\end{theorem}

La deuxième partie de cette thèse est consacrée à l'étude du problème d'existence d'un groupe topologique universel pour la classe des groupes de poids non-dénombrable.\\
Soit $\mathcal{C}$ une classe d'espaces topologiques. $X\in \mathcal{C}$ est dit universel pour cette classe si pour tout $Y\in \mathcal{C}$, il existe un homéomorphisme de $Y$ sur un sous- espace de $X$.\\
L'espace de Cantor $D^{\aleph_{0}}$ par exemple est universel pour la classe des espaces topologiques métrisables séparables et de dimension $0$.\\
Un groupe topologique $G$ est universel pour une classe $\mathcal{C}$ de groupes topologiques si pour tout groupe topologique $H\in \mathcal{C}$, il existe un isomorphisme de groupes topologiques entre $H$ et un sous-groupe de $G$.
En réponse à une question de Ulam (cf. Problème $103$ dans \cite{ulam}), Uspenskij dans \cite{us5} a établi que le groupe $Iso(\mathbb{U})$ des isométries de l'espace d'Urysohn $\mathbb{U}$ sur lui-même muni de la topologie de la convergence simple est universel pour la classe des groupes métrisables et séparables. Quelques années avant, Uspenskij avait déja établi dans \cite{us4} que le groupe $Homeo(I^{\aleph_{0}})$ muni de la topologie compact-ouvert est universel pour la même classe de groupes topologiques.\\ Pendant longtemps il a été impossible de savoir si les deux groupes universels précédents sont isomorphes en tant que groupes topologiques. Les travaux de Pestov dans \cite{vp4} permettront de répondre à cette question par la négative. En effet, Pestov établit dans \cite{vp4} que le groupe $Iso(\mathbb{U})$ est extrêmement moyennable. En même temps, le groupe $Homeo(I^{\aleph_{0}})$ opère continûment sur l'espace compact $I^{\aleph_{0}}$ sans points fixes.\\
La question d'existence d'un groupe topologique universel pour la classe des groupes topologiques de poids non dénombrables reste ouverte.\\
Le chapitre $4$ est consacré à l'étude des sous-groupes du groupes des isométries de l'espace d'Urysohn-Kat\v etov de densité non-dénombrable $\mathbb{U}_{\mathfrak m}$.\\
Soit $\mathfrak m$ un cardinal infini vérifiant
\begin{equation}
\label{eq:mn}
\sup\left\{{\mathfrak m}^{\mathfrak n}\colon {\mathfrak n}<{\mathfrak m}\right\}={\mathfrak m}.
\end{equation}
Par exemple, $\aleph_{0}$, et tout cardinal fortement inaccessible $\mathfrak m$ vérifient la condition précédente.\\
Si $\mathfrak m$ est un cardinal infini vérifiant la condition précédente, alors il existe à isométrie près un unique espace métrique complet $\mathbb{U}_{\mathfrak m}$ de poids $\mathfrak m$ qui contient une copie isométrique de tout autre espace métrique de poids $\leq \mathfrak m$ et est
$ \mathfrak m$-homogène, autrement dit, toute isométrie entre deux sous-espaces métriques de densité $< \mathfrak m$ se prolonge en une isométrie globale de $\mathbb{U}_{\mathfrak m}$ sur lui-même. Par exemple, $\mathbb{U}_{\aleph_{0}}$ est l'espace classique d'Urysohn.\\
Un candidat naturel pouvant être universel dans la classe des groupes topologiques de poids $\mathfrak m>\aleph_{0}$ serait le groupe $Iso(\mathbb{U}_{\mathfrak m})$ des isométries de l'espace $\mathbb{U}_{\mathfrak m}$ (la version non-séparable de l'espace d'Urysohn construit par Kat\v etov \cite{kat}) sur lui-même muni de la topologie de la convergence simple. Mais de façon surprenante, nous observons dans ce chapitre que ce n'est pas le cas.\\
Un groupe topologique possède la propriété (OB)(\cite{ros}) si toutes les orbites d'une action continue de $G$ par isométries sur un espace métrique $(X,d)$ sont bornées.\\
 Après avoir présenté la construction de Kat\v etov de l'espace $\mathbb{U}_{\mathfrak m}$, nous caractérisons les voisinages de l'élément neutre dans un sous-groupe $G$ de $Iso(\mathbb{U}_{\mathfrak m})$ de densité $< \mathfrak{m}$. Cette caractérisation nous permet d'établir:
\begin{theorem}
Tout sous-groupe $G$ de $Iso(\mathbb{U}_{\mathfrak m})$ possédant la propriété (OB) et ayant une densité $<\mathfrak m$ est nécessairement FSIN: toute fonction bornée uniformément continue à gauche est uniformément continue à droite.
\end{theorem}
En particulier, si $G$ est métrisable ou localement connexe, alors $G$ est SIN: les structures uniformes gauche et droite coïncident sur $G$ .
Ceci est une restriction sérieuse qui montre qu'en particulier, pour un cardinal non-dénombrable $\mathfrak m$, le groupe $Iso(\mathbb{U}_{\mathfrak m})$ ne contient pas une copie isométrique du groupe $Iso(\mathbb{U})$.\\
D'autre part, nous faisons quelques observations dans le sens de caractériser les sous-groupes topologiques de $Iso(\mathbb{U}_{\mathfrak m})$. Dans cette direction, nous faisons l'observation suivante:
 \begin{theorem}
Tout groupe métrisable SIN de poids $\leq \mathfrak m$ se plonge dans $Iso(\mathbb{U}_{\mathfrak m})$.
\end{theorem}
La preuve de ce théorème utilise la même technique que Uspenkij (\cite{us5},\cite{us6}) pour établir l'universalité du groupe $Iso(\mathbb{U})$ pour la classe des groupes topologiques vérifiant le deuxième axiome de dénombrabilité.
En même temps, nous ignorons si tous les groupes SIN de poids $\leq \mathfrak m$ se plongent dans $Iso(\mathbb{U}_{\mathfrak m})$. Par contre, nous obervons que tous les sous-groupes de $Iso(\mathbb{U}_{\mathfrak m})$ ne sont pas SIN. De façon précise, nous obtenons:
\begin{theorem}
 Si $G$ est un groupe de poids $\leq \mathfrak m$ tel que toute intersection d'une famille d'ouverts de $G$ de cardinalité $< \mathfrak m$ est un ouvert ($P_{\mathfrak m}$-groupe\cite{vp2}), alors $G$ est isomorphe à un sous-groupe du groupe $Iso(\mathbb{U}_{\mathfrak m})$.
 \end{theorem}
La caractérisation complète des sous-groupes de $Iso(\mathbb{U}_{\mathfrak m})$ reste une question ouverte.\\
Les principaux résultats de cette thèse ont fait l'objet de publications: \cite{brice1} pour la première partie et \cite{brice2} pour la deuxième partie. Ce travail s'achève par un annexe où sont étudiées les notions de base utilisées dans notre travail.

\chapter{Systèmes dynamiques topologiques abstraits et généralités sur la moyennabilité}
Ce chapitre regroupe des éléments de bases de la théorie des systèmes dynamiques topologiques abstraits. Ces notions seront d'une grande importance pour la compréhension de la suite de la thèse. Pour plus de détails sur les systèmes dynamiques abstraits, le lecteur pourra consulter \cite{aus} ou \cite{kpt}. Le présent chapitre regroupe également des généralités sur les groupes moyennables.
\section{Systèmes dynamiques topologiques abstraits}
\subsection{Actions et Représentations}
Rappelons qu'un groupe $G$ opère sur un ensemble $E$ s'il existe une application
\[\begin{array}{lll}
\tau:& G\times E &\longrightarrow E\\
&(s,\,x) &\longmapsto \tau(s,x)=s.x\\
\end{array}\]
telle que  $s.(t.x)=(st).x$ et $e.x=x$ pour tous  $s,\,t\in
G$ et $x\in E$.\\
On dit encore que $\tau$ est une action de $G$ sur $E$.\\
 Si $E$ un espace topologique, $G$ un groupe topologique et l'application $\tau$ est continue ($G\times E$ muni de la topologie produit), on dit que $G$ opère continûment sur $E$.
\begin{definition}
Un système dynamique topologique abstrait\index{système!dynamique topologique} est un triplet $(X,\,G,\,\tau)$ tel que:
\begin{enumerate}
  \item X est un espace topologique,
  \item $G$ un groupe topologique,
  \item $\tau$ une action\index{action} continue de $G$ sur $X$.
\end{enumerate}
\end{definition}
Si $(X,G,\tau)$ est un système dynamique topologique abstrait, nous dirons que $X$ est un flot du groupe $G$\index{flot}.
\begin{remark}
\begin{enumerate}
  \item Si l'espace $X$ est compact, alors une action continue de $G$ sur $X$ peut être identifiée à un homomorphisme de groupe topologique $G\longrightarrow Homeo_{c}(X)$ où $Homeo_{c}(X)$ désigne le groupe des homéomorphismes de $X$ sur lui-même muni de la topologie compact-ouvert. De façon précise, à tout $g\in G$, on associe l'application
      \[ \xymatrix@C=2cm@R=0.3em{
 \Theta:  G  \ar[r] &  Homeo(X) &  \\
 \,\,\,\,\,\,\,\,\,\,g \ar@{|->}[r] & \Theta(g) : X \ar[r] &  X \\
& \,\,\,\,\,\,\,\,\,\,\,\,\,\,\,\,\,\,\,\,x \ar@{|->}[r] &  \Theta(g)(x)=gx}
\]

Cette application a bien sûr un sens même si $X$ n'est pas compact. Cependant, on montre que l'homomorphisme $G\longrightarrow Homeo_{c}(X)$ est continue seulement si $X$ est compact. Réciproquement, tout homomorphisme continu d'un groupe topologique $G$ vers le groupe des homéomorphismes $Homeo_{c}(X)$ d'un espace compact $X$ sur lui-même détermine de façon unique une action continue de $G$ sur $X$(\cite{vp1}, Proposition $2.1.4$).
  \item L'espace topologique $X$ peut avoir une structure supplémentaire:
  \begin{enumerate}
    \item Si $X=(X,d)$ est un espace métrique, alors l'action du groupe $G$ sur $X$ est dite par isométries\index{action!par isométries} si pour tout $g\in G$, l'application
    \[\begin{array}{lll}
& X &\longrightarrow X\\
&x &\longmapsto gx\\
\end{array}\]
est une isométrie de $X$ sur lui-même. Dans ce cas, la continuité de l'action est équivalente à la continuité de l'homomorphisme de groupe $G\longrightarrow Iso(X)$ où $Iso(X)$ désigne le groupe de toutes les isométries de $X$ sur lui-même muni de la topologie de la convergence simple(i.e celle induite par la topologie produit $X^{X}$).
    \item Si $X=E$ est un espace de Banach, alors l'action $\tau$ est une représentation\index{representation} de $G$ dans $E$ si pour tout $g\in G$, l'application
     \[\begin{array}{lll}
& E &\longrightarrow E\\
&x &\longmapsto gx\\
\end{array}\]
     est un opérateur linéaire.
    \item Si $X=\mathcal{H}$ est un espace de Hilbert, alors la représentation $\tau$ de $G$ dans $\mathcal{H}$ est dite unitaire\index{representation!unitaire} si pour tout $g\in G$, l'application

         \[\begin{array}{lll}
& \mathcal{H} &\longrightarrow \mathcal{H}\\
&x &\longmapsto gx\\
\end{array}\]
est un opérateur unitaire: $(\tau_{g}x,y)=(x,\tau_{g^{-1}}y)$ pour tout $x,y\in \mathcal{H}$. C'est un cas particulier d'une représentation par isométries.
  \end{enumerate}
\end{enumerate}
\end{remark}

\subsection{Flots et compactifié équivariant}
\begin{definition}
Soit $X$ un flot compact d'un groupe $G$ et soit $x\in X$. On appelle orbite\index{orbite} de $x$ l'ensemble $G.x=\{g.x:\,\,g\in G\}$. On
notera $\overline{G.x}$ son adhérence dans $X$.
\end{definition}
\begin{remark}
\begin{enumerate}
  \item $\overline{G.x}$ est un sous-espace compact et $G$-invariant de $X$.\\
  En effet, si on note $\pi:G\times X\longrightarrow X$ l'action continue de $G$ sur $X$, alors pour tout $g\in G$, l'application $\pi_{g}:X\longrightarrow X$ est également continue et $\pi_{g}(G.x)\subset G.x$ pour tout $x\in X$ et pour tout $g\in G$. Ainsi, on a:
    $\pi_{g}(\overline{G.x})\subset \overline{\pi_{g}(G.x)}\subset \overline{G.x}$

 \item En général, Si $Y\subseteq X$ est une partie $G$-invariante compacte non vide de $X$, on peut définir une action de $G$ sur $Y$ en
 considérant la restriction de l'action de $G$ sur $X$.
\end{enumerate}
\end{remark}
\begin{definition}
Soient $X$ et $Y$ deux flots du même groupe topologique $G$. Un morphisme de $X$ dans $Y$ (ou une application équivariante\index{application!équivariante}) est une application continue $\pi:X\longrightarrow Y$ telle que $\pi(g.x)=g.\pi(x)$ pour tout $x\in X$ et pour tout $g\in G$.
\end{definition}

\begin{definition}
Un flot $X$ d'un groupe topologique $G$ est dit minimal\index{$G$-espace!minimal} s'il ne contient pas de sous-flot propre de $G$.
\end{definition}
\begin{lemma}
Un flot $X$ de $G$ est minimal si et seulement si pour tout $x\in X$, l'orbite de $x$ est dense dans $X$.
\end{lemma}
\begin{prof}
Soit $X$ un flot minimal et soit $x\in X$. $\overline{G.x}$ est fermé dans $X$, donc compact. Comme $\overline{G.x}$ est de plus invariant, on a: $\overline{G.x}=X$ car $X$ est minimal.\\
Supposons maintenant que $X$ est non minimal et soit $N$ un sous-ensemble fermé de $X$ tel que $\emptyset \neq N \subsetneq X$.
Si $x\in N$, alors $\overline{G.x}\subset N$ alors $\overline{G.x}\neq X$.
\end{prof}
\begin{theorem}
Tout flot $X$ d'un groupe topologique $G$ contient un sous-flot minimal $Y\subseteq X$.
\end{theorem}
\begin{prof}
Posons $\mathcal{M}=\{\emptyset \neq N \subseteq X\,\,\,\, \text{tel que}\,\, N\,\,\, \text{est fermé et invariant}\}$.\\
Puisque $X\in \mathcal{M}$, alors $\mathcal{M}\neq \emptyset$.\\
$\mathcal{M}$ est partiellement ordonné par l'inclusion.\\
Soit $\{M_{\alpha}\}$ une chaîne de sous-ensemble de $\mathcal{M}$. Il est clair que la famille $\{M_{\alpha}\}$ possède la propriété d'intersection finie. Comme $X$ est compact, $M^{\star}=\cap M_{\alpha}\neq \emptyset$.\\
Ainsi, $M^{\star}\in \mathcal{M}$ et par le Lemme de Zorn $\mathcal{M}$ contient un élément minimal qui est un sous-flot minimal de $X$.
\end{prof}\\
Nous démontrerons plus loin le théorème:

\begin{theorem}\label{flotmi}
Si $G$ est un groupe topologique, alors il existe un flot minimal $M(G)$ de $G$ vérifiant la propriété suivante:\\
Pour tout flot minimal $X$ de $G$, il existe un morphisme $\pi: M(G)\longrightarrow X$. De plus, $M(G)$ est uniquement déterminé à isomorphisme près par cette propriété.
\end{theorem}
\begin{definition}
Le flot minimal $M(G)$ du théorème précédent est appelé flot minimal universel de $G$\index{$G$-espace!minimal universel}.
\end{definition}
Introduisons à présent une autre notion importante de la théorie des systèmes dynamiques topologiques abstraits.
\begin{definition}
Soit $G$ un groupe topologique.
\begin{enumerate}
\item On appelle compactifié équivariant \index{compactifié équivariant} de $G$, tout couple $(X,x_{0})$ où $X$ est un flot de $G$ et $x_{0}\in X$ possède une orbite dense dans $X$.
\item Soient $(X,x_{0})$ et $(Y,y_{0})$ deux compactifiés équivariants de $G$. Un morphisme de $(X,x_{0})$ dans $(Y,y_{0})$ est un morphisme de flots $\pi:X\longrightarrow Y$ vérifiant $\pi(x_{0})=y_{0}$.
 \end{enumerate}
\end{definition}
Le théorème suivant est un résultat important de la théorie des systèmes dynamiques topologiques abstraits.
\begin{theorem}\label{grea}
Si $G$ est un groupe topologique; alors il existe un compactifié équivariant $(Y,y_{0})$ vérifiant la propriété suivante:\\
Pour tout compactifié équivariant $(X,x_{0})$, il existe un morphisme de compactifié équivariant $\pi: Y\longrightarrow X$.
De plus, $(Y,y_{0})$ est uniquement déterminé à isomorphisme près par cette propriété.
\end{theorem}
\begin{definition}
Le compactifié équivariant $(Y,y_{0})$ du théorème précédent est appelé compactifié équivariant de Samuel ou compactifié équivariant universel.
\end{definition}
Dans le prochain paragraphe, nous allons construire cet important objet afin de démontrer le théorème \ref{grea}.

\subsection{Compactifié équivariant de Samuel \index{compactifié équivariant!de Samuel}}

Soit $G$ un groupe topologique. Notons $E=RUCB(G)$ l'ensemble des fonctions $x:G\longrightarrow \mathbb{C}$ bornées et uniformément continues à droite\index{fonction!uniformément continue à droite} (pour tout $\varepsilon > 0$, il existe un voisinage $V$ de $e$ tel que $gh^{-1}\in V\Rightarrow |x(g)-x(h)|< \varepsilon$ pour tout $g,h\in G$).\\
En munissant $E$ de l'addition et de la multiplication point par point, prenant l'involution comme la conjugaison, et avec la norme
 $\|x\|_{\infty}=\sup \{|x(g)|:\,\,g\in G\}$, $E$ est une $C^{\star}$-algèbre commutative et unifère.\\
 Notons $S(G)$ l'ensemble des caractères non nuls de la $C^{\star}$-algèbre commutative et unifère $E$ muni de la topologie faible-étoile que nous appelerons ici topologie de Gelfand.\\
 D'après le Théorème de Gelfand(Theorème \ref{gel}), $S(G)$ est compact et $E$ s'identifie via la tansformation de Gelfand à l'algèbre $C(S(G))$ des fonctions continues sur $S(G)$.\\De façon précise, $S(G)$ est l'ensemble des morphismes d'algèbre non identiquement nuls et continues $\varphi: E\longrightarrow \mathbb{C}$ muni de la topologie engendrée par la famille d'applications\index{topologie!vague}

 \[\begin{array}{lll}
\widehat{x}:& S(G) &\longrightarrow \mathbb{C}\\
&\varphi &\longmapsto \widehat{x}(\varphi)=\varphi(x)\\
\end{array}\]
 Nous allons identifier $x$ à $\widehat{x}$ en cas de besoin.
 \begin{lemma}\label{lemgrea}
 Si $G$ est un groupe topologique, alors $G$ opère continûment sur $S(G)$.
 \end{lemma}
 \begin{prof}
 Soit $G$ un groupe topologique. Notons comme précédement $E=RUCB(G)$. L'application
 \[ \xymatrix@C=2cm@R=0.3em{
 \eta:G\times E \ar[r] &  E &  \\
 (g,x) \ar@{|->}[r] & g.x=(^{g}x) : G \ar[r] &  \mathbb{C} \\
& \,\,\,\,\,\,\,\,\,\,\,\,\,\,\,\,\,\,\,\,\,\,\,\,\,\,\,\,\,\,h \ar@{|->}[r] &  ^{g}x(h)=x(g^{-1}h) }
\]
 est bien définie.\\
  En effet, si $g\in G$ et $x\in E$, alors $\|^{g}x\|_{\infty}=\|x\|_{\infty}$. \\
  Soit $\varepsilon> 0$, il existe un voisinage $V$ de $e$ dans $G$ tel que $bc^{-1}\in V\Rightarrow |x(b)-x(c)|\leq \varepsilon$. Posons $W=gVg^{-1}$. Si $bc^{-1}\in W$ alors $(g^{-1}b)(g^{-1}c)^{-1}\in V$. Ainsi, $$|^{g}x(b)-(^{g}x)(c))|=|x(g^{-1}b)-x(g^{-1}c)|< \varepsilon.$$ Donc $^{g}x$ est uniformément continue à droite. Il est facile de montrer que application $\eta$ définit une action de $G$ sur $E$.\\
  Pour tout $g\in G$, l'application
  \[\begin{array}{lll}
& E &\longrightarrow E\\
&x &\longmapsto g.x\\
\end{array}\]
 laisse invariante les fonctions constantes.  Ainsi nous pouvons définir l'application
  \[ \xymatrix@C=2cm@R=0.3em{
 \xi:G\times S(G) \ar[r] &  S(G) &  \\
 (g,\varphi) \ar@{|->}[r] & g.\varphi=(^{g}\varphi) : E \ar[r] &  \mathbb{C} \\
& \,\,\,\,\,\,\,\,\,\,\,\,\,\,\,\,\,\,\,\,\,\,\,\,\,\,\,\,\,x \ar@{|->}[r] &  ^{g}\varphi(x)=\varphi(g^{-1}.x) }
\]
On a: $^{g}\varphi\in S(G)$ pour tout $g\in G$ et $\varphi\in S(G)$. Ainsi l'application $\xi$ est bien définie et définie clairement une action de $G$ sur $S(G)$.\\
   Montrons que cette action est continue.\\
    Soit $g_{0}\in G,\,\varphi_{0}\in S(G)$ et soit $W$ un voisinage de $g_{0}.\varphi_{0}$ dans $S(G)$. Puiqu'on peut identifier $E$ à $C(S(G))$ et $S(G)$ est complètement régulier, il existe $\widehat{x}\in C(S(G))$ que nous allons indentifier à $x\in E$ par le théorème de Gelfand tel que $(g_{0}.\varphi_{0})(x)=1$ et $\varphi(x)=0$ pour tout $\varphi \in S(G)\setminus W$. Posons $$W_{0}=\{\varphi\in S(G):\,\,\,\varphi(g_{0}^{-1}.x)> \frac{1}{2}\}.$$ Observons que $\varphi_{0}(g_{0}^{-1}.x)=1$, donc $W_{0}$ est un ouvert de $S(G)$ contenant $\varphi_{0}$. Maintenant, posons $$V_{0}=\{g\in G:\,\,\|g^{-1}.x-g_{0}^{-1}.x\|_{\infty}\leq \frac{1}{2}\}.$$ Il existe un voisinage $V$ de $e$ dans $G$ tel que $$bc^{-1}\in V\Rightarrow |x(b)-x(c)|\leq \frac{1}{2}.$$ Si $g\in Vg_{0}$, alors $gb(g_{0}b)^{-1}=gg_{0}^{-1}\in V$. Ainsi $$|(g^{-1}.x)(b)-(g_{0}^{-1}.x)(b)|=|x(gb)-x(g_{0}b)|\leq \frac{1}{2}$$ pour tout $b\in G$ et $g\in Vg_{0}$. Donc $V_{0}\supseteq Vg_{0}$ est un voisinage de $g_{0}$. \\
   Si $g\in V_{0}$ et $\varphi \in W_{0}$, alors

   $|\varphi(g^{-1}x)-\varphi(g_{0}^{-1}x)|\leq\|g^{-1}x-g_{0}^{-1}x\|_{\infty}\leq \frac{1}{2}$. Car $\varphi$ est un caractère.
    Donc

   $(g.\varphi)(x)=\varphi(g^{-1}.x)\geq \varphi(g_{0}^{-1}.x)-\frac{1}{2}>0$. Ainsi $g.\varphi\in W$.
 \end{prof}

\begin{lemma}
 Si $G$ est un groupe topologique, alors $G$ est homéomorphe à un sous espace dense de $S(G)$.
\end{lemma}
\begin{prof}
Si $g\in G$, alors on peut définir un élément $\varphi_{g}\in S(G)$ par $\varphi_{g}(x)=x(g)$ pour tout $x\in E=RUCB(G)$.
 Ceci nous permet de définir une application
 \[\begin{array}{lll}
\Psi:& G &\longrightarrow S(G)\\
&g &\longmapsto \varphi_{g}\\
\end{array}\]
Puisque la topologie de $G$ coïncide avec la topologie induite par la structure uniforme droite sur $G$, l'application $\Psi$ est continue.\\
 En effet, pour tout $x\in RUCB(G)$ et pour tout $g\in G$, on a:
 $$\widehat{x}\circ \Psi(g)=\widehat{x}(\varphi_{g})=\varphi_{g}(x)=x(g).$$
Pour montrer que $\{\varphi_{g}:\,\,\,g\in G\}$ est dense dans $S(G)$, nous allons utiliser un corollaire du théorème de Hahn-Banach(Voir le lemme \ref{hahnban}).\\
Supposons le contraire. D'après le lemme \ref{hahnban}, il existe $h\in C(S(G))$ non nul tel que $h(\varphi_{g})=0$ pour tout $g\in G$.\\
 D'après le théorème de Gelfand, il existe $x\in RUCB(G)$ tel que $z(x)=h(z)$ pour tout $z\in S(G)$. Dans ce cas, $x$ ne peut pas être la fonction nulle.\\Mais $x(g)=\varphi_{g}(x)=h(\varphi_{g})=0$ pour tout $g\in G$. Ce qui est absurde.\\
Soit $g_{0},h_{0}\in G$ tel que $g_{0}\neq h_{0}$. Il existe une pseudométrique continue bornée et invariante à droite $d_{r}$ sur $G$ telle que $d_{r}(g_{0},h_{0})\neq 0$ (voir \cite{hewit}, $8.2$). En posant $x(g)=d_{r}(g,h_{0})$, nous avons :
$$\varphi_{g_{0}}(x)=x(g_{0})=d_{r}(g_{0},h_{0})\neq 0 = d_{r}(h_{0},h_{0})=x(h_{0})=\varphi_{h_{0}}(x).$$ Ainsi, $\varphi_{g_{0}}\neq \varphi_{h_{0}}$. Ce qui montre l'injectivité.
\end{prof}

\begin{remark}
\begin{enumerate}
  \item On peut également voir $S(G)$ comme le compactifié de Samuel de $G$ par rapport à la structure uniforme droite sur $G$ (\cite{brook}). Dans la terminologie anglo-saxone, on l'appelle "greatest ambit". Nous l'appelerons ici tout simplement compactifié de Samuel équivariant.
  \item Si $G$ est un groupe dénombrable et discret, alors le compactifié de Samuel équivariant $S(G)$ coïncide avec le compactifié de Stone-\v Cech $\beta G$.\cite{eng}
\end{enumerate}
\end{remark}
Revenons maintenant à la preuve du théorème \ref{grea}\\
\begin{prof}(du théorème \ref{grea})
Prenons $Y=S(G)$ et $y_{0}=e$. Puisque l'orbite de $e$ dans $S(G)$ est $G$ qui est dense, $(S(G),e)$ est clairement un compactifié équivariant de $G$.\\
Par le lemme \ref{lemgrea}, $G$ opère sur $E$ par l'action: $g.x(h)=x(g^{-1}h)$, puis canoniquement continûment sur $S(G)$ par l'action $g.\varphi(x)=\varphi(g^{-1}.x)$.\\
Considérons maintenant un compactifié équivariant arbitraire $(X,x_{0})$ de $G$.\\
Soit $f\in C(X)$. Définissons l'application
 \[\begin{array}{lll}
f^{\ast}:& G &\longrightarrow \mathbb{C}\\
&g &\longmapsto f^{\ast}(g)=f(gx_0)\\
\end{array}\]
Montrons que $f^{\ast} \in E$.\\
 Pour tout $x\in X$, il existe un voisinage $U_{x}$ de $x$ dans $X$ tel que $|f(y)-f(x)|<\varepsilon$ pour tout $y\in U_{x}$ car $f\in C(X)$.\\
Puisque l'action de $G$ sur $X$ est continue, et $ex=x$, il existe un voisinage $V_{x}$ de $x$ dans $X$ et un voisinage symétrique $O_{x}$ de $e$ dans $G$ tel que $O_{x}V_{x}\subset U_{x}$. Puisque $e\in O_{x}$, on a: $V_{x}\subset U_{x}$.\\
Comme $X$ est compact, il existe $F\subset X$ fini tel que $X\subset \underset{x\in F}{\overset{}{\bigcup}}V_{x}$.\\
Posons $O=\underset{x\in F}{\overset{}{\bigcap}}O_{x}$. Il est clair que $O$ est un voisinage symétrique de $e$ dans $G$.\\
Si $g\in O$ et $y\in X$, alors il existe $x\in F$ tel que $y\in V_{x}$. Ainsi, $g.y\in OV_{x}\subset O_{x}V_{x}\subset U_{x}.$\\
 On peut donc conclure que pour tout $\varepsilon> 0$, il existe un voisinage $V$ de $e$ dans $G$ tel que $g\in V\,\Longrightarrow\,\,|f(g.x)-f(x)|<\varepsilon$ pour tout $x\in X$.\\
Si $gh^{-1}\in V$, alors $$|f^{\ast}(g)-f^{\ast}(h)|=|f(g.x_{0})-f(h.x_{0})|=|f(gh^{-1}(h.x_{0}))-f(h.x_{0})|<\varepsilon.$$
Donc $f^{\ast}\in E$. Puisque $f^{\ast}$ est clairement bornée.\\
En identifiant comme souvant $E=RUCB(G)$ avec $C(S(G))$, nous avons donc un monomorphisme de $C^{\star}$-algèbres
\[\begin{array}{lll}
\pi:& C(X) &\longrightarrow C(S(G))\\
&f &\longmapsto f^{\ast}\\
\end{array}\]
Il est connu (voir par exemple \cite{const}, $2.4.3.6$) que tout monomorphisme unitaire de $C^{\star}$-algèbres $\pi:C(K)\longrightarrow C(L)$ où $K$ et $L$ sont des espaces compacts non vides est de la forme $\pi(f)=f\circ \Pi$ pour une surjection unique $\Pi:L\longrightarrow K$. De plus, si $K$ et $L$ sont des $G$-espaces, si nous faisons agir $G$ sur $C(K)$ et $C(L)$ par l'action $g.f(x)=f(g^{-1}.x)$ et si $\pi$ est équivariante, alors $\Pi$ est équivariante. En appliquant ceci à l'application \[\begin{array}{lll}
\pi:& C(X) &\longrightarrow C(S(G))\\
&f &\longmapsto f^{\ast}\\
\end{array}\] précédente, il existe un unique homomorphisme de $G$-flot $\Pi:S(G)\longrightarrow X$ avec $f^{\ast}=f\circ \Pi$.
Il nous reste juste à montrer que $\Pi(e)=x_{0}$.\\
Pour tout $f\in C(X)$, on a: $f^{\ast}(e)=f(x_{0})=f(\Pi(e))$. Ainsi, on a: $\Pi(e)=x_{0}$.
\end{prof}
\begin{definition}
\begin{enumerate}
  \item Un semi-groupe\index{semi-groupe} est la donnée d'un ensemble non vide $X$ muni d'une opération associative
  \item Un semi-groupe $X$ est dit semi-topologique à gauche \index{semi-groupe! semi-topologique} s'il existe une topologie sur $X$ tel que pour tout $y\in X$, l'application
      \[\begin{array}{lll}
& X &\longrightarrow X\\
&x &\longmapsto xy\\
\end{array}\]
  est continue.
\end{enumerate}
\end{definition}
\begin{theorem}
Pour tout groupe topologique $G$, le compactifié équivariant de Samuel $S(G)$ de $G$ possède une structure de semi-groupe semi-topologique à gauche et la multiplication $S(G)\times S(G)\longrightarrow S(G)$ prolonge l'action  $G\times S(G)\longrightarrow S(G)$.
\end{theorem}
\begin{prof}
Soit $x,y\in S(G)$, par la propriété universelle de $S(G)$, il existe un morphisme de compactifié équivariant $r_{y}:S(G)\longrightarrow S(G)$ tel que $r_{y}(e)=y$.\\
Définisons l'opération $xy=r_{y}(x)$.\\Montrons que l'opération $(x,y)\longmapsto xy$ définie une structure de semi-groupe semi-topologique sur $S(G)$. Soit $y$ fixé, l'application $x\longmapsto xy$ coïncide avec $r_{y}$. Donc est continue.\\
Soit $y,z\in S(G)$. Par la partie unicité du théorème \ref{grea}, les applications $r_{z}r_{y}$ et $r_{yz}$ définies de $S(G)$ dans $S(G)$ coïncident. D'où l'associativité de l'opération. Puisque $ex=r_{x}(e)=x$ et $xe=r_{e}(x)=x$, $e$ est l'unité de $S(G)$.\\
Pour terminer, Soit $g\in G$ et $x\in S(G),\,\,gx$ peut être compris de deux manières:
\begin{itemize}
  \item Comme l'action de $G$ sur $S(G)$
  \item Comme le produit dans $S(G)$
\end{itemize}
En effet, $gx=r_{x}(g)=r_{x}(ge)=gr_{x}(e)=gx$. Autrement dit, l'opération de semi-groupe sur $S(G)$ prolonge l'action continue de $G$ sur $S(G)$.
\end{prof}
\begin{definition}
\begin{enumerate}
  \item Une partie $I$ de $S(G)$ est un idéal à gauche si $S(G)I\subset I$.
  \item Un élément $x$ d'un semi-groupe est dit idempotent si $x^{2}=x$
\end{enumerate}
\end{definition}
\begin{remark}
\begin{enumerate}
  \item Tout sous-espace $G$-invariant fermé de $S(G)$ est un idéal à gauche.
En effet, soit $X$ un sous-espace $G$-invariant fermé de $S(G)$. Pour $a\in X$, on a: $r_{a}(G)=\{ga:\,\,\,g\in G\}\subseteq X$. Puisque $r_{a}$ est continue, $r_{a}(S(G))=r_{a}(\overline{G})\subseteq \overline{r_{a}(G)}\subseteq X$. Donc $X$ est un idéal à gauche
  \item Dans la suite,  $M(G)$ sera un
  sous-flot minimal de $S(G)$ un sous-espace compact et $G$-invariant de $S(G)$). \end{enumerate}

\end{remark}
Démontrons le théorème \ref{flotmi}, en établissant le corollaire suivant du théorème \ref{grea}
\begin{corollary}\label{coroflot}
Pour tout $G$-espace minimal $X$, il existe un homomorphisme $\pi$ de $M(G)$ dans $X$.
\end{corollary}
\begin{prof}

Soit $X$ un flot minimal de $G$. Fixons $x_{0}\in X$. Alors, $(X,x_{0})$ est un compactifié équivariant de $G$. D'après le théorème \ref{grea}, il existe un morphisme de compactifiés équivariants $\pi:(S(G),e)\longrightarrow (X,x_{0})$. Il est clair que la restriction de $\pi$ à $M(G)$ reste encore un homomorphisme.
\end{prof}\\
Montrons pour terminer que $M(G)$ est unique à isomorphisme près.\\De façon précise, montrons la dernière partie du Théorème \ref{grea} i.e le théorème suivant:
\begin{theorem}\label{unicité}
Tout flot minimal compact de $G$ vérifiant la propriété du corollaire \ref{coroflot} est isomorphe à $M(G)$.
\end{theorem}
Avant de procéder à la preuve de cet important théorème, faisons les rappels suivants:
Pour $a\in S(G)$, nous noterons $r_{a}$ l'application $x\longmapsto xa$ de $S(G)$ dans lui-même
\begin{lemma}
Si $f:S(G)\longrightarrow S(G)$ est un $G$-morphisme et $a=f(e)$, alors $f=r_{a}$.
\end{lemma}
\begin{prof}
Pour tout $x\in G$, nous avons:$f(x)=f(xe)=xf(e)=xa=r_{a}(x)$
\end{prof}\\
Nous aurons besoin du théorème suivant:
\begin{theorem}(Ellis \cite{vp2})\label{thellis}
Tout semi-groupe semi-topologique compact non vide $K$ contient un élément idempotent.
\end{theorem}
\begin{lemma}\label{iso}
Tout $G$-morphisme $f:M(G)\longrightarrow M(G)$ est un $G$-isomorphisme
\end{lemma}
\begin{prof}
Comme $M(G)$ est un sous-semi-groupe fermé, d'après le théorème de Ellis(théorème \ref{thellis}), $M(G)$ contient un élément idempotent $p$.
Comme $M(G)$ est minimal, on a: $M(G)p=M(G)$. Ainsi: $r_{p}(x)=r_{p}(yp)=yp^{2}=yp=x$ pour tout $x\in M(G)$.\\
Si $f:M(G)\longrightarrow M(G)$ est un morphisme de $G$-espaces, alors $f\circ r_{p}$ est un morphisme de $S(G)$ dans $M(G)\subseteq S(G)$.
\\Puisque $f\circ r_{p}(e)=f(r_{p}(e))=f(p)$, on a $f\circ r_{p}=r_{b}$ pour $b=f(p)$.\\
Puisque la restriction de $f\circ r_{p}$ à $M(G)$ coïncide avec $f$, on a: $f(x)=xb$ pour tout $x\in M(G)$.\\
Une fois encore $M(G)b=M(G)$ donc $p=cb$ avec $c\in M(G)$.\\ Ainsi, le morphisme $g=r_{c}:M(G)\longrightarrow M(G)$ est l'inverse à droite de $f$.\\En effet, $fg(x)=xcb=xp=x$.\\
Nous venons de montrer que dans le semi-groupe $S$ des $G$-homomorphismes, chaque élément possède un inverse à droite. Donc $S$ est un groupe.
\end{prof}\\
Nous pouvons à présent démontrer le théorème \ref{unicité}\\
\begin{prof}(du théorème \ref{unicité})\\
Soit $M^{\prime}$ un autre flot compact universel de $G$, alors il existe des $G$-morphismes $f:M(G)\longrightarrow M^{\prime}$ et
$g:M^{\prime}\longrightarrow M(G)$. Puisque $M^{\prime}$ est minimal, $f$ est surjective. Par le lemme \ref{iso}, l'application $gf:M(G)\longrightarrow M(G)$ est bijective. Ainsi $f$ est injective
\end{prof}\\
D'après ce qui précède, on peut associer à tout groupe topologique $G$ un flot compact minimal $M(G)$ unique à un isomorphisme près.\\
Il est possible de caractériser cet espace pour certains groupes particuliers comme le montre le résultat suivant:
\begin{theorem}(Pestov \cite{vp2}, Théorème $6.6$)\label{exapes}
Si $\mathbb{S}^{1}$ désigne le cercle unité de $\mathbb{C}$ et $Homeo(\mathbb{S}^{1})$ le groupe de tous les homéomorphismes de $\mathbb{S}^{1}$ dans $\mathbb{S}^{1}$ muni de la topologie compact-ouvert, alors $\mathbb{S}^{1}$ est le flot compact universel du groupe $G=Homeo(\mathbb{S}^{1})$.
\end{theorem}
Il est naturel de se poser la question de savoir si tout groupe topologique admet une action continue sur un espace compact? Cette question fondamentale a été répondu par Teleman \cite{tel} de la manière suivante:
\begin{theorem}(Teleman \cite{tel})\label{tel}\index{théorème!de Teleman}
Tout groupe topologique $G$ opère effectivement:
\begin{enumerate}
  \item sur un espace de Banach par isométries
  \item sur un espace compact.
\end{enumerate}
\end{theorem}
\section{Généralités sur la moyennabilité}
Dans ce cette section, on rappelle la définition de la moyennabilité et les propriétés de base sur la moyennabilité.
\subsection{Moyenne invariante}
Soit $G$ un groupe topologique séparé. Si $f$ est une fonction sur $G$ à valeurs complexes, et si $s\in G$, posons
\[\begin{array}{lll}
^{s}f:& G &\longrightarrow \mathbb{C}\\
&x &\longmapsto\, ^{s}f(x)=f(s^{-1}x)\\
\end{array}\]
et notons comme précédement, $RUCB(G)$ l'espace de Banach de toutes les fonctions sur $G$ à valeurs complexes bornées et uniformément continue à droite.
\begin{definition}\cite{ey}
Soit $\mathcal{E}$ un sous-espace de Banach de $RUCB(G)$ tel que:
\begin{enumerate}
  \item $1\in \mathcal{E}$
  \item $f\in \mathcal{E}$ implique $\overline{f}\in \mathcal{E}$
\end{enumerate}
une moyenne\index{moyenne} sur $\mathcal{E}$ est une forme linéaire $m$ sur $\mathcal{E}$ telle que:
\begin{enumerate}
  \item $m(1)=1$
  \item Pour toute $f\in \mathcal{E}$, on a: $m(\overline{f})=\overline{m(f)}$
  \item $m(f)\geq 0$ pour toute $f\geq 0$
\end{enumerate}
\end{definition}
\begin{remark}
Une moyenne $m$ sur $\mathcal{E}$ est automatiquement continue.\\
En effet, $-\|f\|_{\infty}1\leq f\leq \|f\|_{\infty}1$. Ainsi, $|m(f)|\leq \|f\|_{\infty}$.
\end{remark}
Si nous supposons de plus que $\mathcal{E}$ est stable par translation, c'est-à-dire que $f\in \mathcal{E}$ et $s\in G$ impliquent $^{s}f\in \mathcal{E}$, alors une moyenne $m$ sur $\mathcal{E}$ est dite invariante si pour tous $f\in \mathcal{E}$ et $s\in G$, on a: $$m(^{s}f)=m(f)$$

\begin{definition}
Un groupe topologique $G$ est dit \index{groupe!moyennable}moyennable s'il existe une moyenne invariante sur $RUCB(G)$.
\end{definition}
Si $G$ est un groupe compact, alors $RUCB(G)=C(G)$ et une moyenne invariante $m$ sur $C(G)$ est une mesure borélienne invariante et régulière sur $G$ avec $m(1_{G})=1$. Ainsi, la mesure normalisée de Haar est l'unique moyenne invariante sur $C(G)$. En particulier, tout groupe compact est moyennable.

\subsection{Moyennabilité et propriété du point fixe}
\subsubsection{Mesures boréliennes}
Soit $X$ un espace topologique. La tribu borélienne de $X$ que l'on notera $\mathcal{B}(X)$ est la plus petite tribu contenant tous les ouverts de $X$. Une mesure définie sur $\mathcal{B}(X)$ et prenant des valeurs finies sur les compacts est dite borélienne. Une mesure borélienne $\mu$ est dite régulière si elle est intérieurement régulière et extérieurement régulière, c'est-à-dire si pour tout borélien $B\in \mathcal{B}(X)$, on a: $$\mu(B)=\inf\{\mu(V):\,\,\,B\subset V,\,\,V\,\,\text{est ouvert}\}$$
 et
$$\mu(B)=\sup\{\mu(K):\,\,\,K\subset B,\,\,K\,\,\text{est compact}\}$$

Supposons $X$ compact. On note $\mathbb{P}(X)$ l'ensemble des mesures de probabilité boréliennes régulières sur $X$.\\
 Soit $x\in X$. La masse de Dirac en $x$ est la mesure borélienne $\delta_x$ définie pour tout $B\in \mathcal{B}(X)$ par
$$\delta_{x}(B)=\left\{
    \begin{array}{ccccc}
      1 &\text{si}\,\, x\in B&  &  \\
       0 & \text{sinon} & &  \\
    \end{array}\right.$$
    Notons que $\delta_x\in \mathbb{P}(X)$.\\
Soit $C(X)$ l'espace vectoriel réel des fonctions continues de $X$ dans $\mathbb{R}$ que l'on muni de la norme
$\|f\|=\underset{x\in X}{\overset{}{\sup}}|f(x)|$.\\
 Si $\mu \in \mathbb{P}(X)$, alors l'application $L_{\mu}:C(X)\longrightarrow \mathbb{R}$ définie par $L_{\mu}(f)=\int_{X}fd\mu$ pour tout $f\in C(X)$ est une forme linéaire continue positive de norme $\|L_{\mu}\|=1$. Réciproquement, il résulte du théorème de representation de Riesz (voir \cite{ru}, théorème 2.14) que si $L:C(X)\longrightarrow \mathbb{R}$, est une forme linéaire continue positive de norme $1$, alors il existe une unique mesure de probabilité $\mu \in \mathbb{P}(X)$ vérifiant $L=L_{\mu}$. On peut donc identifier $\mathbb{P}(X)$ à l'espace $\{L\in C(X)^{\star}:\,\,L\geq 0,\,\,\,L(1_{X})=1\}$. L'espace $\mathbb{P}(X)$ s'identifie donc à un sous-ensemble convexe de la boule unité $B=\{L\in C(X)^{\star}:\,\|L\|\leq 1\}$ de $C(X)^{\star}$. Où $C(X)^{\star}$ désigne le dual topologique de $C(X)$. \\
En effet, si $\mu\in \mathbb{P}(X)$, alors: $$-\mu(|f|)=\mu(-|f|)\leq \mu(f)\leq \mu(|f|),$$ ainsi $|\mu(f)|\leq \mu(|f|)$. Avec $|f|(x)=|f(x)|$ pour tout $x\in X$.\\
De plus, nous avons:
$$\mu(|f|)\leq \mu(\|f\|1_{X})=\|f\|\mu(1_{X})=\|f\|,$$ et $|\mu(f)|\leq \|f\|$. Donc $\mu$ est continue et $\|\mu\|\leq 1$.\\
Rappelons que la topologie vague sur $C(X)^{\star}$ est la topologie la moins fine rendant continue toutes les applications $\Psi_f:C(X)^{\star}\longrightarrow\mathbb{R},\,\,f\in C(X)$ où $\Psi_f(L)=L(f)$ pour tout $L \in C(X)^{\star}$. Muni de cette topologie, l'espace $C(X)^{\star}$ devient un espace vectoriel topologique séparé localement convexe dans lequel $B$ est compact. On munit $\mathbb{P}(X)\subset C(X)^{\star}$ de la topologie induite par la topologie vague de $C(X)^{\star}$. Alors $\mathbb{P}(X)$ est compact puisque c'est un sous ensemble fermé du compact $B$. Notons $\mathbb{P}_{0}(X)$ l'enveloppe convexe des masses de Dirac. En application du théorème de Krein-Milman (théorème \ref{krein}), $\mathbb{P}_{0}(X)$ est un sous-espace convexe dense de $\mathbb{P}(X)$.\\
Supposons que $X$ soit muni d'une action continue d'un groupe topologique $G$. On dit qu'une mesure $\mu \in \mathbb{P}(X)$ est $G$-invariante si $\mu(gB)=\mu(B)$ pour tout borélien $B$ et pour tout $g\in G$.

\subsubsection{Action sur $\mathbb{P}(X)$}
\begin{definition}
Soit $X$ un sous-ensemble convexe d'un espace localement convexe.
  \begin{enumerate}
    \item Une application $\alpha:X\longrightarrow X$ est dite \index{affine} affine si: $$\alpha(tx+(1-t)y)=t\alpha(x)+(1-t)\alpha(y)$$ pour tous $x,y\in X,\,\,\,0\leq t\leq 1$.
    \item Soit $G$ un groupe topologique. Une action continue $\tau:G\times X\longrightarrow X$ est dite affine\index{action!affine} si pour tout $g\in G$ l'application orbite
        \[\begin{array}{lll}
\tau_{g}:& X &\longrightarrow X\\
&x &\longmapsto \tau_{g}(x)=gx\\
\end{array}\] est affine.
  \end{enumerate}
  \end{definition}
Soient $X$ et $Y$ deux espaces topologiques, $\alpha: X\longrightarrow Y$ une application continue. Pour toute mesure borélienne $\mu$ sur $X$, on définit la mesure borélienne image sur $Y$ par $\mu^{\star}(B)=\mu(\alpha^{-1}(B))$ pour tout ensemble borélien $B$ de $Y$. \\
Soient $X$ et $Y$ deux espaces compacts, $\alpha : X\longrightarrow Y$ une application continue. $\alpha$ induit une application linéaire continue

\[ \xymatrix@C=2cm@R=0.3em{
 \widetilde{\alpha}:  C(Y)  \ar[r] &  C(X) &  \\
 f \ar@{|->}[r] & \widetilde{\alpha}(f) : X \ar[r] &  \mathbb{C} \\
& \,\,\,\,\,\,\,\,\,\,\,\,\,\,\,\,\,\,x \ar@{|->}[r] &  \widetilde{\alpha}(f)(x)=f(\alpha(x)) }
\]
Remarquons que $\widetilde{\alpha}(f)=f\circ \alpha$ pour $f\in C(X)$.\\
 L'application $\widetilde{\alpha}$ permet de définir une application $\alpha^{\star}:C(X)^{\star}\longrightarrow C(Y)^{\star}$ définie par
  $\alpha^{\star}(\mu)=\mu\circ \widetilde{\alpha}$ pour $\mu \in C(X)^{\star}$. Plus précisement, nous avons: $$\alpha^{\star}(\mu)(f)=\mu(f\circ \alpha)$$ pour tout $f\in C(Y)$ ou de manière équivalente, $$\int_{Y}fd(\alpha^{\star}\mu)=\int_{X}(f\circ \alpha)d\mu$$ pour tout $f\in C(Y)$. Ceci permet de conclure que pour tout sous-ensemble borélien $A$ de $Y$, on a: $\alpha^{\star}(\mu)(A)=\mu(\alpha^{-1}(A))$. Ainsi, $\alpha^{\star}(\mu)$ coïncide avec la mesure image $\mu^{\star}$ définie précédement.\\
  De plus, nous avons les propriétés suivantes:
  \begin{enumerate}
    \item Pour tout $\mu\in C(X)^{\star}$, on a: $Supp(\alpha^{\star}(\mu))=\overline{\alpha[Supp(\mu)]}$. (\cite{vries}, Appendix C.$9$)
    \item Si $\alpha$ est une injection continue, alors $\alpha^{\star}:C(X)^{\star}\longrightarrow C(Y)^{\star}$ est injective. (\cite{vries}, Appendix C.$9$)
    \item Si $Z$ est un espace topologique et $\beta:Y\longrightarrow Z$ est une application continue, alors $(\alpha\circ \beta)^{\star}=\alpha^{\star}\circ \beta^{\star}$. De plus, $(id_X)^{\star}=id_{C(X)^{\star}}$. (\cite{vries}, Appendix C.$9$)
    \item L'injection continue
    \[ \xymatrix@C=2cm@R=0.3em{
 \delta^{X}:  X  \ar[r] &  C(X)^* &  \\
 x \ar@{|->}[r] & \delta^{X}(x)=\delta_x : C(X) \ar[r] &  \mathbb{R} \\
& \,\,\,\,\,\,\,\,\,\,\,\,\,\,\,\,\,\,\,\,\,\,\,\,\,\,\,\,\,f \ar@{|->}[r] &  f(x) }
\]
vérifie $\delta^{(Y)}\circ\alpha=\alpha^{\star}\circ \delta^{(X)}$. Autrement dit, le diagramme suivant

        \[\xymatrix{X\ar[rr]^{\alpha}\ar[dd]_{\delta^{(X)}}&&Y\ar[dd]^{\delta^{(Y)}}\\
\\
C(X)^{\star}\ar[rr]^{\alpha^{\star}}&&C(Y)^{\star}}\] est commutatif:

   \item L'application   $\alpha^{\star}:C(X)^{\star}\longrightarrow C(Y)^{\star}$ est affine et continue par-rapport à la topologie vague respectivement sur $C(X)^{\star}$ et $C(Y)^{\star}$.\\
       En effet, soit $\mu \in C(X)^{\star} $. Rappelons qu'une base voisinage de $\mu$ est formée des ensembles de la forme
  $$V_{\mu}=\{\nu \in C(X)^{\star},\,\,\,\,\,\,|\int_X f_{i} d\nu-\int_X f_{i} d\mu|< \varepsilon,\,\,i=1,2,...,k\}$$
   où $f_{1},f_{2},...f_{k}\in C(X)$ et $\varepsilon> 0$. \\
   Soit $\mu_{0}\in C(X)^{\star}$. Montrons que $\alpha^{\star}$ est continue en $\mu_{0}$.
   Soit $$V_{\alpha^{\star}(\mu_{0})}=\{\nu \in C(Y)^{\star},\,\,\,|\int_Y h_{i} d\nu-\int_Y h_{i} d\alpha^{\star}(\mu_{0})|< \varepsilon,\,\,i=1,2,...,k\}$$ un voisinage de $\alpha^{\star}(\mu_{0})$.
   Cherchons $U_{\mu_{0}}=\{\nu \in C(X)^{\star},\,\,\,\,\,\,|\int_X f_{i} d\nu-\int_X f_{i} d\mu_{0}|< \delta,\,\,i=1,2,...,l\}$ tel que:\\
   $$\alpha^{\star}(U_{\mu_{0}})\subset V_{\alpha^{\star}(\mu_{0})}.$$
   Soit $\lambda\in U_{\mu_{0}}$. Si $\alpha^{\star}(\lambda) \in V_{\alpha^{\star}(\mu_{0})}$, alors
   $$|\int_Y h_{i} d\alpha^{\star}(\lambda)-\int_Y h_{i} d\alpha^{\star}(\mu_{0})|< \varepsilon,$$ pour tout $i=1,2,...,k$.
   Donc $$|\int_X h_{i}\circ \alpha d\lambda-\int_X h_{i}\circ \alpha d\mu_{0}|< \varepsilon$$ pour tout $i=1,2,...,k$.\\
   Posons $f_{i}=h_{i}\circ \alpha$ pour tout $i=1,2,...,k$ et choisissons $\delta<\varepsilon$. On a le résultat.
   \item La restriction de  $\alpha^{\star}:C(X)^{\star}\longrightarrow C(Y)^{\star}$ à $\mathbb{P}(X)$ est continue, affine et envoie $\mathbb{P}(X)$ sur $\mathbb{P}(Y)$.
  \end{enumerate}
On en déduit le lemme suivant:
\begin{lemma}\label{prolon}
Soit $G$ un groupe topologique. Toute action continue de $G$ sur un espace compact $X$, se prolonge en une action affine continue de $G$ sur $\mathbb{P}(X)$.
\end{lemma}
\begin{prof}
Soit $X$ un $G$-espace compact. Notons
\[\begin{array}{lll}
\pi:& G\times X &\longrightarrow X\\
&(g,x) &\longmapsto gx\\
\end{array}\]
l'action continue de $G$ sur $X$. Pour tout $g\in G$, l'application

\[\begin{array}{lll}
\pi_{g}:& X &\longrightarrow X\\
&x &\longmapsto gx\\
\end{array}\]
est un homéomorphisme. On déduit de ce qui précède que pour tout $g\in G$, il existe un homéomorphisme $\pi_{g}^{\star}:\mathbb{P}(X)\longrightarrow \mathbb{P}(X)$. D'après la propriété $3$ précédente, l'application

\[\begin{array}{lll}
\pi^{\star}:& G\times \mathbb{P}(X) &\longrightarrow \mathbb{P}(X)\\
&(g,\mu) &\longmapsto \pi_{g}^{\star}(\mu)\\
\end{array}\]
est une action de $G$ sur $\mathbb{P}(X)$. Comme toutes les applications $\pi_{g}^{\star}:\mathbb{P}(X)\longrightarrow \mathbb{P}(X)$ sont continues, l'action $\pi^{\star}$ de $G_{d}$ dans $\mathbb{P}(X)$ est continue, $G_d$ étant le groupe $G$ muni de la topologie discrète. Nous écrirons dans la suite tout simplement $g\mu$ pour désigner $\pi_{g}^{\star}\mu$.
Par définition de $\pi_{g}^{\star}\mu$, on a:
$$(g\mu)(f)=\mu(f\circ \pi_g)=\int_{X}f(gx)d\mu(x)$$ pour tout $g\in G$ et $f\in C(X)$ et $$(g\mu)(B)=\mu(g^{-1}B)$$ pour tout $g\in G$ et $B$ un sous-ensemble borélien de $X$. L'action précédente de $G$ sur $\mathbb{P}(X)$ est continue même si $G$ est muni de sa topologie initiale.\\
En effet, il nous suffit en vertu des propriétés de définition de la topologie vague sur $\mathbb{P}(X)$ de montrer que pour tout $f\in C(X)$, l'application \[\begin{array}{lll}
\Psi:& G\times \mathbb{P}(X) &\longrightarrow \mathbb{R}\\
&(g,\mu) &\longmapsto \mu(f\circ \pi_g)\\
\end{array}\] est continue. \\
Soit $(g,\mu)\in G\times\mathbb{P}(X),\,\,\,f\in C(X)$ et soit $\varepsilon> 0$. Considérons le voisinage  $$V=\{\nu \in \mathbb{P}(Y),\,\,\,|\int_Y f_{i} d\nu-\int_Y f_{i} d\alpha^{\star}(\mu)|< \frac{\varepsilon}{2}\,\,i=1,2,...,k\}$$ où $f_{1},f_{2},...f_{k}\in C(X)$ de $\mu$ dans $\mathbb{P}(X)$.

Pour tout $x\in X$, il existe un voisinage $U_{x}$ de $x$ dans $X$ tel que $|f(y)-f(x)|<\frac{\varepsilon}{2}$ pour tout $y\in U_{x}$.
Puisque l'action de $G$ sur $X$ est continue, et $ex=x$, il existe un voisinage $V_{x}$ de $x$ dans $X$ et un voisinage symétrique $O_{x}$ de $e$ dans $G$ tel que $O_{x}V_{x}\subset U_{x}$. Puisque $e\in O_{x}$, on a: $V_{x}\subset U_{x}$.\\
Comme $X$ est compact, il existe $F\subset X$ fini tel que $X\subset \underset{x\in F}{\overset{}{\bigcup}}V_{x}$.\\
Posons $O=\underset{x\in F}{\overset{}{\bigcap}}O_{x}$. Il est clair que $O$ est un voisinage symétrique de $e$ dans $G$.\\
Si $g\in O$ et $y\in X$, alors il existe $x\in F$ tel que $y\in V_{x}$. Ainsi, $$g.y\in OV_{x}\subset O_{x}V_{x}\subset U_{x}.$$
Si $(h,\nu)\in O\times V$, on a:

 \[\begin{array}{llllll}
 |\nu(f\circ \pi_{h})-\mu(f\circ \pi_{g})|&\leq
|\nu(f\circ \pi_{h})-\nu(f\circ \pi_{g})|+|\nu(f\circ \pi_{g})-\mu(f\circ \pi_{g})|
\\
&\leq \|f\circ \pi_{h}-f\circ \pi_{g}\|+|\nu(f\circ \pi_{g})-\mu(f\circ \pi_{g})|\\
\\
&\leq \frac{\varepsilon}{2}+\frac{\varepsilon}{2}=\varepsilon
\end{array}\] \\
Ainsi l'application $\Psi$ est continue.
\end{prof}\\
En général, on a le lemme suivant:
\begin{lemma}\label{lemaction}
Soient $X$ et $Y$ deux $G$-espaces compacts et soit $\alpha : X\longrightarrow Y$ une application continue et équivariante, alors $\mathbb{P}(X)$ et $\mathbb{P}(Y)$ sont également des $G$-espaces et l'application $\alpha^{\star}:\mathbb{P}(X)\longrightarrow \mathbb{P}(Y)$ définie précédemment est équivariante.
\end{lemma}

\begin{prof}
Il nous reste seulement à vérifier que $\alpha^{\star}:\mathbb{P}(X)\longrightarrow \mathbb{P}(Y)$ est équivariante.\\
Soient $g\in G,\,\,\mu \in \mathbb{P}(X)\,\,f\in C(Y)$, montrons que $g(\alpha^{\star}\mu)=\alpha^{\star}(g\mu)$. Notons $\lambda$ l'action de $G$ sur $Y$. Nous avons d'une part:
$$g\alpha^{\star}\mu(f)=(\alpha^{\star}\mu)(f\circ \lambda_g)=\mu(f\circ\lambda_{g}\circ\alpha)$$ et d'autre part $$\alpha^{\star}(g\mu)(f)=g\mu(f\circ\alpha)=\mu(f\circ\alpha\circ \pi_g).$$ Comme $\alpha : X\longrightarrow Y$ est équivariante, on a: $\alpha\circ\pi_{g}=\lambda_{g}\circ\alpha$. D'où le résultat.

\end{prof}
Le lemme suivant est bien connu et est d'une importance capitale pour la suite.
\begin{lemma}\label{lembary}(\cite{krieger}, Lemme $1.32$)
Soit $K$ un compact convexe non vide d'un espace vectoriel topologique séparé localement convexe $E$. Soit $\mu \in \mathbb{P}(K)$. Alors
\begin{enumerate}
  \item Il existe un unique $b=b(\mu)\in K$ vérifiant $f(b)=\int_{K}f(k)d\mu(k)$ quel que soit $f\in E^{\star}$. Le point $b(\mu)$ est appelé $\mu$-barycentre de $K$.
  \item L'application \[\begin{array}{lll}
b:& \mathbb{P}(K) &\longrightarrow K\\
&\mu &\longmapsto b(\mu)\\
\end{array}\]
   est continue.
  \item On a $A(b(\mu))=b(A\ast\mu)$ pour toute application affine continue $A:K\longrightarrow K$, où $A\ast\mu$ désigne la mesure image de $\mu$ par $A$, c'est-à-dire $A\ast\mu(B)=\mu(A^{-1}(B))$.
\end{enumerate}
\end{lemma}
Pour plus de détails et résultats concernant la notion de \index{barycentre} barycentre voir \cite{cho}.

\begin{lemma} \label{th1}
Pour une action affine et continue d'un groupe topologique $G$ sur un sous-ensemble compact et convexe $X$ d'un espace localement convexe $E$, les propositions suivantes sont équivalentes:
\begin{enumerate}
  \item L'action possède un point fixe.
  \item L'action possède une mesure de probabilité invariante.
\end{enumerate}
\end{lemma}
\begin{prof}
\begin{itemize}
  \item [$1)\Rightarrow 2)$] Soit $x_{0}\in X$ un point fixe pour l'action continue de $G$ sur $X$.
 La mesure de Dirac $\delta_{x_{0}}$ de support $x_{0}$ est une mesure invariante.
  \item [$2)\Rightarrow 1)$]
  Supposons que $G$ agit continûment et affinement sur un compact convexe $X$ d'un espace vectoriel topologique séparé localement convexe $E$. D'après $2.$ l'espace $X$ admet une mesure de probabilité $\mu \in \mathbb{P}(X)$ qui est $G$-invariante. Notons $b=b(\mu)$ le $\mu$-barycentre de $X$ (voir lemme \ref{lembary}) et montrons que $b$ est un point fixe pour l'action continue de $G$ sur $X$. Notons $A_g:X\longrightarrow X$ l'application définie par $A_{g}(x)=gx$ quel que soient $x\in X$ et $g\in G$. Puisque $\mu$ est $G$-invariante, on a: $A_{g}\ast\mu=\mu$. Il en résulte en utilisant le lemme \ref{lembary} $$A_{g}(b(\mu))=b(A_{g}\ast\mu)=b(\mu)$$ pour tout $g\in G$. Donc $b(\mu)$ est un point fixe pour l'action de $G$ sur $X$.
\end{itemize}
\end{prof}

\begin{theorem}
Soit $G$ un groupe topologique. Les propriétés suivantes sont équivalentes:
\begin{enumerate}
\item Toute action affine et continue de $G$ sur un sous-espace compact et convexe $K$ d'un espace localement convexe possède un point fixe.
\item Toute action continue de $G$ sur un espace compact $X$ possède une mesure de probabilité invariante.
\item Il existe une mesure de probabilité invariante sur le compactifié de Samuel équivariant $S(G)$ de $G$.
\item $G$ est moyennable
\end{enumerate}
\end{theorem}
\begin{prof}
\begin{itemize}
\item [$1)\Longrightarrow 2)$] Soit $X$ un $G$-espace compact. Notons $\pi:G\times X\longrightarrow X$ l'action continue de $G$ sur $X$. Par le lemme \ref{prolon}, l'action continue $\pi$ se prolonge en une action affine et continue $\pi^{\star}$ de $G$ sur $\mathbb{P}(X)$. D'après $1$, cette dernière action possède un point fixe qui est une mesure de probabilité invariante pour l'action initiale de $G$ sur $X$.
  \item [$2)\Longrightarrow 1)$] est une conséquence du lemme \ref{th1}.
  \item [$2)\Longrightarrow 3)$] évident puisque $S(G)$ est compact et $G$ opère continûment sur $S(G)$
  \item [$3)\Longrightarrow 4)$] D'après le théorème de dualité de Gelfand (théorème \ref{gel}), il existe un isomorphisme isométrique de $E=RUCB(G)$ sur l'algèbre $C(S(G))$ des fonctions continues sur $S(G)$. Ainsi une mesure sur $S(G)$ correspond à une moyenne sur $E$. Donc $G$ est moyennable
  \item [$4)\Longrightarrow 2)$] Supposons $G$ moyennable. Nous noterons par $m$ une moyenne invariante sur $RUCB(G)$. Soit $X$ un $G$-espace compact. Soit $x_{0}\in X$ fixé et soit

      \[\begin{array}{lll}
t:& G &\longrightarrow X\\
&g &\longmapsto gx_{0}\\
\end{array}\]
 l'application orbital correspondante. Pour tout $f\in C(X)$, considérons l'application $f^{\ast}=f\circ t$. On a
      $f^{\ast}\in RUCB(G)$. \\
      En effet, l'action $G\times X\longrightarrow X$ étant continue, il existe comme dans la preuve du théorème \ref{grea} un voisinage $V$ de $e$ dans $G$ tel que $g\in V\,\Longrightarrow\,\,|f(g.x)-f(x)|<\varepsilon$ pour tout $x\in X$ et pour tout $\varepsilon> 0$.\\
 Si $gh^{-1}\in V$, alors $$|f^{\ast}(g)-f^{\ast}(h)|=|f(g.x_{0})-f(h.x_{0})|=|f(gh^{-1}(h.x_{0}))-f(h.x_{0})|<\varepsilon.$$
Donc $f^{\ast}\in E$ puisque $f^{\ast}$ est clairement bornée car $X$ est compact.
La mesure de probabilité $\mu_{m}$ définie sur $X$ par $$\mu_{m}(f)=m(f^{\ast})$$ est invariante.
\end{itemize}
\end{prof}
\subsection{Exemples de Groupes Moyennables}
Le théorème suivant apparait dans \cite{frem}.
\begin{theorem}\label{moyexample}
Soit $G$ un groupe topologique.
\begin{enumerate}
  \item Si $G$ est moyennable et si $H$ est un groupe topologique tel qu'il existe un homomorphisme surjectif continu de $G$ sur $H$, alors $H$ est moyennable.
  \item Si $A$ un sous-espace dense de $G$ tel que tout sous-ensemble fini de $A$ est contenu dans un sous-groupe moyennable de $G$, alors $G$ est moyennable.
   \item Si $H$ est un sous-groupe normal de $G$ et si $H$ et $G/ H$ sont moyennables, alors $G$ est moyennable.
  \item Si $G$ est abélien, alors $G$ est moyennable.
\end{enumerate}
\end{theorem}

\begin{prof}
\begin{enumerate}
  \item Soit $\phi:G\longrightarrow H$ un homomorphisme surjectif continu. Soient $X$ un espace compact et $\bullet:H\times X\longrightarrow X$ une action continue de $H$ sur $X$. Pour $g\in G$ et $x\in X$, posons $$g\bullet_{1}x=\phi(g)\bullet x.$$ $\bullet_{1}$ est une action continue de $G$ sur $X$. Par hypothèse, il existe une mesure de probabilité borélienne invariante $\mu$ sur $X$. Puisque $\phi(G)=H,\,\,\mu$ est aussi $H$-invariante. Donc $H$ est moyennable.
  \item Soit $X$ un $G$-espace compact. Pour $g\in G$ et $\mu\in \mathbb{P}(X)$, notons $g\bullet\mu$ l'action continue de $G$ sur $\mathbb{P}(X)$. L'espace $$Q_{g}=\{\mu:\,\mu\in \mathbb{P}(X),\,\,g\bullet \mu=\mu\}$$
      est un fermé de $\mathbb{P}(X)$ pour tout $g\in G$ et $$G_{\mu}=\{g:\,g\in G,\,\,g\bullet \mu=\mu\}$$
      est un fermé de $G$ pour tout $\mu\in \mathbb{P}(X)$.\\
      Pour tout sous-ensemble fini $I$ de $A$, il existe un sous-groupe moyennable $H_{I}$ de $G$ contenant $I$. La restriction $H_{I}\times X\longrightarrow X$ de l'action continue de $G$ sur $X$ est continue. Ainsi, il existe une mesure borélienne $H_{I}$-invariante. De plus, $\underset{a\in I}{\overset{}{\bigcap}}Q_{a}\supseteq \underset{a\in H_{I}}{\overset{}{\bigcap}}Q_{a}$. Donc $\underset{a\in I}{\overset{}{\bigcap}}Q_{a}$ est non vide. Comme $\mathbb{P}(X)$ est compact, $\underset{a\in A}{\overset{}{\bigcap}}Q_{a}$ est non vide. Soit $\mu\in \underset{a\in A}{\overset{}{\bigcap}}Q_{a}$. Puisque $G_{\mu}$ contient le sous-espace dense $A$, il coïncide avec $G$. Donc $\mu$ est $G$-invariant.

  \item  Soit $X$ un $G$-espace compact. Notons
\[\begin{array}{lll}
\pi:& G\times X &\longrightarrow X\\
&(g,x) &\longmapsto gx\\
\end{array}\]
l'action continue de $G$ sur $X$, et pour tout $g\in G$, notons

\[\begin{array}{lll}
\pi_{g}:& X &\longrightarrow X\\
&x &\longmapsto gx\\
\end{array}\]
homéomorphisme de $X$ sur lui-même.\\
Pour $g\in G$ et $\mu\in \mathbb{P}(X)$, notons comme précédement $g\bullet\mu$ l'action continue de $G$ sur $\mathbb{P}(X)$. Pour $h\in H$, posons
      $$Q=\{\mu:\,h\bullet \mu=\mu\}.$$
      $Q$ est un sous-espace fermé de $\mathbb{P}(X)$. Puisque $H$ est moyennable, $Q$ est non-vide. On a $g\bullet \mu\in Q$ pour tous $\mu\in Q$ et $g\in G$. En effet, si $h\in H$, on a:

      $$ h\bullet(g\bullet \mu)=(hg)\bullet \mu=(gg^{-1}hg)\bullet \mu=g\bullet(g^{-1}hg)\bullet\mu)=g\bullet \mu,$$
      puisque $H$ est normal, $g^{-1}hg\in H$. Ceci nous permet de définir une action continue de $G$ sur l'espace compact $Q$.\\
      Si $h\in H$ et $g\in G$, alors $g\bullet \mu =(gh)\bullet \mu$ pour tout $\mu\in Q$. Ainsi, l'application
  \[\begin{array}{lll}
& G/H\times Q &\longrightarrow Q\\
&(gH,\mu) &\longmapsto (gH).\mu=g\bullet \mu\\
\end{array}\]
définie une action continue de $G/H$ sur $Q$. En effet, \\
Puisque $G/H$ est moyennable, il existe une mesure de probabilité borélienne $G/H$-invariante $\lambda$ sur $Q$.\\
Pour $f\in C(X)$, posons $$p(f)=\int_{Q}(\int_{X}f(x)d\mu(x))d\lambda(\mu)$$
$p$ est bien défini car l'application $\mu\longmapsto \int_{X}f(x)d\mu(x)$ est continue pour la topologie vague sur $Q$. $p$ est une forme linéaire, $p(f)\geq 0$ pour tout $f\geq 0$ et $p(1_{X})=1$. Ainsi, il existe une mesure borélienne $\nu$ sur $X$ telle que $p(f)=\int fd\nu$ pour tout $f\in C(X)$. Si $g\in G$, alors

 \begin{tabular}{lll}

$\int fd(g\bullet \nu)$
& $=$ &
$\int f\circ \pi_{g}d\nu=p(f\circ \pi_{g})$\\
& $=$ & $ \int_{Q}(\int_{X}(f\circ \pi_{x})d\mu(x))d\lambda(\mu)$ \\
& $= $ & $ \int_{Q}\int_{X}fd\mu(g\bullet x)d\lambda(\mu)$ \\
& $=$ & $\int_{Q}\int_{X}fd\mu((gH).x)d\lambda(\mu)$\\
& $=$ & $\int_{Q}\int_{X}fd\mu d\lambda(\mu)$,\,\,\,($\lambda$ \,\,est\,\,$G/H$-invariante)\\
& $=$ & $\int fd\nu$
\end{tabular}

Pour tout $f\in C(X)$.\\
Ainsi $g\bullet \nu=\nu$ pour tout $g\in G$. Donc $G$ est moyennable.
  \item Commençons par rappeler que le théorème du point fixe de Kakutani (\cite{bourb1}) affirme que si $K$ est un compact convexe d'un espace localement convexe
$E$, alors toute application affine continue $T: K\longrightarrow K$ admet un point fixe.\\
Maintenant, soit $G\times K\longrightarrow K$ une action affine continue de $G$ commutatif sur un compact convexe $K$. L'espace $$K_{g}=\{x\in K:\,\,gx=x\}$$ est non vide (par le théorème du point fixe de Kakutani) et compact pour tout $g\in G$.
Si $g^{\prime}\in G$, alors pour tout $x\in K_{g}$, on a: $$g(g^{\prime}x)=(gg^{\prime})x=(g^{\prime}g)x=g^{\prime}(gx)=g^{\prime}x.$$
Ainsi $K_{g}\cap K_{g^{\prime}}\neq\emptyset$ pour tout $g,g^{\prime}\in G$. \\
Par récurrence on obtient que les intersections finies des $K_{g}, \,g\in G$ sont non vides, et donc par compacité $\underset{g\in G}{\overset{}{\bigcap}}K_{g}\neq \emptyset$. Si $x\in \underset{g\in G}{\overset{}{\bigcap}}K_{g}\neq \emptyset$, alors $x$ est un point fixe pour l'action affine précédente.
\end{enumerate}
\end{prof}

 Terminons cette section par deux exemples de groupes non-moyennables.
 Le premier est un exemple classique d'un groupe discret non-moyennable construit par Von Neumann (\cite{von})
\begin{enumerate}
  \item Le groupe libre à deux générateurs\index{groupe!libre}  non-abélien $F_2$  muni de la topologie discrete n'est pas moyennable.\\
  Commençons par rapeller que le groupe libre non-abélien à deux générateurs est l'ensemble des mots (simplifiés) de longueurs fini construit à partir de l'alphabet $a,a^{-1},b,b^{-1}$ et comprenant comme identité le mot vide $\phi$. L'opération de groupe ici est la concaténation i.e si $w,z\in F_2$ sont deux mots  avec $w=w_1w_2...w_m$ et $z=z_1z_2...z_n$, alors $w.z=w_1w_2...w_mz_1z_2...z_n$ qui est ensuite réduit de façon approprié. Montrons que $F_2$ muni de la topologie discrète est non-moyennable.\\
Supposons le contraire et soit $\mu$ une mesure de probabilité invariante sur $F_2$. Pour $x\in F_2$, notons $E_x=\{y\in F_2:\,y\, \text{est un mot réduit débutant par}\,\, x\}$. Alors $F_2=\{\phi\}\cup E_a\cup E_{a^{-1}}\cup E_b\cup E_{b^{-1}}$. On a aussi, $aE_{a^{-1}}=E_{a^{-1}}\cup E_{b}\cup E_{b^{-1}}\cup\{\phi\}$. Ainsi, $F_2=E_{a}\cup aE_{a^{-1}}$. De la même manière, on a: $F_2=E_{b}\cup bE_{b^{-1}}$. Posons
$A=E_{a}\cup E_{a^{-1}}$ et $B=E_{b}\cup E_{b^{-1}}$. En utilisant l'invariance de $\mu$, nous obtenons:
$\mu(A)=\mu(E_a)+\mu(E_{a^{-1}})=\mu(E_a)+\mu(aE_{a^{-1}})=1$. De la même manière, $\mu(B)=1$. Maintenant, $\mu(F_2)=\mu(A\cup B)=\mu(A)+\mu(B)=2$. Absurde.\\

Le deuxième exemple d'un groupe non-moyennable est plus récent.

  \item Le groupe des automorphismes\index{groupe!des automorphismes} mesurables non-singuliers préservant la mesure $\mu$ d'un espace de Borel standard $(X,\mu)$.
\end{enumerate}

Introduisons les notions nécéssaires pour la bonne compréhension de cet exemple.\\
Soit $(X,\Sigma,\mu)$ un espace de Borel standard muni d'une mesure de probabilité borélienne non-atomique $\mu$.
Rappelons les notions suivantes:
\begin{enumerate}
\item On appelle $\mu$-atome, une partie $A\in \Sigma$ tel que $\mu(A)> 0$ et $\forall\,B\in \Sigma$,
                    $B\subset A\Longrightarrow \mu(B)=0$ ou $\mu(A\setminus B)=0$.
\item une mesure $\mu$ est dite diffuse ou non-atomique si elle n'admet aucun $\mu$-atome.
\item Un automorphisme $T:X\longrightarrow X$ est dit singulier par rapport à $\mu$ si $$\mu(B)=0\Longrightarrow \mu(T(B))=\mu(T^{-1}(B))=0$$

\end{enumerate}
Rappelons également qu'on peut définir sur $(X,\Sigma,\mu)$ la relation d'équivalence:\\
$\mu_{1}\sim \mu_{2}\Longleftrightarrow\mu_{1}$ et $\mu_{2}$ ont les mêmes ensembles de mesures nulles.\\
De façon précise, $\mu_{1}\sim \mu_{2}\Longleftrightarrow(\mu_{1}(A)=0 \Longleftrightarrow\mu_{2}(A)=0)$.\\
Notons $Aut^{\star}(X,\mu)$ le groupe de tous les automorphismes non-singuliers $\tau:X\longrightarrow X$ qui préserve la classe de la mesure $\mu$ ($\tau_{\star}\mu \sim \mu$).

\begin{remark}
Le groupe $Aut(X,\mu)$ de tous les automorphismes mesurables non-singuliers préservant la mesure $\mu$ ($\tau_{\star}\mu=\mu$) est un sous-groupe propre de $Aut^{\star}(X,\mu)$.\\
En effet, prenons $X=[0,1],\,\,\tau(t)=t^{2}$ et $A=[0,\frac{1}{3}]$.
 Nous avons $\tau(A)=[0,\frac{1}{9}]$. Ainsi, $\tau \in Aut^{\star}(X,\mu)$ et $\tau \notin Aut(X,\mu)$
\end{remark}
L'application $d(\tau,\sigma)=\mu\{x\in X:\,\tau(x)\neq\sigma(x)\}$ définie une distance sur $Aut^{\star}(X,\mu)$ invariante à gauche.

La topologie induite par cette distance est compatible avec la structure de groupe de $Aut^{\star}(X,\mu)$. (\cite{vpbresil} page $119$)\\
La topologie induite par cette distance sur $Aut^{\star}(X,\mu)$ est appelée topologie uniforme\index{topologie!uniforme}.
\begin{theorem}(Giordano et Pestov \cite{giopes1,giopes2})
Le groupe $Aut(X,\mu)$ muni de la topologie uniforme n'est pas moyennable.
\end{theorem}
\begin{remark}
Contrairement à l'exemple classique de groupe non-moyennable établi par Von Neumann, l'exemple de Giordano et Pestov n'est pas un groupe localement compact.
\end{remark}

La question suivante reste une question ouverte.
  \begin{question}
  Le groupe $Aut^{\star}(X,\mu)$ muni de la topologie uniforme est-il moyennable? (Giordano et Pestov \cite{giopes1})
\end{question}
\begin{remark}
 \begin{enumerate}
   \item De nombreux auteurs utilisent la terminologie "moyennable" pour désigner un groupe moyennable pour la topologie discrète. Le danger de cette terminologie est que de nombreux résultats concernant les groupes moyennables discrets ne se généralisent pas aux groupes topologiques moyennables, même pas aux groupes localement compacts moyennables comme le montre les observations suivants:
       \begin{enumerate}
        \item Tout sous-groupe d'un groupe discret moyennable est moyennable. (\cite{paterson}, Proposition $0.16$ page $14$)
   \item Tout sous-groupe fermé d'un groupe localement compact moyennable est moyennable. (\cite{paterson}, Proposition $1.12$ page $31$)
   \item En dehors du cas des groupes localement compacts, les sous-groupes fermés d'un groupe localement compact moyennable ne sont en général pas moyennable. En effet, le groupe $Aut\,(\mathbb{Q},\leq)$ des bijections de $\mathbb{Q}$ dans lui-même qui préservent l'ordre muni de la topologie de la convergence simple est moyennable(car il est extrêmement moyennable). De plus, le groupe libre(non-moyennable) $F_{2}$ est isomorphe à un sous-groupe fermé de $Aut\,(\mathbb{Q},\leq)$ (\cite{vp2}).
       \end{enumerate}

    \item Les références \cite{paterson} et \cite{pier} contiennent une revue complète sur les groupes moyennables localement compacts.
 \end{enumerate}
 \end{remark}
\chapter{Espaces test pour la moyennabilité et la moyennabilité extrême}

\section{Espace test pour la moyennabilité des groupes polonais}
En réponse à une question de Grigorchuk, Giordano et de la Harpe (\cite{gio}, $1997$) ont montré qu'un groupe discret dénombrable
$G$ est moyennable si et seulement si toute action continue de $G$ sur l'ensemble
de Cantor $D^{\aleph_{0}}$ possède une mesure de probabilité invariante. Ce résultat permet de détecter la moyennabilité d'une classe de groupes topologiques en utilisant un seul espace compact. D'où la notion d'espace test. On peut dire que l'ensemble de
Cantor est un espace test\index{espace!test} pour la moyennabilité des groupes discrets dénombrables. Dans le
même sens, Bogatyi et Ferdorchuk (\cite{boga}, $2007$) ont répondu à une question de \cite{gio} et ont demontré que le
cube de Hilbert $I^{\aleph_0}$ est également un espace test pour la moyennabilité des groupes discrets
dénombrables.\\
Dans cette section, nous allons établir que le cube de Hilbert reste un espace test pour la moyennabilité des groupes polonais.\\
 Commençons par rappeler le théorème de Keller suivant qui est d'une très grande importance dans nos démonstrations.
\begin{theorem}(Keller\cite{bes} Théorème $3.1$)\index{théorème!de Keller}
Toute partie compacte, convexe, et de dimension infinie(admet une partie libre infini) d'un espace localement convexe métrisable est homéomorphe au cube de Hilbert $I^{\aleph_{0}}$.
\end{theorem}
\begin{lemma}\label{lemgel}
Soient $G$ un groupe topologique et $A$ une algèbre de Banach. Notons $X_A$ le spectre de $A$. Si $G$
 opère continûment sur $A$ par automorphismes, alors $G$ opère continûment sur $X_A$.
\end{lemma}

\begin{prof}
L'application
\[ \xymatrix@C=2cm@R=0.3em{
 \eta:  G\times A^{\star}  \ar[r] &  A^{\star} &  \\
  (g,\phi) \ar@{|->}[r] & g.\phi : A \ar[r] &  \mathbb{C} \\
& \,\,\,\,\,\,\,\,\,\,\,\,\,\,\,x \ar@{|->}[r] &  g.\phi(x)=\phi(g^{-1}.x) }
\]
définie une action continue de $G$ sur $A^{\star}$. $X_A$ étant stable par cette action, on peut considérer la restriction

\[ \xymatrix@C=2cm@R=0.3em{
 \eta_{X_A}:  G\times X_{A}  \ar[r] &  X_{A} &  \\
  (g,\phi) \ar@{|->}[r] & g.\phi : A \ar[r] &  \mathbb{C} \\
& \,\,\,\,\,\,\,\,\,\,\,\,\,\,\,x \ar@{|->}[r] &  g.\phi(x)=\phi(g^{-1}.x) }
\]
 de celle-ci à $X_A$.\\
Montrons qu'en munissant $X_{A}$ de la topologie vague, $\eta_{X_A}$ est continue.\\
Soit $\phi\in X_{A}$, alors pour tout $x\in A$, nous avons: $|\phi(x)|\leq \|x\|$.\\
Comme $x\longmapsto g^{-1}.x$ est continue, il existe un voisinage symétrique de $g$ que nous noterons $U$ tel que:
 $$h\in U \Longrightarrow \|g^{-1}.x-h^{-1}.x\|< \frac{\varepsilon}{2},\,\,\,\,\forall\,x\in A$$
 Prenons $$V=\{f\in X_{A},\,\,\,|(f-\phi)(x_{i})|<\frac{\varepsilon}{2}\,,\,\,i=1,2,...,k\}.\,\,x_{1},...,x_{k}\in A$$
 Pour $h\in U$ et pour $\psi\in V$, nous avons:\\
 \[\begin{array}{llllll}|g.\psi(x_{i})-h.\phi(x_{i})|&=
|\psi(g^{-1}.x_{i})-\phi(h^{-1}.x_{i})|
\\&= |\psi(g^{-1}.x_{i})-\psi(h^{-1}.x_{i})+\psi(h^{-1}.x_{i})-\phi(h^{-1}.x_{i})|\\&
\leq |\psi(g^{-1}.x_{i}-h^{-1}.x_{i})|+|(\psi-\phi)(h^{-1}.x_{i})|
 \\
\end{array}\] \\
Comme $\psi\in V\subset X_{A}$, on a: $$|\psi(g^{-1}.x_{i}-h^{-1}.x_{i})|<\|g^{-1}.x_{i}-h^{-1}.x_{i}\|<\frac{\varepsilon}{2}.$$ De plus,
$$\psi\in V\Longrightarrow|(\psi-\phi)(h^{-1}.x_{i})| < \frac{\varepsilon}{2}.$$
Donc $|g.\psi(x_{i})-h.\phi(x_{i})|<\varepsilon$.
\end{prof}
\begin{lemma}\label{lemstar}
Soient $A$ et $B$ deux $C^{\star}$-algèbres commutatives et unifères et soit $\lambda$ une application de $A$ dans $B$. Notons $X_A$
(respectivement $X_B$) le spectre de l'algèbre $A$(respectivement $B$). Considérons l'application $\psi_{\lambda}:X_B\longrightarrow X_A$ définie par $\psi_{\lambda}(\chi)=\chi\circ \lambda$. Les propositions suivantes sont équivalentes:
\begin{enumerate}
  \item $\lambda$ est injective
  \item $\psi_{\lambda}$ est surjective.
\end{enumerate}
\end{lemma}
\begin{remark}
D'après le théorème de Gelfand, le lemme \ref{lemstar} est équivalent au lemme suivant:
\begin{lemma} Soient $X$ et $Y$ deux espaces compacts et soit $f$ une application continue de $X$ dans $Y$.
Notons $\widehat{f}:C(Y)\longrightarrow C(X)$ l'application définie par: $\widehat{f}(g)=g\circ f$.
Les propositions suivantes sont équivalentes:
\begin{enumerate}
  \item $f$ est surjective
  \item $\widehat{f}$ est injective
\end{enumerate}
\end{lemma}
\end{remark}
\begin{prof}
\begin{enumerate}
\item [$1)\Longrightarrow 2)$] Soient $g,h\in C(Y)$ tel que $h\neq g$. Alors il existe $y\in Y$ tel que $g(y)\neq h(y)$.\\ Comme $f$ est surjective, il existe $x\in X$ tel que $f(x)=y$. On a donc $$\widehat{f}(g(x))=g(y)\neq h(y)=\widehat{f}(h(x)).$$
    \item [$2)\Longrightarrow 1)$] Supposons $f$ non surjective et montrons que $\widehat{f}$ n'est pas injective.\\
    Soit $y\in Y\setminus f(X)$. Comme $f$ est continue, $f(X)$ est compact donc fermé.\\ Comme $Y$ est compact, $Y$ est
     complètement régulier. Ainsi, il existe $h\in C(Y)$ tel que $h\mid_{f(X)}=0$ et $h(y)=1$.
      Considérons l'application $g\in C(Y)$ définie par $g(y)=0$ pour tout $y\in Y$. Nous avons $\widehat{f}(g)=0=\widehat{f}(h)$ sur $X$. Mais $g\neq h$.
 \end{enumerate}
\end{prof}

 \begin{lemma}\label{lemmesurinverse}
S'il existe une mesure de probabilité invariante sur chaque $G$-espace dans un système inverse des $G$-espaces compacts, alors il existe une mesure de probabilité invariante sur la limite inverse correspondante.
\end{lemma}
\begin{prof}
Soit $(X_{\alpha},\pi_{\alpha \beta}, I)$ un système inverse de $G$-espaces compacts. Notons $X=\underset{\longleftarrow}{\overset{}{\lim}}\,X_{\alpha}$. Par le lemme \ref{lemaction}, ce système inverse permet d'obtenir un système inverse $(\mathbb{P}(X_{\alpha}),\pi_{\alpha \beta}^{\star}, I)$. On a: $\mathbb{P}(X)=\underset{\longleftarrow}{\overset{}{\lim}}\,\mathbb{P}(X_{\alpha})$. En effet, notons $\pi_{\alpha}:X\longrightarrow X_{\alpha}$ la restriction de la projection $pr_{\alpha}$ à $X$. Par le lemme \ref{lemaction}, il existe

\[\begin{array}{llll}
\pi_{\alpha}^{\star}: &\mathbb{P}(X) &\longrightarrow & \mathbb{P}(X_{\alpha}) \\
&\mu &\longmapsto & \mu\circ \alpha ^{-1}
\end{array} \]

L'application
 \[\begin{array}{llll}
T: &\mathbb{P}(X) &\longrightarrow & \underset{\longleftarrow}{\overset{}{\lim}}\,\mathbb{P}(X_{\alpha}) \\
&\mu &\longmapsto & (\pi_{\alpha}^{\star}\mu)_{\alpha \in I}
\end{array} \]
est un homéomorphisme. En effet, $T$ est surjective, puisque chaque $\pi_{\alpha}^{\star}$ est surjective.\\
 Si $T(\mu)=T(\nu)$, alors $\int_{X}(f\circ \pi_{\alpha})d\mu=\int_{X} (f\circ \pi_{\alpha}) d\nu$ pour tout $f\in C(X_{\alpha})$. Ainsi
$\mu=\nu$.\\
 Notons $\mathbb{P}_{inv}(X_\alpha)$ l'espace des mesures de probabilités invariantes sur $X_\alpha$ et considérons le sous-système $(\mathbb{P}_{inv}(X_{\alpha}),\pi_{\alpha \beta}^{\star},I)$ de $(\mathbb{P}(X_{\alpha}),\pi_{\alpha \beta}^{\star}, I)$. Notons $\mathbb{P}_{inv}(X)=\underset{\longleftarrow}{\overset{}{\lim}}\,\mathbb{P}_{inv}(X_{\alpha})$. Puisque $\mathbb{P}_{inv}(X_\alpha)$ est compact et non vide pour tout $\alpha$, on a: $\mathbb{P}_{inv}(X)\neq \emptyset$ (propriété d'intersection finie). Si $\mu\in \mathbb{P}_{inv}(X)$, alors $\mu$ est une mesure de probabilité invariante sur $X$.
\end{prof}
\begin{definition}
Un groupe topologique $G$ est dit polonais s'il est séparable et si sa topologie est définie par une métrique complète.
\end{definition}
La décomposition suivante du compactifié de Samuel par un systéme projectif est très utile dans nos démonstrations.
\begin{theorem}\label{poloinverse}
Si $G$ est un groupe polonais, alors il existe un système projectif de $G$-espaces compacts et métrisables
$(X_{\alpha},\pi_{\alpha \beta}, I)$ tel que $S(G)=\underset{\longleftarrow}{\overset{}{\lim}}\,X_{\alpha}$
\end{theorem}
\begin{prof}
Fixons une partie dénombrable dense $D\subseteq G$. Soit $(I,\preceq)$ l'ensemble ordonné suivant: $I$ est l'ensemble de toutes les sous $C^{\star}$-algèbres unifères séparables, fermées et $G$-invariantes de $RUCB(G)$.\\
 Soit $A,B\in I$, définissons l'ordre par $A\preceq B\,\Longleftrightarrow\,A\subseteq B$.\\
 Pour $A\in I$, l'espace de Gelfand $X_A$ de $A$ est compact et métrisable, puisque $A$ est séparable. Par le lemme \ref{lemgel}, $G$ opère continûment sur $X_A$. Nous pouvons identifier tout $g\in G$ à un élément de $X_{A}$ de telle façon que $G$ soit un sous-espace dense de $X_{A}$. Ainsi $(X_{A},1_{G})$ est un compactifié équivariant métrisable de $G$.\\
  Dans la suite, nous allons identifier par le théorème de Gelfand(théorème \ref{gel}) $A$ à $C(X_A)$ via l'application
  \[ \xymatrix@C=2cm@R=0.3em{
   A  \ar[r] &  C(X_A) &  \\
 x \ar@{|->}[r] &  \,\,\,\,\,\,\,\widehat{x}_A: X_A \ar[r] &  \mathbb{C} \\
& \,\,\,\,\,\,\,\,\,\,\,\,\,\,\,f \ar@{|->}[r] &  \widehat{x}_A(f)=f(x) }
\]
dite application de Gelfand.\\
  Si $A\subseteq B$, alors l'injection canonique de $i:A\longrightarrow B$ permet de définir par le lemme \ref{lemstar} une surjection $\pi_{AB}:X_B\longrightarrow X_A$ continue telle que pour $x\in A,\,\,\,\widehat{x}_B=\widehat{x}_A\circ \pi_{AB}$. Ainsi, pour $\varphi \in X_B$, on a: $\widehat{x}_A(\pi_{AB}(\varphi))=\widehat{x}_B(\varphi)$ ou $\pi_{AB}(\varphi)(x)=\varphi(x)$ i.e $\pi_{AB}(\varphi)=\varphi|_{A}$. De la même façon, il existe une surjection continue $\pi_{A}: S(G)\longrightarrow X_A$ définie par $\pi_{A}(\varphi)=\varphi|_{A}$ pour tout $A\in I$. D'où $\pi_{A}=\pi_{AB}\circ \pi_B$ pour $A\preceq B$. Comme $\pi_A(1_G)=\pi_{AB}(1_G)=1_G$, l'application

  \[\begin{array}{lll}
\chi:& (S(G),1_G) &\longrightarrow \underset{\longleftarrow}{\overset{}{\lim}}\,(X_{A},1_G)\\
&\varphi &\longmapsto (\pi_{A}(\varphi))_{A\in I}\\
\end{array}\]
est un homomorphisme de compactifiés équivariants de $G$. \\
  Soit $(\varphi_A)_{A\in I} \in \underset{\longleftarrow}{\overset{}{\lim}}\,X_{A}$, nous avons pour $A\preceq B,\,\,\varphi_{A}=\varphi_{B}|_{A}$, donc il existe un unique $\varphi \in S(G)$ tel que $\pi_{A}(\varphi)=\varphi_A$. Ainsi, $ \varphi\longmapsto (\pi_{A}(\varphi))_{A\in I}$ est un isomorphisme entre les compactifiés équivariants $(S(G),1_G)$ et $\underset{\longleftarrow}{\overset{}{\lim}}\,(X_{A},1_G)$ de $G$.
 \end{prof}

\begin{proposition}\label{pr4}
Un groupe polonais $G$ est moyennable si et seulement si toute action continue de $G$ sur un espace compact et métrisable, possède une mesure de probabilité borélienne invariante.
\end{proposition}

\begin{prof}
La nécessité est évidente.\\ Montrons la suffisance. Par le théorème \ref{poloinverse}, il existe une système projectif de $G$-espaces compacts et métrisables
$(X_{\alpha},\pi_{\alpha \beta}, I)$ tel que $S(G)=\underset{\longleftarrow}{\overset{}{\lim}}\,X_{\alpha}$. Par hypothèse, il existe sur chaque $G$-espace compact métrisable $X_\alpha$ une mesure de probabilité invariante. Donc il existe une mesure de probabilité invariante sur $S(G)$ par le lemme \ref{lemmesurinverse} et $G$ est moyennable.
\end{prof}
\begin{remark}\label{remar}
\begin{enumerate}
  \item Le résultat précédent est bien connu pour la classe des groupes dénombrables et discrets.
  \item Évidemment, on peut supposer sans perte de généralité que tous les $G$-espaces $X$ dans la proposition \ref{pr4} sont infinis. C'est le cas si et seulement si $G$ est infini.
\end{enumerate}
\end{remark}

\begin{theorem}\label{main1}
Un groupe polonais $G$ est moyennable si et seulement si toute action continue de $G$ sur le cube de Hilbert $I^{\aleph_{0}}$ possède une mesure de probabilité borélienne invariante.
\end{theorem}
\begin{prof}La nécessité est évidente. Montrons la suffisance.\\
Soit $X$ un $G$-espace compact et métrisable. Par la remarque \ref{remar}, on peut supposer que $X$ est infini. $\mathbb{P}(X)$ est donc un sous-espace compact métrisable de dimension infinie de $\mathbb{R}^{C(X)}$. Par le théorème de Keller (voir \cite{bes}), $\mathbb{P}(X)$ est homéomorphe au cube de Hilbert $I^{\aleph_{0}}$. Ainsi l'action de $G$ sur $\mathbb{P}(X)$ possède une mesure de probabilité borélienne invariante. L'action de $G$ sur $\mathbb{P}(X)$ étant affine, elle possède par le lemme \ref{th1} un point fixe $\mu \in \mathbb{P}(X)$, qui est une mesure de probabilité borélienne invariante pour l'action initiale de $G$ sur $X$.
\end{prof}
\begin{remark}
L'idée d'utiliser le théorème de Keller dans le contexte dynamique est de Uspenskij, qui était  le premier à l'employer dans \cite{us4}.
\end{remark}

\section{Espace test pour la moyennabilité des groupes polonais non archimédiens}
Dans cette section, nous établissons que l'ensemble de Cantor reste un espace test pour la moyennabilité des groupes polonais non archimédiens. Nous généralisons ainsi le résultat de Giordano et de la Harpe (\cite{gio}) à la classe des groupes polonais non archimédiens.
\subsection{Généralités sur les groupes polonais non archimédiens}
\begin{definition}
Un groupe topologique $G$ est dit {\em non archimédien}\index{groupe!non archimédien} s'il est séparé et s'il possède une base de voisinages de l'élément neutre formé de sous-groupes ouverts.
\end{definition}
L'ensemble des groupes polonais non archimédiens comprend:
\begin{enumerate}
\item les groupes localement compacts et totalement discontinus (\cite{hewit}),
  \item le groupe symétrique infini $S_{\infty}$ i.e le groupe de toutes les bijections de $\mathbb{N}$ dans $\mathbb{N}$ muni de la topologie de la convergence simple,
  \item le groupe $Homeo(D^{\aleph_{0}})$ des homéomorphismes de l'ensemble de Cantor muni de la topologie de la convergence uniforme,
  \item le groupe $Homeo(X)$ des homéomorphismes de $X$ sur lui-même. $X$ étant un espace compact de dimension de Lebesgue $0$ (\cite{megre}).
\end{enumerate}
Les groupes polonais non archimédiens peuvent être caractérisés par la dimension de Lebesgue de leur compactifié de Samuel équivariant de la manière suivante.
\begin{theorem}(Pestov \cite{vp2})
Soit $G$ un groupe polonais. $G$ est non archimédien si et seulement si les applications continues de $S(G)$ dans $D =\{0,1\}$ séparent les points de $S(G)$.
\end{theorem}
\begin{remark}
\begin{enumerate}
\item Les exemples $2.$ et $3.$ précédents sont traités avec plus de détails en annexe. Ces deux exemples permmettent de conclure que la classe des groupes discrets dénombrables est strictement incluse dans celle des groupes polonais non archimédiens.
  \item Les groupes polonais non archimédiens jouent un rôle important en logique où ils sont les groupes des automorphismes des structures de Fraïssé \cite{becker}.\\
  \item La classe des groupes polonais non archimédiens est stable par passage au sous-groupe, au produit et au quotient \cite{morris}.
\end{enumerate}
\end{remark}
Rappelons qu'un groupe topologique $G$ est universel pour une classe $\mathcal{C}$ de groupes topologiques si pour tout groupe topologique $H\in \mathcal{C}$, il existe un isomorphisme de groupes topologiques de $H$ sur un sous-groupe de $G$.\\
Les groupes polonais $S_{\infty}$ et $Homeo(D^{\aleph_{0}})$ sont en réalité universels pour la classe des groupes polonais non archimédiens.
\begin{theorem}
Soit $G$ un groupe polonais. Les propositions suivantes sont équivalentes:
\begin{enumerate}
  \item $G$ est non archimédien.
  \item  $G$ est isomorphe à un sous-groupe topologique de $S_{\infty}$.
\end{enumerate}
\end{theorem}
\begin{prof}
\begin{enumerate}
  \item [$1)\Longrightarrow 2)$]
  Soient $G$ un groupe polonais non archimédien et $(H_{i})_{i\in I}$ une base dénombrable de voisinage de $e$ formé de sous-groupes ouverts.
  Pour tout $i\in I$, $G$ opère continûment sur l'espace quotient $G/H_{i}$. Pour tout $i\in I$, l'espace quotient $G/H_{i}$ est discret, car $H_i$ est ouvert. Ainsi, l'espace $X=\underset{i=1}{\overset{\infty}{\bigcup}}(G/H_{i})$ est discret et dénombrable. De plus,
  $G$ opère continûment sur $X$. Cette action continue entraîne un plongement $G\hookrightarrow S_{X}=S_{\infty}$.
  \item [$2)\Longrightarrow 1)$] Puisque $S_{\infty}$ est non archimédien, alors tout sous-groupe $G^{\prime}$ de $S_{\infty}$ est non archimédien. Si $G$ est isomorphe à $G^{\prime}$, alors $G$ est également non archimédien.
 \end{enumerate}
\end{prof}
\begin{remark}
\begin{enumerate}
  \item Le résultat précédent apparait aussi dans \cite{becker}.
  \item Rappelons pour la suite que le théorème de Cantor (\cite{eng}) affirme que tout espace compact, métrisable, totalement discontinu et sans point isolé est homéomorphe à l'ensemble de Cantor. L'ensemble de Cantor désignera donc pour nous un espace possédant les propriétés précédentes.
  \item Le résultat suivant apparait dans \cite{megre}. Mais nous proposons ici une preuve différente.
\end{enumerate}

\end{remark}

\begin{theorem}
 Soit $G$ un groupe polonais. Les propositions suivantes sont équivalentes:
\begin{enumerate}
  \item $G$ est non archimédien.
  \item  $G$ est isomorphe à un sous-groupe topologique de $Homeo(D^{\aleph_{0}})$.
\end{enumerate}
\end{theorem}
\begin{prof}
\begin{enumerate}
\item [$1)\Longrightarrow 2)$] Soit $G$ un groupe polonais non archimédien. Il nous suffit de montrer que $S_{\infty}$ est isomorphe à un sous-groupe de $Homeo(D^{\aleph_{0}})$ pour conclure. Rappelons que $S_{\infty}$ opère continûmement sur $\mathbb{N}$.\\
    Notons $\alpha \mathbb{N}=\mathbb{N}\cup\{\infty\}$ le compactifié d'Alexandroff de $\mathbb{N}$.
    Tout $n\in \mathbb{N}$ est un point isolé de $\alpha \mathbb{N}$.\\
    Si $\infty$ est un point isolé de $\alpha \mathbb{N}$, alors $\{n\}_{n\in \mathbb{N}}\cup\{\infty\}$ est un recouvrement ouvert de $\alpha \mathbb{N}$ sans sous-recouvrement fini. Ainsi $\infty$ est le seul point non-isolé de $\alpha \mathbb{N}$.
    Comme $|\alpha \mathbb{N}|\geq 2$, l'espace produit $(\alpha \mathbb{N})^{\mathbb{N}}$ est sans-point isolé.\\
    Montrons que $\alpha \mathbb{N}$ est totalement discontinue:\\
    Soit $A\subseteq \alpha \mathbb{N}$ tel que $|A|\geq 2$. Alors il existe $n\in A$ tel que $n\in \mathbb{N}$. Le singleton $\{n\}$ est à la fois ouvert et fermé, i.e $A$ n'est pas connexe. Donc $\alpha \mathbb{N}$ est totalement discontinue. Ainsi l'espace produit
    $(\alpha \mathbb{N})^{\mathbb{N}}$ est totalement discontinu.\\
    $\alpha \mathbb{N}$ est métrisable par le théorème de métrisabilité d'Urysohn (Tout espace normal $T_{1}$ vérifiant le deuxième axiome de dénombrabilité est métrisable). Ainsi l'espace produit
    $(\alpha \mathbb{N})^{\mathbb{N}}$ est métrisable. Donc $(\alpha \mathbb{N})^{\mathbb{N}}$ est homéomorphe à l'ensemble de Cantor $D^{\aleph_{0}}$.
    L'action continue de $S_{\infty}$ sur $\mathbb{N}$ se prolonge en une action continue de $S_{\infty}$ sur $\alpha \mathbb{N}$. Ainsi, $S_{{\infty}}$ opère continûment sur $(\alpha \mathbb{N})^{\mathbb{N}}$ par l'action produit.
    Cette dernière action continue entraine le plongement $S_{\infty}\hookrightarrow Homeo ((\alpha \mathbb{N})^{\mathbb{N}})=Homeo(D^{\aleph_{0}})$ de $S_{\infty}$ dans $Homeo(D^{\aleph_{0}})$.
  \item [$2)\Longrightarrow 1)$] $Homeo(D^{\aleph_{0}})$ étant non archimédien, alors tout sous-groupe $G^{\prime}$ de $S_{\infty}$ est non archimédien. Si $G$ est isomorphe à $G^{\prime}$, alors $G$ est également non archimédien.
 \end{enumerate}
\end{prof}\\
Une question naturelle serait de savoir si les deux groupes précédents sont isomorphes en tant que groupes topologiques?
Nous apporterons une réponse par la négative à cette question dans la suite du paragraphe:\\

Rappelons le résultat suivant:
\begin{theorem}(\cite{bekka})
Soit $(G_{i})_{i\in I}$ une famille ordonnée de sous-groupes fermés d'un groupe topologique $G$ tel que
$\underset{i\in I}{\overset{}{\bigcup}}G_{i}$ soit dense dans $G$. Si $G_{i}$ est moyennable pour tout $i\in I$, alors $G$ est moyennable.
\end{theorem}
\begin{corollary}
Le groupe $S_{\infty}$ est moyennable
\end{corollary}
\begin{prof}
Soit $n\in \mathbb{N}$. Notons $S_{n}$ le groupe symétrique d'ordre $n$. $(S_{n})_{n\in \mathbb{N}}$ s'identifie à une famille de sous-groupes fermés de $S_{\infty}$ croissante par l'inclusion. La croissance ici est comprise dans le sens suivant: Tout élément $f:\{1,2,...,n\}\longrightarrow \{1,2,...,n\}$ de $S_{n}$ peut s'identifier à un élément

\[\begin{array}{lll}
\widetilde{f}:& \mathbb{N}^{\star} &\longrightarrow \mathbb{N}^{\star}\\
&k &\longmapsto \widetilde{f}(k)=\left\{
    \begin{array}{ccccc}
      f(k) &\text{si}\,\, x\in \{1,2,...,n\}&  &  \\
       k & \text{sinon} & &  \\
    \end{array}\right.\\
\end{array}\]

de $S_{\infty}$. Chaque groupe $S_{n}$ étant fini, est moyennable.\\ Montrons pour terminer que
$S_{\infty}=\overline{\underset{n\in \mathbb{N}}{\overset{}{\bigcup}}S_{n}}$.\\ Posons $S=\underset{n\in \mathbb{N}}{\overset{}{\bigcup}}S_{n}$ et soit $f:\mathbb{N}^{\star}\longrightarrow \mathbb{N}^{\star}$ une bijection. Pour tout $n\in \mathbb{N}^{\star}$, posons: $$M_{n}=n+f(1)+f(2)+...+f(n).$$ Considérons \\$A=\{i\in \mathbb{N}:\,\,\,n< i\leq M_{n}\},\,\,B=\{1,2,...,M_{n}\}$ et $C=\{f(1),f(2),...,f(n)\}$. Posons $D=B\setminus C$. Il est clair que $|A|=|D|$. Ainsi il existe une bijection $g_{n}:A\longrightarrow D$ . Considérons la suite de fonctions $(f_{n})_{n\in \mathbb{N}}$ de $\mathbb{N}$ dans $\mathbb{N}$ définie par:
\[ f_{n}(k)= \left\{ \begin{array}{ll}
f(k) & \textrm{si $ k \leq n$}\\
g_{n}(k)& \textrm{si $n< k \leq M_{n}$}\\
k&\textrm{si $k> M_{n}$}\end{array}\right.\]
 La suite $(f_{n})$ est une suite d'éléments de $S$ et elle converge simplement vers $f$. Ainsi $S_{\infty}$ est moyennable.
\end{prof}\\
Le lemme suivant est crucial pour établir la non moyennabilité du groupe $Homeo(D^{\aleph_{0}})$.
\begin{lemma}\label{lemtrans}Le groupe
$Homeo(D^{\aleph_{0}})$ opère d'une manière transitive à la fois sur les points de $D^{\aleph_{0}}$ et sur les parties propres de $D^{\aleph_{0}}$ à la fois ouvertes et fermées.
\end{lemma}
\begin{theorem}
Le groupe $Homeo(D^{\aleph_{0}})$ n'est pas moyennable.
\end{theorem}
\begin{prof}
Supposons le contraire et soit $\mu$ une mesure de probabilité invariante pour l'action continue de $Homeo(D^{\aleph_{0}})$ sur $D^{\aleph_{0}}$. Puisque $D^{\aleph_{0}}$ est totalement discontinue, il existe une partition $\mathcal{C}=(C_{n})_{n\in\mathbb{N}}$ de $D^{\aleph_{0}}$ constituée d'ensembles à la fois ouverts et fermés. Soit $C_{i_{0}}\in \mathcal{C}$ tel que $0<\mu(C_{i_{0}})<\frac{1}{3}$. D'après le lemme \ref{lemtrans}, il existe
$f\in Homeo(D^{\aleph_{0}})$ tel que $f(C_{i_{0}})=D^{\aleph_{0}}\setminus C_{i_{0}}$ et  $f(D^{\aleph_{0}}\setminus C_{i_{0}})=C_{i_{0}}$. Donc $\mu(D^{\aleph_{0}}\setminus C_{i_{0}})=\mu(C_{i_{0}})$ par l'invariance de $\mu$. Mais
$$\mu(D^{\aleph_{0}}\setminus C_{i_{0}})=1-\mu(C_{i_{0}})>\frac{2}{3}> \mu(C_{i_{0}})$$  Ce qui est absurde. D'où $\mu$ n'est pas préservé par $Homeo(D^{\aleph_{0}})$.
\end{prof}
\begin{corollary}
Les groupes topologiques $Homeo(D^{\aleph_{0}})$ et $S_{\infty}$ ne sont pas isomorphes.
\end{corollary}

\subsection{Décomposition du compactifié équivariant de Samuel en limite inverse}
Commençons par rappeler le lemme suivant:
\begin{lemma}(\cite{Bou} chap. $10$, prop. $8$)\label{theopropre}
Soient $X$ un espace compact, $\mathcal{R}$ une relation d'équivalence sur $X$ et $C$ son graphe dans $X\times X$. $C$ est fermée dans $X\times X$ si et seulement si l'espace quotient $X/\mathcal{R}$ est séparé.
\end{lemma}
A présent, soient $G$ un groupe topologique, $Z$ un espace topologique, et $f\colon G\to Z$ une application uniformément continue à droite de $G$ telle que $X=f(G)$ soit compact, alors $f$ se prolonge en une application continue $f : S(G)\longrightarrow X$, encore notée $f$.\\
Définissons sur $S(G)$ la relation d'équivalence:
$$x\mathcal{R}y \,\,\text{si} \,\,f(gx) = f(gy)\,\,\text{pour tout}\,\, g\in G.$$
Notons $C$ le graphe de $\mathcal{R}$. \\
Soit $(x,y)\in (S(G)\times S(G)\setminus C)$, il existe $g\in G$ tel que $f(gx)\neq f(gy)$. $X$ étant séparé, il existe deux ouverts disjoints $U$ et $V$ de $X$ tel que  $f(gx)\in U$ et $f(gy)\in V$. En posoant $W=g^{-1}(f^{-1}(U))\times g^{-1}(f^{-1}(V))$, on a $W\subset (S(G)\times S(G)\setminus C)$. Donc $C$ est fermé dans $S(G)\times S(G)$ et l'espace quotient $X_{f}=S(G)/ \mathcal{R}$ est séparé par le lemme \ref{theopropre}. Puisque la surjection canonique $\pi_f:S(G)\longrightarrow X_{f}$ est continue, l'espace quotient $X_{f}$ est compact.
De plus, l'application
\[\begin{array}{lll}
\bar f:& X_{f} &\longrightarrow X\\
&[x] &\longmapsto f(x)\\
\end{array}\]
 est continue (par définition de la topologie quotient) et rend commutatif le diagramme suivant:
   \begin{center}
    \begin{tabular}{c}
     \xymatrix@M=0pt@R=27pt@C=27pt{
    S(G) \ar@{->}[rr]^{f}\ar@{->}[dd]_{\pi_{f}}&& X  \\
    &&\\
    X_{f}\ar@{->}[uurr]_{\bar f}\\ }
    \end{tabular}
    \end{center}
Soit $F$ une famille d'applications uniformément continues à droite de $G$ dans un espace compact $Y$ (vues comme applications $S(G)\longrightarrow Y$ ). Posons $X = Y^F$ et soit $f:S(G)\longrightarrow X$ le produit diagonal de la famille $F$. Notons $X_F = X_f$ dans ce cas-là, et $\pi_F\colon S(G)\to X_F$ la surjection canonique. Le lemme suivant est immédiat:
\begin{lemma}\label{lemm2.8}
Si $F$ sépare les points de $S(G)$, c'est à dire, pour tous $x,y\in S(G)$ avec $x\neq y$, il existe $h\in F$ tel que $h(x)\neq h(y)$, alors $\pi_F$ est un homéomorphisme de $S(G)$ sur $X_F$.
\end{lemma}
\begin{prof}
$\pi_F$ est surjective et continue par définition. Rappelons que $\pi_F$ est la surjection canonique associée à la relation d'équivalence:\\
 $x \mathcal{R}_{F} y\,\Longleftrightarrow$ pour tout $g\in G$ et pour tout $h\in F$, on a $h(gx)=h(gy)$.  Comme $F$ sépare les points de $S(G)$, pour tout $x,y\in S(G)$ avec $x\neq y$, il existe $h\in F$ tel que $h(x)\neq h(y)$ et alors $\pi_{F}(x)\neq \pi_F(y)$. Donc $\pi_F$ est injective. $\pi_F$ est donc un homéomorphisme.
\end{prof}\\
Notons que pour toute famille d'applications $F$, $G$ opère continûment sur l'espace quotient $X_{F}$ par l'action quotient: $$g[x]_{\mathcal{R}_{F}}=[gx]_{\mathcal{R}_{F}}.$$
\begin{lemma}\label{pif}
Si $F\subseteq F^{\prime}$, alors il existe une surjection continue et $G$-équivariante canonique $\pi^{F^{\prime}}_F$ de $X_{F^{\prime}}$ sur $X_F$.
\end{lemma}
\begin{prof}
Soit $x\in S(G)$. Pour tout $F\subseteq F^{\prime}$, définissons l'application
\[\begin{array}{lll}
\pi^{F^{\prime}}_{F}:& X_{F^{\prime}} &\longrightarrow X_{F}\\
&[x]_{\mathcal{R}_{F^{\prime}}} &\longmapsto [x]_{\mathcal{R}_{F}}\\
\end{array}\]

Montrons que $\pi^{F^{\prime}}_{F}$ est bien définie. Soit $x,y \in S(G)$ tel que $[x]_{\mathcal{R}_{F^{\prime}}}=[y]_{\mathcal{R}_{F^{\prime}}}$, alors pour tout $g\in G$ et pour tout $h\in F^{\prime}$, on a $h(gx)=h(gy)$ et en particulier, pour tout $g\in G$ et pour tout $h\in F$, on a $h(gx)=h(gy)$. Donc $[x]_{\mathcal{R}_{F}}=[y]_{\mathcal{R}_{F}}$. De même, pour tout $F\subseteq F^{\prime}$, on a: $\pi^{F^{\prime}}_{F}\circ \pi_{F^{\prime}}=\pi_{F}$. Donc $\pi^{F^{\prime}}_{F}$ est continue.  $\pi^{F^{\prime}}_{F}$  est équivariante par rapport à l'action quotient sur $X_{F^{\prime}}$ et sur $X_{F}$. En effet, pour tout $x\in S(G)$, on a:
$$g\pi^{F^{\prime}}_{F}([x]_{\mathcal{R}_{F^{\prime}}})=g[x]_{\mathcal{R}_{F}}=[gx]_{\mathcal{R}_{F}}=
\pi^{F^{\prime}}_{F}([gx]_{\mathcal{R}_{F^{\prime}}})=\pi^{F^{\prime}}_{F}(g[x]_{\mathcal{R}_{F^{\prime}}}).$$
\end{prof}

\begin{lemma}Sous les hypothèses du lemme \ref{pif}, et si $G$ est dénombrable, $X$ compact métrisable, et $F=(f_i)_{i\in\mathbb{N}}$ est dénombrable, alors $X_F$ est compact et métrisable.
\end{lemma}
\begin{prof}
Pour tout $g\in G$, notons $g:S(G)\longrightarrow S(G)$ l'action de $G$ sur $S(G)$. Pour tout $f_i\in F$, on a: $f_i:S(G)\longrightarrow X$. Considérons le produit diagonal $\Delta_{i}f_{i}:S(G)\longrightarrow X^{|F|}$ et posons $f^{\prime}=\Delta_{i}f_{i}\circ g:S(G)\longrightarrow X^{|F|}$. Comme $X$ est compact et métrisable, $X^{|F|}$ est aussi compact et métrisable car $F$ est dénombrable. Notons encore $f^{\prime}$ l'application $f^{\prime}:S(G)\longrightarrow f^{\prime}(S(G))$. Par le même raisonement que précédement, il existe une application $\overline{f^{\prime}}:X_{F}\longrightarrow X=f^{\prime}(S(G))$ continue telle que $\overline{f^{\prime}}=\bar f\circ \pi_{F}$. Comme $F$ sépare les points de $S(G),\,\,\overline{f^{\prime}}$ est injective. $X_F$ étant compact, $\overline{f^{\prime}}$ est un homéomorphisme de $X_F$ sur $f^{\prime}(S(G))$. Ainsi $X_F$ est compact et métrisable.
\end{prof}
 \begin{lemma}     Soit $\Phi$ une collection de familles de fonctions de $S(G)$ dans $X$. Supposons que $\Phi$ est dirigée par l'inclusion, c'est à dire, quels que soient $F,F^{\prime}\in\Phi$, il existe $F^{\prime\prime}\in\Phi$ tel que $F\subseteq F^{\prime\prime}, F^{\prime}\subseteq F^{\prime\prime}$. Alors le système $(X_F,\pi^{F^{\prime}}_F,\Phi)$ forme un système projectif de $G$-espaces.
 \end{lemma}
 \begin{prof}
Pour tout $F\in \Phi,\,\,G$ opère continûment sur $X_{F}$ par l'action quotient. Pour tout $F\subseteq F^{\prime}$, les applications $\pi^{F^{\prime}}_F$ sont équivariantes par le lemme \ref{pif}. Soit $F,F^{\prime}, F^{\prime\prime} \in\Phi$ tel que $F\subseteq F^{\prime}\subseteq F^{\prime\prime}$ et soit $x\in S(G)$, on a: $$(\pi^{F^{\prime}}_{F}\circ\pi^{F^{\prime \prime}}_{F^{\prime}})([x]_{\mathcal{R}_{F^{\prime \prime}}})=\pi^{F^{\prime}}_{F}([x]_{\mathcal{R}_{F^{\prime}}})=[x]_{\mathcal{R}_{F}}=\pi^{F^{\prime \prime}}_{F}([x]_{\mathcal{R}_{F^{\prime \prime}}}).$$ Donc $$\pi^{F^{\prime}}_{F}\circ\pi^{F^{\prime \prime}}_{F^{\prime}}=\pi^{F^{\prime \prime}}_{F}.$$
\end{prof}
 \begin{remark}\label{rem1}
Supposons que la réunion $\cup \{F\colon F\in\Phi\}$ sépare les points de $S(G)$. Alors la limite projective du système $(X_F,\pi^{F^{\prime}}_F)$ est isomorphe, en tant que $G$-espace compact à $S(G)$. En effet, considérons l'application
\[\begin{array}{lll}
\Psi:& S(G) &\longrightarrow \underset{\longleftarrow}{\overset{}{\lim}}\,X_{F}\\
&x &\longmapsto (\pi_F(x))_{F\in \Phi}\\
\end{array}\]

Soit $x,y\in S(G)$ tel que $x\neq y$. puisque $\cup \{F\colon F\in\Phi\}$ sépare les points de $S(G)$, il existe $F\in\Phi$ tel que $\pi_F(x)\neq \pi_F(y)$. Donc $\Psi(x)\neq\Psi(y)$ et $\Psi$ est injective. $\Psi$ est surjective par définition. Soit $x\in S(G)$, on a:
$$g\Psi(x)=g(\pi_F(x))_F=(g\pi_F(x))_F=(\pi_F(gx))_F=\Psi(gx).$$ Donc $\Psi$ est équivariante.
\end{remark}

\begin{lemma}
Si $X$ est un compact de dimension zéro (au sens de Lebesgue), alors $X_{F}$ est de dimension zéro.
\end{lemma}
\begin{prof}
Soit $[x],[y]\in X_{F}$ tel que $[x]\neq[y]$, alors $x\neq y$. Puisque $F$ sépare les points de $S(G)$, il existe $f_{i}\in F$ tel que $f_{i}(x)\neq f_{i}(y)$. Comme $X$ est totalement discontinu, il existe une décomposition $U\cup V$ de $X$ telle que $f_{i}(x)\in U$ et $f_{i}(y)\in V$. Posons  $\dot{U}=\overline{f_{i}}^{-1}(U)$ et $\dot{V}= \overline{f_{i}}^{-1}(V)$. Clairement, $\dot{U}\cup \dot{V}$ est une décomposition de $X_F$ et $[x]\in \dot{U},\,\,[y]\in \dot{V}$. Ainsi $X_{F}$ est totalement discontinu, donc de dimension $0$.
\end{prof}
\begin{remark}
Si $X_F$ contient un point isolé $x_0$, alors le sous-groupe ouvert \\
$H=St_{x_0}=\{g\in G:\,\,gx_0=x_0\}$ est tel que pour tout $h\in F$, on a:$h|_{gH}$ est constante.\\
 Pour garantir la non-existence des points isolés, nous allons supposer que la famille $F$ vérifie la condition $(\star)$ suivante:\\
  Pour tout voisinage $V$ de $e$, il existe $h\in F$ et $x\in V$ tel que $h(x)\neq h(e)$.
\end{remark}
Dans ces conditions, on a le lemme:
\begin{lemma}
$G$ étant un groupe topologique, $X$ un espace compact, et $X_{F}$ définie comme précédemment avec $F$ vérifiant la condition $(\star)$, alors ou bien $X_{F}$ est fini, ou bien $X_{F}$ ne contient aucun point isolé.
\end{lemma}
\begin{prof}
Supposons $X_{F}$ infini et soit $x_{0}$ un point isolé de $X_{F}$. Le sous-groupe ouvert $H=St_{x_0}=\{g\in G:\,\,gx_0=x_0\}$ est donc un voisinage de $e$ tel que pour tout $h\in F,\,\,\,h|_{gH}$ est constante. Ce qui est absurde puisque $F$ vérifie la condition $(\star)$. En effet, vérifié la propriété $(\star)$ revient à dire qu'il existe $h\in F$ non-constante sur $H$ en particulier.
\end{prof}
\begin{corollary}\label{corcantor}
Si $G$ est infini dénombrable, $X$ compact de dimension zéro et $F$ dénombrable et vérifie $(\star)$, alors $X_{F}$ est homéomorphe à l'espace de Cantor.
\end{corollary}
On déduit de tout ce qui précède le théorème fondamental suivant:
\begin{theorem}\label{coroinverse}
Si $G$ est non archimédien et polonais, alors $S(G)$ se développe en limite projective d'un système de $G$-espaces compacts métrisables de dimension zéro. Si de plus, $G$ est infini, alors il existe un développement constitué de $G$-espaces homéomorphes à l'espace de Cantor.
\end{theorem}
\begin{theorem}
 L'ensemble de Cantor $D^{\aleph_{0}}$ est un espace test pour la moyennabilité des groupes polonais non archimédien. Autrement dit, un groupe polonais non archimédien $G$ est moyennable si et seulement si toute action continue de $G$ sur $D^{\aleph_{0}}$ possède une mesure de probabilité invariante.
\end{theorem}
\begin{prof}
La nécessité est évidente. Montrons la suffisance. $G$ étant un groupe polonais non archimédien, il existe par le théorème \ref{coroinverse}, un système projectif de $G$-espaces $(X_{\alpha},\pi_{\alpha \beta}, I)$ tel que pour tout $\alpha \in I$, $X_{\alpha}$ est homéomorphe à l'ensemble de Cantor $D^{\aleph_{0}}$ et $S(G)=\underset{\longleftarrow}{\overset{}{\lim}}\,X_{\alpha}$. Par hypothèse, il existe sur chaque $G$-espace $X_{\alpha}$ une mesure de probabilité invariante $\mu_{\alpha}$. Par le lemme \ref{lemmesurinverse}, il existe une mesure de probabilité borélienne invariante $\mu$ sur $S(G)$. Donc $G$ est moyennable.
\end{prof}

\section{Espace test pour la moyennabilité extrême}
Dans cette section, nous utilisons la décomposition du compactifié équivariant de Samuel en limite inverse d'espaces de Cantor obtenue dans la section précédente pour établir que l'espace de Cantor est un espace test pour la moyennabilité extrême des groupes polonais non archimédiens.
\subsection{Généralités sur la moyennabilité extrême}
\begin{definition}
Un groupe topologique $G$ est dit \index{groupe!extrêmement moyennable}{\em extrêmement moyennable} si toute action continue de $G$ sur un espace compact $K$ possède un point fixe.
\end{definition}
\begin{remark}
\begin{enumerate}
 \item Il est clair que tout groupe $G$ extrêmement moyennable est  moyennable. L'extrême moyennabilité est donc une propriété plus forte que la moyennabilité d'où la terminologie extrêmement moyennable introduite par Granirer \cite{granirer}.
  \item Soit $G$ un groupe topologique. On note $M(G)$ son espace compact universel. $G$ est extrêmement moyennable si et seulement si $M(G)$ est un singleton.\\
En effet, Si $M(G)=\{x\}$, soit $X$ un flot compact de $G$, alors il existe un morphisme de flots $\pi:M(G)\longrightarrow X$ et $\pi(x)$ est un point fixe pour l'action continue de $G$ sur $X$.\\
Si $G$ est extrêmement moyennable, comme $G$ opère sur l'espace compact $M(G)$, il existe un point fixe $x_{0}$ pour l'action de $G$ sur $M(G)$. Autrement dit, on a: $G.x_{0}=\{x_{0}\}$. Puisque $M(G)$ est minimal, $\overline{G.x_{0}}=M(G)$ et $M(G)=\{x_{0}\}$.
  \item En général, $G$ est extrêmement moyennable si et seulement si tout flot minimal de $G$ est un singleton.
\end{enumerate}
\end{remark}
Le résultat suivant est immédiat.
\begin{proposition}\label{pr3}
Un groupe topologique $G$ est extrêmement moyennable si et seulement si l'action canonique de $G$ sur son compactifié équivariant de Samuel $S(G)$ possède un point fixe.
\end{proposition}
La technique pour établir le résultat suivant est identique à celle utilisée dans la preuve du thèoreme \ref{moyexample}.
\begin{proposition}
Soit $G$ un groupe topologique.
\begin{enumerate}
  \item Si $H$ est un groupe topologique tel qu'il existe un homomorphisme continu surjectif de $G$ sur $H$ et si $G$ est extrêmement moyennable, alors $H$ est extrêmement moyennable.
  \item S'il existe un sous-espace dense $A$ de $G$ tel que tout sous-ensemble fini $I$ de $A$ est contenu dans un sous-groupe extrêmement moyennable, alors $G$ est extrêmement moyennable.
  \item Si $H$ est un sous-groupe normal extrêmement moyennable de $G$ et le groupe quotient $G/H$ est extrêmement moyennable, alors $G$ est extrêmement moyennable.
\end{enumerate}
\end{proposition}

Les groupes de Lévy\index{groupe!de Lévy} constituent une classe importante de groupes extrêmement moyennables.
\subsection{Groupes de Lévy et moyennabilité extrême}
\begin{definition}
On appelle {\em mm-espace}\index{mm-espace} la donnée d'un triplet $(X,\,d,\,\mu)$ où:
\begin{enumerate}
\item[1)] $X$ est un ensemble non vide,
\item[2)] $d$ une distance sur $X$,
\item[3)]  $\mu$ une mesure de probabilité sur $X$.
 \end{enumerate}
\end{definition}
\begin{definition}
\begin{enumerate}
  \item Soit $(X,\,d,\,\mu)$ un mm-espace. On appelle {\em fonction de concentration}\index{fonction!de concentration} de
$X$ la fonction notée $\alpha_{X}$ définie sur $\mathbb{R}_{+}$ par:
\[ \alpha_{X}(\varepsilon)= \left\{ \begin{array}{ll}
\dfrac{1}{2} & \textrm{si $ \varepsilon =0$}\\
& \\
1- inf\{\mu(B_{\varepsilon}):\,\,\,\, B\subseteq X ;\,\,\,\,\, \mu (B)\geq \dfrac{1}{2}\} & \textrm{si $\varepsilon> 0$}\end{array}\right.\]
Où $B_{\varepsilon}$ désigne le voisinage d'ordre $\varepsilon$ de $B$.
  \item Une famille $(X_{n},d_{n},\mu_{n})_{n\in \mathbb{N}}$ de mm- espace est dite de Lévy si la suite de fonctions
$(\alpha_{X_{n}})_{n\in \mathbb{N}}$ converge simplement vers $0$ sur $]0,\,+\infty[$.
\end{enumerate}
\end{definition}
Quelques propriétés des mm-espaces et des familles de Lévy sont condensées dans le résultat suivant:
\begin{proposition}\cite{ledoux}\label{levyfam}
\begin{enumerate}
  \item Pour tout mm-espace $(X,d,\mu)$, on a: $\underset{\varepsilon \rightarrow +\infty}{\overset{}{\lim}}\,\alpha_{X}(\varepsilon)=0$.

  \item Une famille $(X_{n},d_{n},\mu_{n})_{n\in \mathbb{N}}$ de mm- espace est de Lévy si et seulement si chaque fois qu'un sous-ensemble borelien $A_{n}$ de $X_{n}$ vérifie: $\mathrm{\liminf}\,\mu_{n}(A_{n})> 0$ on a pour tout $\varepsilon> 0\,\,\,\underset{n \rightarrow +\infty}{\overset{}{\lim}}\,\mu_{n}((A_{n})_{\varepsilon})=1$.

\end{enumerate}
\end{proposition}

\begin{example}
\begin{enumerate}
  \item Sur le  groupe symétrique $S_{n}$ on défini la distance normalisée de Hamming
$d(\sigma,\tau)=\dfrac{1}{n}|\{i:\,\,\,\sigma(i)\neq\tau(i)\}|$
et la mesure uniforme $\mu(A)=\dfrac{|A|}{n}$. La famille $(S_{n},\,d,\,\mu)$ est une famille de Lévy. (Maurey \cite{mau}).
  \item Soit $E_{n}=\{0,1\}^{n}$ le cube de Hamming. On défini sur $E_{n}$ la mesure normalisée
$P_{n}(A)=|A|.2^{-n}$ et la distance normaliée de Hamming:
$d(x,y)=\dfrac{1}{n}|\{i:\,x_{i}\neq y_{i}\}|$. la famille $(E_{n},\,d,\,P_{n})$ est une famille de Lévy (\cite{grommil1}).
\end{enumerate}
\end{example}

\begin{definition}
Soit $G$ un groupe métrisable opérant continûment par isomorphismes sur un espace métrique $(X,d)$. Le $G$-espace $X$ est dit de Lévy (Gromov et Milman\cite{grommil}) s'il existe une suite $(G_{\alpha})$ de sous-groupes compacts de $G$ ordonné par l'inclusion et une suite $(\mu_\alpha)$ de mesures de probabilités sur $(X,d)$ telle que:
\begin{enumerate}
  \item $\underset{\alpha}{\overset{}{\bigcup}}G_{\alpha}$ est dense dans $G$
  \item $\mu_\alpha$ est $G_\alpha$-invariant pour tout $\alpha$
  \item $(X,d,\mu_\alpha)$ est une famille de Lévy
\end{enumerate}
\end{definition}
Dans le cas particulier où $X=G$ est muni d'une distance invariante à droite et de l'action à gauche de $G$ sur lui-même, le groupe $G$ est dit de Lévy.

\begin{theorem}(Gromov et Milman\cite{grommil})\\
Tout groupe de Lévy est extrêmement moyennable.
\end{theorem}
La référence \cite{kei} contient des  progrès récents sur la relation entre les familles de Lévy et l'existence de points fixes.\\
À présent, les exemples de groupes de Lévy sont nombreux et comprennent:
\begin{enumerate}
\item Le groupe unitaire $\mathcal{U}(\ell^{2})$, muni de la topologie forte (Gromov et Milman \cite{grommil}).
\item Le groupe $Aut\,(X,\mu)$ des automorphismes mesurables préservant la mesure $\mu$ d'un espace borelien $(X,\mu)$ muni de la topologie faible (Giordano et Pestov \cite{giopes1}). Rappelons que la topologie faible sur $Aut\,(X,\mu)$ est la topologie la moins fine rendant continue toutes les applications
    \[\begin{array}{lll}
\Upsilon_{A}:& Aut\,(X,\mu) &\longrightarrow X\\
&T &\longmapsto \Upsilon_{A}(T)=T(A)\\
\end{array}\] où $A$ est un borélien de $X$.
 \item Le groupe $Iso(\mathbb{U})$ des isométries de l'espace d'Urysohn $\mathbb{U}$ muni de la topologie compact-ouvert (Pestov \cite{vp4}).
 \end{enumerate}
La classe des groupes extrêmement moyennables comprend aussi le groupe $Aut\,(\mathbb{Q},\leq)$ des bijections de $\mathbb{Q}$ sur lui-même qui préservent l'ordre muni de la topologie de la convergence simple (Pestov \cite{vp2}).\\
Il existe également des exemples de groupes non-extrêmement moyennables:
\begin{enumerate}
 \item Les groupes localement compacts (Veech \cite{vp1}) .
\item Le groupe symétrique infini $S_{\infty}$ muni de sa topologie polonaise (Pestov \cite{vp2}).
\end{enumerate}
\subsection{Espace test pour la moyennabilité extrême}
\begin{theorem}\label{main3}
Un groupe polonais non archimédien $G$ est extrêmement moyennable si et seulement si toute action continue de
$G$ sur l'ensemble de Cantor $D^{\aleph_{0}}$ possède un point fixe.
\end{theorem}
\begin{prof}
  La nécessité est évidente. Montrons la suffisance.\\
  Montrons que l'action canonique de $G$ sur $S(G)$ possède un point fixe. Comme
  $G$ est polonais et non-archimédien, il existe par le théorème \ref{coroinverse}, un système projectif de $G$-espaces $(X_{\alpha},\pi_{\alpha \beta}, I)$ avec $X_{\alpha}\cong D^{\aleph_{0}}$ pour tout $\alpha \in I$ tel que $S(G)=\underset{\longleftarrow}{\overset{}{\lim}}\,X_{\alpha}$. Par hypothèse, il existe sur chaque $G$-espace $X_{\alpha}$ un point fixe $x_{\alpha}$. Notons $\pi_{\alpha}$ la restriction de la projection $pr_{\alpha}$ à $S(G)=\underset{\longleftarrow}{\overset{}{\lim}}\,X_{\alpha}$ et posons $M_{\alpha}=\pi_{\alpha}^{-1}(x_{\alpha})$. Les applications $\pi_{\alpha}$ étant surjectives, on a: $M_{\alpha}\neq \emptyset$ pour tout $\alpha$. La famille $(M_{\alpha})_{\alpha \in I}$ est centrée car pour tout $i=1,2,..,n\,\,\,\,\,x=(x_{\alpha_{1}},...,x_{\alpha_{n}})\in \underset{i=1}{\overset{n}{\cap}}\,M_{\alpha_{i}}$.
   $S(G)$ étant compact, $\underset{\alpha \in I}{\overset{}{\bigcap}}\,M_{\alpha}\neq\emptyset$. Si $x\in \underset{\alpha \in I}{\overset{}{\bigcap}}\,M_{\alpha}\neq\emptyset$, alors $x$ est un point fixe pour l'action continue de $G$ sur $S(G)$.
\end{prof}

La question suivante reste une question ouverte. Nous y apportons une réponse partielle.
\begin{question}
Existe-t-il un espace test compact et métrisable pour les groupes polonais extrêmement moyennables?
\end{question}
Le théorème du point fixe de Schauder (\cite{schauder}) affirme que toute fonction continue de $I^{\aleph_{0}}$ dans $I^{\aleph_{0}}$ possède un point fixe. En particulier, toute action continue du groupe discret $\mathbb{Z}$ sur $I^{\aleph_{0}}$ par les homéomorphismes possède un point fixe. Ceci permet de conclure que le cube de Hilbert $I^{\aleph_{0}}$ ne peut pas être un espace test pour la moyennabilité extrême des groupes polonais.\\En effet, le théorème de Ellis \cite{el} affirme que tout groupe discret agit librement sur un espace compact et par conséquent, n'est pas extrêmement moyennable.\\
On peut néanmoins observer qu'il existe un espace test compact séparable non nécessairement métrisable pour les groupes polonais extrêmement moyennables.\\
En effet, notons $\mathcal{P}_{0}$ l'ensemble de tous les groupes polonais non-extrêmement moyennables deux à deux non-isomorphes et choisissons pour tout $G\in \mathcal{P}_{0}$ un $G$-espace compact et métrisable $X_{G}$ sans point fixe. L'espace $X=\underset{G\in \mathcal{P}_{0}}{\overset{}{\prod}}X_{G}$ est un espace test séparable compact (non néccessairement métrisable) pour la moyennabilité extrême des groupes polonais .\\
Il est clair que $G$ opère continûment sur  $X$ sans point fixe par l'action produit.\\
Il nous suffit pour conclure  par le célèbre théorème de Hewitt \cite{hewitt} et Pondiczery \cite{pondi} de montrer que
 $|\mathcal{P}_{0}|\leq 2^{\aleph_0}$. Notons $F_{\infty}$ le groupe libre avec un nombre infini dénombrable de générateurs. Notons $\mathcal{P}$ l'ensemble des groupes  polonais et $\mathcal{D}$ l'ensemble de toutes les pseudo-métriques sur $F_{\infty}$. Il est clair que $|\mathcal{D}|\leq |\mathbb{R}^{\mathbb{Z}}|=2^{\aleph_0}$. Nous allons montrer que $|\mathcal{P}| \leq |\mathcal{D}|$.\\
Soit $d$ une pseudo-métrique sur $F_{\infty}$ invariante à gauche. $H_{d}=\{x\in F_{\infty}, d(x,e)=0\}$ est un sous-groupe de $F_{\infty}$. La distance défini sur $F_{\infty}/H_{d}$ par $\widehat{d}(xH_{d},yH_{d})=d(x,y)$ est invariante par translation à gauche.
Notons $G_{d}$ le completé de l'espace métrique $(F_{\infty}/H_{d},\widehat{d})$. Si  $H_{d}$ est un sous-groupe normal de $F_{\infty}$, le groupe topologique $G_{d}$ est un groupe polonais et chaque groupe polonais est de la forme $G_{d}$. Notons $\mathcal{D}_{N}$ le sous-ensemble de $\mathcal{D}$ constitué des pseudo-métriques $d$ telles que $H_{d}$ soit normal. Donc l'application
\[\begin{array}{lll}
& \mathcal{D}_{N} &\longrightarrow \mathcal{P}\\
&d &\longmapsto G_{d}\\
\end{array}\]

est surjective. Ainsi $|\mathcal{P}|\leq |\mathcal{D}_{N}|\leq |\mathcal{D}|$.\\

\begin{definition}
Un groupe topologique $G$ est dit \index{groupe!monothétique} {\em monothétique}, s'il existe un sous-groupe $H$ de $G$ qui est à la fois cyclique et dense
\end{definition}
\begin{example}(Glasner \cite{gla3})\\
Soit $\Gamma$ un groupe localement compact abélien(Par exemple $\Gamma=\mathbb{R}$). Nous notons $\widetilde{\Gamma}$ le groupe de tous les caractères continus de $\Gamma$ i.e le groupe des homomorphismes continus de $\Gamma$ dans le cercle unité $\mathbb{T}$. Soit $E$ un espace localement compact. Nous notons par $S(E)$ le groupe des fonctions $f:E\longrightarrow \mathbb{C}$ continues vérifiant $|f|=1$. Si $E\subset \Gamma$, nous notons par $SU(E)$ le sous-groupe de $S(E)$ constitué des fonctions qui sont uniformément continues (par rapport à la restriction de la structure uniforme canonique du groupe $\Gamma$ sur $E$.)

\begin{definition}
Un sous-ensemble fermé $E\subset \Gamma$ est dit de \index{Kronecker}Kronecker si pour tout $\varepsilon>0$, et pour tout $f\in SU(E)$, il existe $\chi \in \widetilde{\Gamma} $ tel que $\underset{x\in E}{\overset{}{\sup}}|f(x)-\chi(x)|\leq \varepsilon.$
\end{definition}
Soit $\Omega$ un sous-ensemble de Kronecker du cercle unité $\mathbb{T}$ et soit $(\Omega,\mathcal{B},\mu)$ un espace mesuré. Notons $G=\{f:\Omega \longrightarrow \mathbb{C},\,\,f\,\text{mesurable},\,\,|f|=1\}$. Munissons $G$ de la multiplication point par point et de la distance $d(f,g)=\int|f-g|d\mu$. Alors $G$ est un groupe polonais monothétique. En effet, par définition de $\Omega$, l'ensembles des restrictions à $\Omega$ des caratères continus sur $\mathbb{T}$ est uniformément dense dans l'ensemble des fonctions continues sur $\Omega$ à valeurs complexes et de module $1$. Donc le sous-groupe $G_{0}=\{\chi:\Omega \longrightarrow \mathbb{C},\,\,\chi\,\text{est un caractère continu de}\,\,\mathbb{T}\}$ est dense dans G. Puisque le groupe des caractères continus de $\mathbb{T}$ est isomorphe à $\mathbb{Z}$, on conlut que $G$ est monothétique.
\end{example}
\begin{definition}
 Un groupe topologique $G$ est dit \index{groupe!solénoïde}{\em solénoïde}, s'il existe un morphisme continu $f$ de $\mathbb{R}$ dans $G$ dont l'image est partout dense dans $G$.
\end{definition}
\begin{remark}
Il est clair que tout groupe monothétique ou solénoïde est abélien, donc moyennable.
\end{remark}
\begin{theorem}\label{actR}
Toute action continue de $\mathbb{R}$ sur $I^{\aleph_0}$ possède un point fixe.
\end{theorem}
\begin{prof}
Soit $n\in \mathbb{N}^{\star}$. Posons $$G_{n}=\frac{1}{2^{n}}\mathbb{Z}=\{\frac{k}{2^{n}}:\,\,\,\,k\in \mathbb{Z}\}.$$ Toute action continue de $\mathbb{R}$ sur $I^{\aleph_0}$ induit une action continue de $G_{n}$ sur $I^{\aleph_0}$. D'après le théorème du point fixe de Schauder, il existe un point fixe pour l'action de $G_{n}$ sur $I^{\aleph_0}$. Posons $$F_{n}=\{x_{n}\in I^{\aleph_0}:\,\,\,\,G_{n}.x_{n}=x_{n}\}.$$ Pour tout $n,m\in \mathbb{N}$ avec $n< m$, on: $G_{n}\subset G_{m}$ et $F_{m}\subseteq F_{n}$. Ainsi, $(F_{n})_{n\in \mathbb{N}^{\star}}$ est une famille centrée de $I^{\aleph_0}$. Donc $\underset{n=1}{\overset{\infty}{\bigcap}}F_{n}\neq \emptyset$.\\
 Si $x\in \underset{n=1}{\overset{\infty}{\bigcap}}F_{n}$, alors $x$ est un point fixe pour l'action de $\underset{n=1}{\overset{\infty}{\bigcup}}G_{n}$ sur $I^{\aleph_0}$. Si $x\in \mathbb{R}$, alors la suite $(2^{-n}E(x2^{-n}))_{n\in \mathbb{N}}$ est une suite de points de $\underset{n=1}{\overset{\infty}{\bigcup}}G_{n}$ qui converge vers $x$. Donc $\underset{n=1}{\overset{\infty}{\bigcup}}G_{n}$ est dense dans $\mathbb{R}.$ Ainsi, $x$ est un point fixe pour l'action de $\mathbb{R}$ sur $I^{\aleph_0}$.
\end{prof}\\
Par la même démarche, on a le résultat suivant:
\begin{theorem}
Soit $G$ un groupe monothétique. Toute action continue de $G$ sur le cube de Hilbert $I^{\aleph_0}$ possède un point fixe.
\end{theorem}
De même nous avons le résultat suivant:
\begin{theorem}
Soit $G$ un groupe solénoïde. Toute action continue de $G$ sur le cube de Hilbert $I^{\aleph_0}$ possède un point fixe.
\end{theorem}
\begin{prof}
Soit $G$ un groupe solénoïde opérant continûment sur le cube de Hilbert $I^{\aleph_0}$ et soit $f:\mathbb{R}\longrightarrow G$ un morphisme continue à image dense. $\mathbb{R}$ opère continûment sur $I^{\aleph_0}$ par l'action $(r,x)\longmapsto r.x=f(r)x$. Par le théorème \ref{actR}, il existe $\xi\in I^{\aleph_0}$ tel que $r.\xi=\xi$ pour tout $r\in \mathbb{R}$. Soit $g\in G$ et soit $(r_{n})_{n\in \mathbb{N}}$ une suite de nombres réels tel que $(f(r_{n}))_{n\in \mathbb{N}}$ converge vers $g$. D'une part, $f(r_n)\xi$ converge vers $g\xi$ et d'autre part $f(r_n)\xi=r_{n}.\xi$. Donc $g\xi=\xi$
\end{prof}\\
La question suivante reste une question ouverte.
\begin{question}
 Une action continue d'un groupe polonais moyennable sur le cube de Hilbert $I^{\aleph_{0}}$ possède-t-elle un point fixe?
   Même question pour un groupe moyennable discret.
\end{question}

\chapter{Espaces test pour la moyennabilité topologique}
\section{Généralités sur la moyennabilité topologique}
Dans toute la suite, $Z$ désignera un espace topologique discret dénombrable. Notons $\mathbb{C}^{Z}=\{u:Z\longrightarrow \mathbb{C}\}$ l'espace vectoriel réel des fonctions de $Z$ dans $\mathbb{C}$. Rappelons que $\ell^{\infty}(Z)$ est le sous-espace de $\mathbb{C}^{Z}$ formé des fonctions bornées dans $\mathbb{C}$
et que $\ell^{1}(Z)$ est le sous-espace vectoriel de $\ell^{\infty}(Z)$ défini par
 $$\ell^{1}(Z)=\{u\in \mathbb{C}^{Z}:\,\|u\|_1=\underset{z\in Z}{\overset{}{\sum}}|u_z|<\infty\}.$$ On a un plongement isométrique d'espaces vectoriels normés $$j:(\ell^{1}(Z),\|.\|_1)\longrightarrow (\ell^{\infty}(Z)^{\star},\|.\|_{\ell^{\infty}(Z)^{\star}})$$ défini par:
 $j(v)(u)=\underset{z\in Z}{\overset{}{\sum}}u_zv_z$ pour tout $v\in \ell^{1}(Z)$ et $u\in \ell^{\infty}(Z)$. L'ensemble $\mathbb{P}(Z)$ des mesures de probabilités sur $Z$ est la partie de $\ell^{1}(Z)$ définie par $$\mathbb{P}(Z)=\{v\in \ell^{1}(Z):\,v\geq 0\,\text{et}\,\,\|v\|_1=1\}.$$
En d'autres termes $\mathbb{P}(Z)$ est l'ensemble des fonctions $b:Z\longrightarrow [0,1]$ telles que: $\underset{z\in Z}{\overset{}{\sum}}b(z)=1$. En général, nous verrons $\mathbb{P}(Z)$ comme une partie de $\ell^{1}(Z)$ que nous le munirons de la topologie vague.\\
Si $G$ opère sur $Z$(dans la suite, $Z$ sera généralement $G$ et l'action sera l'action canonique par multiplication de $G$ sur $G$), alors $G$ opère aussi sur $\mathbb{P}(Z)$ par l'action $gb(z)=b(g^{-1}z)$
\begin{definition}\label{def1}
Soit $G$ un groupe discret dénombrable opérant par homéomorphismes sur un espace compact $X$. L'action de $G$ sur $X$ est moyennable s'il existe une suite $(b^{n})_{n\in \mathbb{N}}$ d'applications de $X$ dans $\mathbb{P}(G)$ telle que $b^{n}$ est continue pour tout $n\in \mathbb{N}$ pour la topologie vague sur $\mathbb{P}(G)$ et $\underset{n\longrightarrow \infty}{\overset{}{\lim}}\,\,\underset{x\in X}{\overset{}{\sup}}\|gb^{n}_{x}-b^{n}_{gx}\|_{1}=0$ pour tout $g\in G$.
\end{definition}
\begin{definition}
Un groupe discret dénombrable $G$ est dit \index{groupe!moyennable à l'infini}topologiquement moyennable s'il existe un espace compact $X$ tel que:
\begin{enumerate}
  \item $G$ opère par homéomorphisme sur $X$,
  \item l'action de $G$ sur $X$ est moyennable.
\end{enumerate}
\end{definition}

Rappelons le critère de moyennabilité classique suivant dit condition de Reiter:
\begin{theorem}(Reiter \cite{paterson})\index{théorème!de Reiter}
Soit $p$ un nombre réel tel que $1\leq p< \infty$. Un groupe discret dénombrable $G$ est moyennable si et seulement s'il vérifie la condition de Reiter: Pour tout compact $C\subset G$ et $\varepsilon> 0$, il existe $h\in \{f\in L^{p}(G):\,f\geq 0,\,\|f\|_{p}=1\}$ tel que $\|gh-h\|_{p}< \varepsilon\,\,\,\,\,\,\forall g\in C$.
\end{theorem}
Nous allons utiliser ce théorème pour établir le lien entre la moyennabilité et la moyennabilité topologique
\begin{theorem}\label{theoexamp}
Soit $G$ un groupe discret dénombrable. Les propositions suivantes sont équivalentes:
\begin{enumerate}
  \item $G$ est moyennable.
  \item L'action triviale de $G$ sur tout singleton est moyennable.
\end{enumerate}
\end{theorem}
\begin{prof}
\begin{enumerate}
  \item [$1\Longrightarrow 2$] Soit $G$ un groupe discret dénombrable moyennable. Supposons que $G=\{g_{1},g_{2}...\}$. Par la condition de Reiter, pour tout $n\in \mathbb{N}^{\star}$ et $F_n=\{g_{1},g_{2}...,g_n\}$, il existe une application $b^{n}:X\longrightarrow \mathbb{P}(G)$, où $X=\{x\}$ tel que $\|g_{i}b_{x}^{n}-b_{x}^{n}\|_{1}< \frac{1}{n}\,\,\,\,\,\,\forall i\leq n$. Ainsi, pour tout $g\in G$, il existe $i\in \mathbb{N}^{\star}$ tel que $g=g_i$ et ainsi $g\in F_n$ pour tout $n\geq i$. Ceci implique que $\|g_{i}b_{x}^{n}-b_{x}^{n}\|_{1}< \frac{1}{n}$. Donc $\underset{n\longrightarrow \infty}{\overset{}{\lim}}\,\,\underset{x\in X}{\overset{}{\sup}}\|gb^{n}_{x}-b^{n}_{gx}\|_{1}=0$
  \item [$2\Longrightarrow 1$] Soit $G$ un groupe dénombrable tel que son action sur un singleton $X=\{x\}$ est moyennable. Alors il existe une suite d'applications ($b^{n})_{n\in \mathbb{N}}$ de $X$ dans $\mathbb{P}(G)$ continue pour la topologie vague sur $\mathbb{P}(G)$ et telle que $\underset{n\longrightarrow \infty}{\overset{}{\lim}}\,\,\underset{x\in X}{\overset{}{\sup}}\|gb^{n}_{x}-b^{n}_{gx}\|_{1}=0$. Soit $F$ un sous-ensemble fini de $G$ et soit $\varepsilon> 0$. Puisque l'action est moyennable, il existe $N$ tel que pour tout $n> N$ et pour tout $g\in F$, on a:  $\|gb_{x}^{n}-b_{gx}^{n}\|_{1}< \varepsilon$. En d'autres mots, $G$ vérifie la condition de Reiter pour $p=1$.
\end{enumerate}
\end{prof}
\begin{example}
\begin{enumerate}
  \item Tout groupe moyennable $G$ est topologiquement moyennable. Ceci est une conséquence immédiate du théorème \ref{theoexamp}.
  \item Le groupe libre à deux générateurs $F_2$ est topologiquement moyennable (\cite{anan} exemple $2.7$, \cite{ozawa} exemple $2.2$). La moyennabilité topologique est donc une propriété plus générale que la moyennabilité au sens classique.
\item Soit $G$ un groupe topologique. Si $G$ est compact, alors l'action à gauche de $G$ sur lui même est moyennable. En effet, pour tout $n\geq 1$, définissons  $b^{n}:G\longrightarrow \mathbb{P}(G)$ par $b^{n}(x)=\delta_x$ pour tout $x\in G$. Comme $tb^{n}(x)=b^{n}(tx)$ pour tous $t,x\in G$, la suite $(b^{n})_{n\geq 1}$ vérifie les conditions de la défintion \ref{def1}.
\end{enumerate}
\end{example}

\begin{remark}
Certains auteurs utilisent moyennabilité à l'infini ou encore groupe de Higson-Roe pour désigner les groupes topologiquement moyennables.
\end{remark}

\section{Moyennabilité topologique et Compactifié de Stone-\v Cech}

\begin{lemma}\label{lemmoytopo1}
Soient $X$ et $Y$ deux $G$-espaces compacts. Si $G$ opère moyennablement sur $X$ et $f:Y\longrightarrow X$ est une application équivariante, alors l'action de $G$ sur $Y$ est moyennable.
\end{lemma}
\begin{prof}
Soit $(b^{n})_{n\in \mathbb{N}}$ une suite d'applications de $X$ dans $\mathbb{P}(G)$ continue pour la topologie vague sur $\mathbb{P}(G)$ et telle que  $\underset{n\longrightarrow \infty}{\overset{}{\lim}}\,\,\underset{x\in X}{\overset{}{\sup}}\|gb^{n}_{x}-b^{n}_{gx}\|_{1}=0$ pour tout $g\in G$ et soit $f:Y\longrightarrow X$ une application équivariante. Posons
\[\begin{array}{lll}
c^{n}:& Y &\longrightarrow \mathbb{P}(G)\\
&y &\longmapsto c^{n}_{y}=b^{n}_{f(y)}\\
\end{array}\].

Alors, on a:
 \begin{tabular}{lll}
$\underset{y\in Y}{\overset{}{\sup}}\|gc^{n}_{y}-c^{n}_{gy}\|_{1}$
& $=$ &
$\underset{y\in Y }{\overset{}{\sup}}\|gb^{n}_{f(y)}-b^{n}_{f(gy)}\|_{1}$ \\
& $=$ & $ \underset{y\in Y }{\overset{}{\sup}}\|gb^{n}_{f(y)}-b^{n}_{gf(y)}\|_{1}$ \\
& $= $ & $ \underset{x\in X}{\overset{}{\sup}}\|gb^{n}_{x}-b^{n}_{gx}\|_{1}$
\end{tabular}

Donc $\underset{n\longrightarrow \infty}{\overset{}{\lim}}\,\underset{y\in Y}{\overset{}{\sup}}\|gc^{n}_{y}-c^{n}_{gy}\|_{1}=0$
\end{prof}
\begin{lemma}\label{lemmoytopo2}
Si $G$ possède une action moyennable sur un espace compact, alors $G$ possède une action moyennable sur son compactifié de Stone-\v Cech $\beta G$
\end{lemma}
\begin{prof}
Soit $g\in G$, notons

\[\begin{array}{lll}
L_{g}:& G &\longrightarrow G\\
&h &\longmapsto L_{g}(h)=gh\\
\end{array}\]

l'action continue à gauche de $G$ sur $G$. $L_{g}$ se prolonge de manière unique en $\widetilde{L}_{g}:\beta G\longrightarrow \beta G$ de telle manière que le diagramme suivant:
\[\xymatrix{G\ar[rr]^{L_{g}}\ar[dd]_{}&&G\ar[dd]^{}\\
\\\beta G\ar[rr]^{\widetilde{L}_{g}}&&\beta G}\]
 soit commutatif.\\
 L'application
 \[\begin{array}{lll}
& G &\longrightarrow \beta G\\
&g &\longmapsto \widetilde{L}_{g}\\
\end{array}\]
permet de définir une action par homéomorphismes de $G$ sur $\beta G$.\\
 Soit $X$ un $G$-espace compact tel que l'action de $G$ sur $X$ est moyennable. L'action $\tau:G\times X\longrightarrow X$ permet de définir une application équivariante $ \widetilde{\tau}:\beta G\longrightarrow X$.
 L'action de $G$ sur $\beta G$ est donc moyennable par le Lemme \ref{lemmoytopo1}
\end{prof}\\
Nous allons utiliser les lemmes suivants:
\begin{lemma}(\cite{higson},  Lemme $3.7$)\label{lem1}
Un groupe discret dénombrable $G$ admet une action moyennable sur son compactifié de Stone-\v Cech $\beta G$ si et seulement s'il existe une suite d'applications $(b^{n})_{n\in \mathbb{N}}$ de $G$ dans $\mathbb{P}(G)$ telle que:
\begin{enumerate}
  \item Pour tout $n$, l'image de l'application $b^{n}$ est contenu dans un sous-ensemble compact  de $\mathbb{P}(G)$ pour la topologie vague.
  \item Pour tout $g\in G$, on a: $\underset{n\longrightarrow \infty}{\overset{}{\lim}}\,\,\underset{h\in G}{\overset{}{\sup}}\|gb^{n}_{h}-b^{n}_{gh}\|_{1}=0$
\end{enumerate}
\end{lemma}
\begin{remark}
\begin{enumerate}
  \item Si $b\in\ell^{1}(Z)$ et $F$ est un sous-ensemble de $Z$, alors nous noterons $b|_{F}$ la fonction définie sur $Z$ par
$$b|_{F}(z)=\left\{
    \begin{array}{ccccc}
      b(z) &\text{si}& z\in F & &  \\
       0 & \text{si} & z\notin F&  \\
    \end{array}\right.$$
  \item Le résultat suivant apparait dans \cite{higson}. Mais l'auteur ne fournit pas une preuve complète.
\end{enumerate}
\end{remark}
\begin{lemma}\label{lem2}
Soit $Z$ un ensemble discret. Pour tout sous-ensemble $B$ de $\mathbb{P}(Z)$ compact pour la topologie vague et pour tout $\varepsilon> 0$, il existe $F\subset Z$ fini, tel que $\|b-b|_{F}\|<\varepsilon$ pour tout $b\in B$.
\end{lemma}
\begin{prof}
Fixons $\varepsilon> 0$. Pour tout sous-ensemble fini $H\subset Z$, posons $$U_{H}=\{b\in \mathbb{P}(Z):\,\,\|b|_{H}\|_{1}>1-\varepsilon\}.$$ Nous allons montrer que les ensembles $U_H$ forment un recouvrement étoile-faible ouvert de $\mathbb{P}(Z)$. Fixons $p\in Z$ et soit $g_{p}$ la fonction définie sur $c_{0}(Z)$ par

\[ g_{p}(z)= \left\{ \begin{array}{ll}
1 & \textrm{si $z=p$}\\
0&\textrm{sinon}\end{array}\right.\]

Alors, $$b(p)=<b,\,g_p>=\underset{z\in Z}{\overset{}{\sum}}b(z)g_{p}(z).$$ Ainsi, l'application $b\longmapsto b(p)$ qui coïncide avec l'application $b\longmapsto <b,\,g_p>$ est continue pour la topologie vague pour tout point $p\in Z$.\\
Maintenant pour un sous-ensemble fini $H=\{p_{1},p_{2},...p_{k}\}$ de $Z$, les applications
\[\begin{array}{lll}
& \mathbb{P}(Z) &\longrightarrow \mathbb{R}\\
&b &\longmapsto b(p_{i})\\
\end{array}\]

sont continues pour la topologie vague pour $i=1,2,...,k$. Ainsi l'application

\[\begin{array}{lll}
& \mathbb{P}(Z) &\longrightarrow \mathbb{R}\\
&b &\longmapsto b(p_{1})+b(p_{2})+...+b(p_{k})\\
\end{array}\]

est continue pour la topologie vague. Comme $$b(p_{1})+b(p_{2})+...+b(p_k)=\|b|_H\|_1,$$ l'application $b\longmapsto\|b|_H\|_1$ est continue. Donc pour tout $r\in \mathbb{R}$, l'ensemble $$\{b\in \mathbb{P}(Z):\,\,\|b|_{H}\|_{1}>1-r\}$$ est un ouvert. Ainsi, $U_H$ est un ouvert pour tout $H\subset Z$ fini.\\
Finalemement pour tout $b\in \mathbb{P}(Z)$, on a: $\underset{z\in Z}{\overset{}{\sum}}b(z)=1$. Comme cette série est à termes positifs, il existe un sous-ensemble fini $H$ de $Z$ tel que $\underset{z\in Z\setminus H}{\overset{}{\sum}}b(z)<\varepsilon$. Ainsi  $\underset{z\in H}{\overset{}{\sum}}b(z)>1-\varepsilon$ et $\|b|_{H}\|_{1}>1-\varepsilon$ i.e $b\in U_H$. Donc les ensembles $U_H$ forment un recouvrement étoile-faible de $\mathbb{P}(Z)$ lorsque $H$ parcourt l'ensemble des sous-ensemble fini de $Z$ et $\varepsilon$ les nombres réels positifs. Ainsi, il existe des sous-ensembles finis $H_{1},H_{2},...H_{N}$ de $Z$ tel que $B\subset \underset{i=1}{\overset{N}{\bigcup}}U_{H_i}$. Prendre $F=\underset{i=1}{\overset{N}{\bigcup}}H_i$.
\end{prof}
\begin{lemma}\label{lemhig}
Un groupe discret dénombrable $G$ admet une action moyennable sur son compactifié de Stone-\v Cech $\beta G$ si et seulement s'il existe une suite $(b^{n})_{n\in\mathbb{N}}$ d'applications de $G$ dans $\mathbb{P}(G)$ telle que:
\begin{enumerate}
  \item Pour tout $n$, il existe $F_{n}\subset G$ fini tel que $supp(b^{n}_{g})\subset F_{n}$ pour tout $g\in G$
  \item $\underset{n\longrightarrow \infty}{\overset{}{\lim}}\,\,\underset{h\in G}{\overset{}{\sup}}\|gb^{n}_{h}-b^{n}_{gh}\|_{1}=0$ pour tout $g\in G$
\end{enumerate}
\end{lemma}

\begin{prof}
\begin{enumerate}
  \item [$\Longrightarrow$] Supposons $G$ topologiquement moyennable. Il existe par le lemme \ref{lem1} une suite $(b_{n})_{n\in\mathbb{N}}$ d'applications de $G$ dans $\mathbb{P}(G)$ telle que pour tout $n\in \mathbb{N},\,\,\,b^{n}(G)$ est contenu dans une partie compacte $K_{n}$ de $\mathbb{P}(G)$ et $\underset{n\longrightarrow \infty}{\overset{}{\lim}}\,\,\underset{h\in G}{\overset{}{\sup}}\|gb^{n}_{h}-b^{n}_{gh}\|_{1}=0$ pour tout $g\in G$.\\
      $K_n$ étant une partie compacte de $\mathbb{P}(G)$, d'après le lemme \ref{lem2}, pour tout $\frac{\varepsilon}{5}>0$, il existe un sous-ensemble fini $F_{n}$ de $G$ tel que $\|b-b|_{F_{n}}\|_{1}<\frac{\varepsilon}{5}$ pour tout $b\in b^{n}(G)$. Définissons $\forall\,n\in \mathbb{N}$, les applications $w^{n}$ de $G$ dans $\mathbb{P}(G)$ par $$w_{g}^{n}=\dfrac{b_{g}^{n}|F_{n}}{\|b_{g}^{n}|_{F_{n}}\|_{1}}.$$ Ainsi, nous avons les observations suivantes:
  \begin{enumerate}
    \item $supp(w^{n}_{g})\subset F_{n}$ pour tout $g\in G$ et
    \item $\underset{n\longrightarrow \infty}{\overset{}{\lim}}\,\,\underset{h\in G}{\overset{}{\sup}}\|gw^{n}_{h}-w^{n}_{gh}\|_{1}=0$ pour tout $g\in G$.\\
        En effet, à partir de $$\|b_{g}^{n}|_{F_{n}}\|_{1}> 1-\frac{\varepsilon}{5},$$ on a:  $$\|gw^{n}_{h}-gb^{n}_{h}|_{F_{n}}\|_{1}\leq \frac{\varepsilon}{5}.$$ En prenant $n$ suffisamment grand, nous obtenons:
        $$\|gb^{n}_{h}-b^{n}_{gh}\|_{1}\leq \frac{\varepsilon}{5}$$(puisque $\underset{n\longrightarrow \infty}{\overset{}{\lim}}\,\,\underset{h\in G}{\overset{}{\sup}}\|gb^{n}_{h}-b^{n}_{gh}\|_{1}=0$ pour tout $g\in G$).\\
        De plus, on a:

        $\|gw^{n}_{h}-w^{n}_{gh}\|_{1}$

        \begin{tabular}{lllll}
& $=$ & $ \|gw^{n}_{h}- gb^{n}_{h}|_{F_{n}}+gb^{n}_{h}|_{F_{n}}-gb^{n}_{h}+gb^{n}_{h}-b^{n}_{gh}+b^{n}_{gh}-b^{n}_{gh}|_{F_{n}}+b^{n}_{gh}|_{F_{n}}-w^{n}_{gh}\|_{1}$ \\
& $\leq$ & $\|gw^{n}_{h}- gb^{n}_{h}|_{F_{n}}\|_{1}+\|gb^{n}_{h}|_{F_{n}}-gb^{n}_{h}\|_{1}+\|gb^{n}_{h}-b^{n}_{gh}\|_{1}$\\&
$+$ & $ \|b^{n}_{gh}-b^{n}_{gh}|_{F_{n}}\|_{1}+
\|b^{n}_{gh}|_{F_{n}}-w^{n}_{gh}\|_{1} $ \\
& $\leq $ & $ \frac{\varepsilon}{5}+\frac{\varepsilon}{5}+\frac{\varepsilon}{5}+\frac{\varepsilon}{5}+\frac{\varepsilon}{5}=\varepsilon$
\end{tabular}

  \end{enumerate}
  \item [$\Longleftarrow$] Supposons qu'il existe une suite $(w^{n})_{n\in\mathbb{N}}$ d'applications de $G$ dans $\mathbb{P}(G)$ telle que les conditions $1$ et $2$ sont remplies. Pour tout $n$, il existe $F_n$ fini tel que $supp(b^{n}_{g})\subset F_{n}$ pour tout $g\in G$. Puisque l'ensemble des mesures sur $G$ à support fini est un sous-espace compact de $\mathbb{P}(G),\,w^{n}(G)$ est contenu dans un sous-espace compact de $\mathbb{P}(G)$. Donc l'action de $G$ sur son compactifié de Stone-\v Cech est moyennable par le lemme \ref{lem1}.
\end{enumerate}
\end{prof}

\section{Espaces test pour la moyennabilité topologique}

En utilisant les notations du chapitre précédent, on a le lemme suivant:
\begin{lemma}\label{lemmoytopo}
Soit $G$ un groupe discret dénombrable moyennable à l'infini. Notons $(b^{n})_{n\in \mathbb{N}}$ la suite d'applications correspondantes de $S(G)=\beta G$ dans $\mathbb{P}(G)$ et $F=\{b^{n}:n\in \mathbb{N}\}$. Alors l'action de $G$ sur $X_F$ est moyennable.
\end{lemma}
\begin{prof}
$G$ étant moyennable à l'infini, alors il existe une suite $(b^{n})_{n\in\mathbb{N}}$ d'applications de $\beta G$ dans $\mathbb{P}(G)$ telle que: $\underset{n\longrightarrow \infty}{\overset{}{\lim}}\,\,\underset{x\in \beta G}{\overset{}{\sup}}\|gb^{n}_{x}-b^{n}_{gx}\|_{1}=0$ pour tout $g\in G$. Rappelons que $G$ étant discret et dénombrable, $S(G)=\beta G$.\\
Pour tout $g\in G$, notons
\[\begin{array}{lll}
\overline{g}:& \beta G &\longrightarrow \beta G\\
&x &\longmapsto gx\\
\end{array}\]
l'homéomorphisme de $\beta G$ sur lui même produit par $g$. Considérons le produit diagonal
$$f=\Delta_{(g,n)\in G\times \mathbb{N}}(b_{n}\circ \overline{g}):\beta G\longrightarrow (\mathbb{P}(G))^{G\times \mathbb{N}}$$ défini par $$f(x)=(b_{gx}^{n})_{(g,n)\in G\times \mathbb{N}}.$$
Il est clair que $f$ est continue. La relation d'équivalence $\mathcal{R}_{F}$ est définie sur $\beta G$ par:

 $$x\mathcal{R}_{F}y\Longleftrightarrow b_{gx}^{n}=b^{n}_{gy}$$ pour tout $n\in \mathbb{N}$ et $g\in G$.
   Notons encore $f$ l'application $f:\beta G\longrightarrow f(\beta G)$, et $\pi_F:\beta G\longrightarrow \beta G/ \mathcal{R}_{F}$ la surjection canonique. Il existe une application continue $\bar f$ telle que: $f=\bar f\circ \pi_F$. Posons $X_F=\beta G/ \mathcal{R}_{F}$.
    Considérons l'application $\widetilde{b}^{n}:f(\beta G)\longrightarrow \mathbb{P}(G)$ définie par $\widetilde{b}^{n}=\pi_{e,n}$ où $\pi_{e,n}$ est la projection. Ainsi, l'application $$c^{n}=\widetilde{b}^{n}\circ \bar f:X_F\longrightarrow \mathbb{ P}(G)$$ est continue pour la topologie vague sur $\mathbb{ P}(G)$.\\
    Soit $g\in G$, on a:

    \begin{tabular}{lll}

$\underset{[x]\in X_F }{\overset{}{\sup}}\|gc^{n}_{[x]}-c^{n}_{g[x]}\|_{1}$
& $=$ &
$\underset{[x]\in X_F }{\overset{}{\sup}}\|g(\widetilde{b}^{n}\circ \bar f)_{[x]}-(\widetilde{b}^{n}\circ \bar f)_{g[x]}\|_{1}$ \\
& $=$ & $ \underset{[x]\in X_F }{\overset{}{\sup}}\|g(\widetilde{b}^{n}(\bar f([x])))-\widetilde{b}^{n}(\bar f(g[x]))\|_{1}$ \\
& $= $ & $ \underset{x\in\beta G}{\overset{}{\sup}}\|g(\widetilde{b}^{n}(\bar f([x])))-\widetilde{b}^{n}(\bar f([gx]))\|_{1}$ \\
& $=$ & $\underset{x\in\beta G}{\overset{}{\sup}}\|g(\widetilde{b}^{n}(f(x)))-\widetilde{b}^{n}(f(gx))\|_{1}$\\
& $=$ & $\underset{x\in\beta G}{\overset{}{\sup}}\|g(\widetilde{b}^{n}(b_{hx}^{n}))_{(h,n)\in G\times \mathbb{N}})-\widetilde{b}^{n}((b_{hgx}^{n}))_{(h,n)\in G\times \mathbb{N}}\|_{1}$\\
& $=$ & $\underset{x\in\beta G}{\overset{}{\sup}}\|gb_{x}^{n}-b^{n}_{gx}\|_{1}$
\end{tabular}

Ainsi, $\underset{n\longrightarrow \infty}{\overset{}{\lim}}\,\underset{[x]\in X_{F} }{\overset{}{\sup}}\|gc^{n}_{[x]}-c^{n}_{g[x]}\|_{1}=0$
\end{prof}
\begin{corollary}
Un groupe discret dénombrable $G$ admet une action moyennable sur un espace compact et métrisable si et seulement si son action sur son Compactifié de Stone-\v Cech $\beta G$ est moyennable.
\end{corollary}
\begin{prof}
Ceci est une conséquence immédiate du lemme \ref{lemmoytopo2} et du lemme \ref{lemmoytopo}.
\end{prof}
\begin{theorem}
Un groupe discret dénombrable $G$ est moyennable à l'infini si et seulement s'il admet une action moyennable sur l'ensemble de Cantor $D^{\aleph_{0}}$.
\end{theorem}

\begin{prof}
La suffisance est évidente. Montrons la necessité.\\Si $G$ est moyennable à l'infini, alors $G$ possède une action moyennable sur son compactifié de Stone-\v Cech $\beta G$. Ainsi, il existe une suite $(b^{n})_{n\in \mathbb{N}}$ d'applications de $\beta G$ dans $\mathbb{P}(G)$ telle que
$\underset{n\longrightarrow \infty}{\overset{}{\lim}}\,\,\underset{x\in \beta G}{\overset{}{\sup}}\|gb^{n}_{x}-b^{n}_{gx}\|_{1}=0$ pour tout $g\in G$\\
Pour tout $g\in G$, notons
\[\begin{array}{lll}
\overline{g}:& \beta G &\longrightarrow \beta G\\
&x &\longmapsto gx\\
\end{array}\]
l'homéomorphisme de $\beta G$ sur lui même. L'espace $b^{n}(\beta G)$ étant compact et métrisable, il existe une surjection continue $f^{n}:D^{\aleph_{0}}\longrightarrow b^{n}(\beta G)$.\\ Soit $g\in G$, alors $b_{g}^{n}\in b^{n}(\beta G)$. Puisque $f^{n}$ est surjective, il existe $c_{g}\in D^{\aleph_{0}}$ tel que $f^{n}(c_{g})=b_{g}^{n}$. On a ainsi une application

\[\begin{array}{lll}
T^{n}:& G &\longrightarrow D^{\aleph_{0}}\\
&g &\longmapsto c_{g}\\
\end{array}.\]
Cette application se prolonge de manière unique en une application continue $$\beta T^{n}:\beta G\longmapsto D^{\aleph_{0}}.$$
Pour tout $g\in G$, on a:
$$(f^{n}\circ \beta T^{n})(g)=f^{n}(\beta T^{n}(g))=f^{n}(c_{g})=b_{g}^{n}.$$
Ainsi, $f^{n}\circ \beta T^{n}=b^{n}$ sur $G$. Puisque $G$ est dense dans $\beta G,\,\,\,f^{n}\circ \beta T^{n}=b^{n}$ sur $\beta G$.\\
Posons: $c^{n}=\beta T^{n}\circ \overline{g}:\beta G\longrightarrow D^{\aleph_{0}},\,\,F=\{c^{n}:n\in \mathbb{N}\}$. Notons encore $f^{n}:D^{\aleph_{0}}\longrightarrow f^{n}(D^{\aleph_{0}})$. Posons $\widetilde{c}^{n}=f^{n}\circ\beta T^{n}\circ \overline{g}:\beta G\longrightarrow f^{n}(D^{\aleph_{0}}),\,\,F=\{\widetilde{c}^{n}:n\in \mathbb{N}\}$  et $X=f^{n}(D^{\aleph_{0}})$. Alors $X_{F}\cong D^{\aleph_{0}}$ et l'action de $G$ sur $X_F$ est moyennable à l'infini par le lemme \ref{lemmoytopo}.
\end{prof}
\begin{remark}
Le résultat précédent est le principal résultat de \cite{yousef}. La preuve que nous venons d'effectuer est différente de celle de \cite{yousef}.
\end{remark}
\begin{lemma}\label{lemcantor}
Soit $G$ un groupe discret dénombrable. Si l'action de $G$ sur l'ensemble de Cantor $D^{\aleph_{0}}$ est moyennable, alors les images des $b^{n}$ peuvent être choisis d'images finies.
\end{lemma}
\begin{prof}
Soit $(b^{n})_{n\in \mathbb{N}}$ une suites d'applications de $D^{\aleph_{0}}$ dans $\mathbb{P}(G)$ telle que $\underset{n\longrightarrow \infty}{\overset{}{\lim}}\,\,\underset{x\in D^{\aleph_{0}}}{\overset{}{\sup}}\|gb^{n}_{x}-b^{n}_{gx}\|_{1}=0$ pour tout $g\in G$.
Puisque l'espace de Cantor $D^{\aleph_{0}}$ a la dimension de Lebesgue zéro, il existe une partition finie $\gamma=(A_{i})_{i=1,2,..,k_{n}}$ de $D^{\aleph_{0}}$ en sous-ensembles ouverts et fermés, telle que l'image par $b^{n}$ de chaque élément de $\gamma$ est contenue dans une des boules $B_{\frac{1}{n}}$ de rayon $1/n$ dans $\ell^1(G)$. Notons $c_{i}$ les centres des boules correspondantes.
  Pour tout $A_{i}\in\gamma$, considérons l'application $(b^{n})^\prime$ définie par :
   $$(b^{n})_{x}^\prime=c_{i}\,\, \text{si}\,\, x\in A_{i}.$$ Soit $x\in D^{\aleph_{0}}$, il existe par définition de $b^{n}$ et $(b^{n})^{\prime}$, une boule $B_{\frac{1}{n}}$ telle que:
        $$\|b^n_x-(b^n)'_x\|_{1}\leq \underset{p,q \in B_{\frac{1}{n}}}{\overset{}{\sup}}\|p-q\|_{1}=diam(B_{\frac{1}{n}})=\frac{2}{n}.$$
         D'où $$\underset{n\longrightarrow \infty}{\overset{}{\lim}}\,\,\underset{x\in D^{\aleph_{0}} }{\overset{}{\sup}}\|g(b^{n})^{\prime}_{x}-(b^{n})^{\prime}_{gx}\|_{1}=0.$$
\end{prof}

\begin{theorem}\label{corohilbert}
Un groupe discret dénombrable $G$ est moyennable à l'infini si et seulement s'il admet une action moyennable sur le cube de Hilbert $I^{\aleph_{0}}.$
\end{theorem}

\begin{prof}
La suffisance est évidente. Montrons la necessité .\\
Si $G$ est moyennable à l'infini, alors $G$ admet une action moyennable sur l'ensemble de Cantor. Ainsi, il existe une suite $(b^{n})_{n\in \mathbb{N}}$ d'applications de $D^{\aleph_{0}}$ dans $\mathbb{P}(G)$ telle que $\underset{n\longrightarrow \infty}{\overset{}{\lim}}\,\,\underset{x\in D^{\aleph_{0}}}{\overset{}{\sup}}\|gb^{n}_{x}-b^{n}_{gx}\|_{1}=0$ pour tout $g\in G$.\\
 Par le lemme \ref{lemhig}, on peut supposer sans perte de généralité que pour tout $n$, il existe $F_n \subseteq G$ fini tel que $supp(b^{n}_{x})\subset F_{n}$ pour tout $x\in D^{\aleph_{0}}$. Autrement dit, $b^{n}(D^{\aleph_{0}})$ est contenu dans un sous-espace de dimension finie $V$ de $\mathbb{P}(G)$. Par le lemme \ref{lemcantor}, supposons que les images des $b^n$ sont finies et notons $(c_{i})_{i\in I_{n}}$ où $I_{n}$ est fini les images de tous les $b^n$. Posons $A_{i}=(b^{n})_{c_{i}}^{-1}$ et considérons les applications

 \[\begin{array}{lll}
c^{n}:& \mathbb{P}(D^{\aleph_{0}}) &\longrightarrow \mathbb{P}(G)\\
&\mu &\longmapsto \underset{i=1}{\overset{n}{\Sigma}} \mu(A_{i}) c_{i}\\
\end{array}\]

Les applications $c^{n}$ sont clairement affines. Ce sont précisement les prolongements affines des applications $b^{n}$ sur $\mathbb{P}(D^{\aleph_{0}})$.
Les $c^{n}$ sont continues rapport à la topologie vague sur $\mathbb{P}(D^{\aleph_{0}})$.\\ En effet, les applications $b^{n}:D^{\aleph_{0}}\longrightarrow V$ sont continues.\\
Si $\phi$ une fonctionnelle linéaire sur $V$, alors $\phi\circ b^{n}$ est une fonction continue sur $D^{\aleph_{0}}$. Par définition de la topologie vague sur $P(D^{\aleph_{0}})$, l'extension unique, $\widetilde{\phi\circ b^{n}}$, de $\phi\circ b^{n}$ sur $P(D^{\aleph_{0}})$ est continue. Par unicité du prolongement, $\widetilde{\phi\circ b^{n}}=\phi\circ \widetilde{b^{n}}=\phi\circ c^{n}$ sur $P(D^{\aleph_{0}})$. Ainsi, $c^{n}$ est continue par rapport à la topologie vague sur $P(D^{\aleph_{0}})$.\\
Notons $\mathbb{P}_{0}(D^{\aleph_0})$ le sous-espace de $\mathbb{P}(D^{\aleph_{0}})$ formé des mesures à support fini.\\
     Soit $\mu=\underset{i=1}{\overset{m}{\Sigma}}\alpha_{i}\delta_{x_{i}}\in \mathbb{P}_{0}(D^{\aleph_0})$, on a: $c^{n}_{\mu}=\underset{i=1}{\overset{m}{\Sigma}}\alpha_{i}c^{n}_{\delta_{x_{i}}}$.
     Or $$c^{n}_{\delta_{x_{i}}}=\underset{j=1}{\overset{n}{\Sigma}}\delta_{x_{i}}(A_{j})c_{j}=c_{i}=b^{n}_{x_{i}}.$$ D'où $$c^{n}_{\mu}=\underset{i=1}{\overset{m}{\Sigma}}\alpha_{i}b^{n}_{x_{i}}.$$ De même, $$c^n_{g\mu}=\Sigma_i \alpha_i b^n_{g x_i},$$ car $$g\mu=g\underset{i=1}{\overset{m}{\Sigma}}\alpha_{i}\delta_{x_{i}}=\underset{i=1}{\overset{m}{\Sigma}}\alpha_{i}\delta_{gx_{i}}.$$
      Pour tout $ \varepsilon> 0$, on a: $$\|gc^n_\mu - c^n_{g\mu} \|_{1}= \|\Sigma_i \alpha_i gb^n_{x_i} - \Sigma_i \alpha_i b^n_{g x_i}\|_{1}\leq \Sigma_i \alpha_i \|gb^n_{x_i} - b^n_{g x_i}\|_{1} \leq \Sigma_i \alpha_i \epsilon = \epsilon.$$
      $\mathbb{P}_{0}(D^{\aleph_0})$ étant dense dans $\mathbb{P}(D^{\aleph_{0}})$ pour la topologie vague et les applications $c^{n}$ étoile-faible continues, on a:  $$\underset{\mu \in \mathbb{P}(D^{\aleph_{0}})}{\overset{}{\sup}}\|gc^{n}_{\mu}-c^{n}_{g\mu}\|_{1}<\varepsilon.$$ D'où $$\underset{n\longrightarrow \infty}{\overset{}{\lim}}\,\underset{\mu \in \mathbb{P}(D^{\aleph_{0}})}{\overset{}{\sup}}\|gc^{n}_{\mu}-c^{n}_{g\mu}\|_{1}=0.$$
On conlut que l'action de $G$ sur $\mathbb{P}(D^{\aleph_{0}})$ est moyennable. Par le théorème de Keller, $\mathbb{P}(D^{\aleph_{0}})$ est homéomorphe à $I^{\aleph_{0}}$. Donc l'action de $G$ sur $I^{\aleph_{0}}$ est moyennable.
\end{prof}

\chapter{Groupe des isométries de l'espace d'Urysohn-Kat\v etov}
Dans ce chapitre, on appelle application isométrique ou plongement isométrique toute application $f:(X,d)\longrightarrow (Y,d^{\prime})$ telle que $d^{\prime}(f(a),f(b))=d(a,b)$ pour tous $a,b\in X$. Si $f$ est de plus bijective, on dira que $f$ est une isométrie.
\section{Généralités sur les groupes universels}
\begin{definition}
Soit $\mathcal{C}$ une classe d'espaces topologiques. $X\in \mathcal{C}$ est dit universel\index{espace!universel} pour cette classe si pour tout $Y\in \mathcal{C}$, il existe un homéomorphisme de $Y$ sur un sous- espace de $X$.
\end{definition}
L'espace de Cantor $D^{\aleph_{0}}$ est par exemple universel pour la classe des espaces topologiques métrisables séparables et de dimension de Lebesgue $0$ (\cite{eng}).
\begin{definition}
Un groupe topologique $G$ est dit universel pour une classe $\mathcal{C}$ de groupes topologiques si pour tout groupe topologique $H\in \mathcal{C}$, il existe un isomorphisme de groupes topologiques de $H$ sur un sous-groupe de $G$.
\end{definition}
En réponse à une question de Ulam (cf. Problème $103$ dans \cite{ulam}), Uspenskij montre en $1986$ dans (\cite{us4}) que le groupe $Homeo(I^{\aleph_{0}})$ de tous les homéomorphismes du cube de Hilbert $I^{\aleph_{0}}$ sur lui-même muni de la topologie compact-ouvert est universel pour la classe des groupes métrisables et séparables. Quelques années plutard, il démontre dans \cite{us5} que le groupe $Iso(\mathbb{U})$ des isométries de l'espace universel polonais d'Urysohn $\mathbb{U}$ sur lui-même muni de la topologie de la convergence simple est universel pour la même classe de groupes topologiques. La question d'existence d'un groupe topologique universel pour la classe des groupes de poids non dénombrables reste ouverte. Cette question sera traitée à la fin de ce chapitre.
\section{Le Théorème d'Uspenkij}\index{théorème!d'Uspenkij}
\subsection{Espace Universel d'Urysohn}
L'espace universel polonais d'Urysohn a été construit en réponse à la question suivante posée par Fréchet: existe-t-il un espace métrique séparable $X$ universel pour la classe des espaces métriques séparables? C'est en $1925$ qu'Urysohn apporte une réponse par l'affirmative à cette question dans \cite{ury1}. Cette section regroupe la construction de Kat\v etov(\cite{kat}) de l'espace universel polonais d'Urysohn et les propriétés générales de cet espace. Elle est particulièrment inspirée de Melleray (\cite{mel}) et Pestov(\cite{vp1}).
\begin{definition}
L'espace métrique universel polonais d'Urysohn\index{espace!d'Urysohn} $\mathbb{U}$ est caractérisé à isométrie près par les propriétés suivantes:
\begin{enumerate}
  \item $\mathbb{U}$ est complet et séparable
  \item $\mathbb{U}$ contient une copie isométrique de tout autre espace métrique complet et séparable ($\mathbb{U}$ est universel)
  \item Toute isométrie entre deux sous-ensembles finis de $\mathbb{U}$ se prolonge en une isométrie de $\mathbb{U}$ ($\mathbb{U}$ est $\omega$-homogène)
\end{enumerate}
\end{definition}
Le résultat suivant est démontré par Urysohn dans \cite{ury} (voir aussi \cite{vp1}).
\begin{theorem}
L'espace métrique séparable universel d'Urysohn $\mathbb{U}$ existe et est unique à isométrie près.
\end{theorem}

\begin{remark}
Le résultat précédent a été démontré par Urysohn en 1925. La construction d'Urysohn a été longtemps éclipsée par le célèbre résultat de Banach et Mazur qui établi l'universalité de  $C([0, 1])$ dans la classe des espaces métriques séparables et complets. L'espace $C([0, 1])$ a néanmoins la faiblesse de ne pas être $\omega$-homogène. Ce résultat de Banach et Mazur a contribué à faire tomber dans l'oubli l'espace d'Urysohn. Ce dernier n'a été que très peu étudié pendant $60$ ans, à l'exception d'articles de Sierpinski (\cite{si}) et Huhunaisvili (\cite{hu}).\\
En $1986$, Kat\v etov(\cite{kat}) donne une nouvelle construction de $\mathbb{U}$, qui va considérablement relancé l'intérêt pour cet espace. Le théorème d'Uspenskij utilise la construction de Kat\v etov de l'espace universel polonais d'Urysohn $\mathbb{U}$.
\end{remark}
\begin{definition}
Soit $(X,d)$ un espace métrique. Une fonction $f:X\longrightarrow \mathbb{R}^{+}$ est dite de Kat\v etov\index{fonction!de Kat\v etov} si: $|f(x)-f(y)|\leq d(x,y)\leq f(x)+f(y)$ pour tous $x,y\,\in X$
\end{definition}
\begin{remark}
On note $E(X)$ l'ensemble des fonctions de Kat\v etov sur $X$
\begin{enumerate}
  \item Ces fonctions correspondent en réalité aux extensions métriques de $X$ par un point de la façon suivante: si $f\in E(X)$, on peut définir une distance sur l'espace $X\cup\{f\}$ en posant: $d(x,f)=f(x)$ pour tout $x\in X$.\\
Dire que $f$ est une fonction de Kat\v etov, revient à dire que $d$ ainsi définie est une distance; autrement dit, toutes les extensions métriques de $X$ par un point sont obtenues de cette façon.
\item Le grand intérêt des fonctions de Kat\v etov est qu'il existe une distance naturelle entre fonctions de Kat\v etov.\\
En effet, si $f,g\,\in E(X)$ et $x_{0},x\in X$, alors on a:\\
 $|f(x)-d(x,x_{0})|\leq f(x_{0})$ et $|g(x)-d(x,x_{0})|\leq g(x_{0})$\\
 Par conséquent, $|f(x)-g(x)|\leq f(x_{0})+g(x_{0})$, donc $\underset{x\in X}{\overset{}{\sup}}|f(x)-g(x)|$ est fini.\\
 Toute différence de fonctions de Kat\v etov est donc une fonction borné et on peut donc poser: $$d_{X}^{E}(f,g)=\|f-g\|_{\infty}=\underset{x\in X}{\overset{}{\sup}}|f(x)-g(x)|.$$
  \item  Puisqu'on veut voir $E(X)$ comme l'espace des extensions métriques de $X$ par un point, il est naturel de voir $X$ comme un sous-espace de $E(X)$, celui des extensions triviales par un point(i.e le point qu'on ajoute etait déja dans $X$).\\
 Analytiquement, ceci se fait via l'application de Kuratowski\index{application!de Kuratowski} $x\longmapsto \delta_{x}$ définie par: $\delta_{x}(x^{\prime})=d(x,x^{\prime})$ pour tout $x^{\prime}\in X$.\\
 L'inégalité triangulaire permet de vérifier que $x\longmapsto \delta_{x}$ est un plongement isometrique de $X$ dans $E(X)$.
 \item En général, si $x\in X$ et $f\in E(X)$, alors $d_{X}^{E}(f,\delta_x)=f(x)$.\\ En effet, puisque $f$ est de Kat\v etov, nous avons $f(y)-d(x,y)\leq f(x)$ et $d(x,y)-f(y)\leq f(x)$ pour tout $y\in X$. Ainsi, $d_{X}^{E}(f,\delta_x)=\underset{y\in X}{\overset{}{\sup}}|f(y)-d(x,y)|\leq f(x)$ et en faisant $y=x$, on obtient l'égalité.
  \item Dans la suite, lorsqu'on écrira $X\subset E(X)$, on identifiera toujours $X$ à son image dans $E(X)$ par l'application de Kuratowski
\end{enumerate}
\end{remark}
\begin{definition}
Soient $(X,d)$ un espace métrique et $Y\subset X$ une partie de $X$. Pour toute $f\in E(Y)$, on définit son extension\index{extension de Kat\v etov} de Kat\v etov $k_{Y}(f)$ à $X$ par: $$k_{Y}(f)(x)=\underset{y\in Y}{\overset{}{\inf}}\{f(y)+d(y,x)\},$$ pour tout $x\in X$
\end{definition}
\begin{remark}
\begin{enumerate}
  \item Si $f\in E(X)$ et $Y\subset X$ vérifient: $\forall\,x\in X,\,f(x)=\underset{y\in Y}{\overset{}{\inf}}\{f(y)+d(y,x)\}$, on dira que $f$ est contrôlée par $Y$ ou que $Y$ est un support pour $f$.\\
  Par exemple, la fonction de distance
  \[\begin{array}{lll}
f_{x_{0}}:& X &\longrightarrow \mathbb{R}\\
&x &\longmapsto f_{x_{0}}(x)=d(x,x_{0})\\
\end{array}\]
est contrôlée par le singleton $\{x_{0}\}$.

\item En général, $E(X)$ n'est pas séparable même si $X$ l'est. C'est pourquoi, on considère le sous-espace
$E(X,\omega)=\{f\in E(X):\,\,\,\text{f a un support fini}\}$ de $E(X)$.

\item Dans \cite{mel}, Melleray caractérise les espaces métriques $(X,d)$ tels que $E(X)$ est séparable

\end{enumerate}
\end{remark}

Nous allons rappeler les conditions obtenues par Melleray sur $(X,d)$ pour la séparabilité de $E(X)$.
\begin{definition}
Un espace métrique a la propriété de Heine-Borel si ses sous-ensembles bornés sont précompacts.
\end{definition}
\begin{proposition}(Proposition $1$ dans \cite{mel})
Si $X$ est polonais mais n'a pas la propriété de Heine-Borel, alors $E(X)$ n'est pas séparable.
\end{proposition}
\begin{definition}
Soit $(X,d)$ un espace métrique.
\begin{enumerate}
  \item Si $ \varepsilon> 0$, on dit qu'une suite $(u_{n})_{n\in \mathbb{N}}$ d'éléments de $X$ est $\varepsilon$-bien-alignée si on
a, pour chaque $r\geq0,\,\,\,\,\underset{i=0}{\overset{r}{\sum}}d(u_{i},u_{i+1})\leq d(u_{0},u_{r+1})+\varepsilon$.
  \item Une suite $(u_{n})_{n\in \mathbb{N}}$ d'éléments de $X$ est dite alignée si pour tout $\varepsilon > 0$, il existe $N>0$ tel que $(u_{0},u_{N},u_{N+1},...)$ est $\varepsilon$-bien-alignée.
\end{enumerate}
\end{definition}
\begin{theorem}(Théorème $2$ dans \cite{mel})
Soit $X$ un espace polonais. Les assertions suivantes sont équivalentes :
\begin{enumerate}
  \item $E(X) = \overline{E(X, \omega)}$
  \item $E(X)$ est séparable.
  \item $\forall\,\delta > 0\,\,\forall\,(x_{n})\,\exists\,N\in \mathbb{N}\,\,\forall\,n\geq N,\,\,\exists i\leq N\,\,\,d(x_{0},x_{n})\geq d(x_{0},x_{i})+d(x_{i},x_{n})-\delta$
  \item De toute suite d'éléments de $X$ on peut extraire une sous-suite alignée.
\end{enumerate}
\end{theorem}
Rappelons cette terminologie empruntée à \cite{kal}.

\begin{definition}
\begin{enumerate}
  \item Soit $\varepsilon>0$. Un triplet ordonné $\{x_{1},x_{2},x_{3}\}$ de points de $X$ est $\varepsilon$-colinéaire si $d(x_{1},x_{3})\geq d(x_{1},x_{2})+d(x_{2},x_{3})-\varepsilon$.
  \item Un espace polonais $X$ a la propriété de colinéarité si: Pour tout sous-ensemble infini $A$ de $X$ et tout $\varepsilon> 0$, il existe $x_{1},x_{2},x_{3}\in A$ (deux à deux distincts) tels que $\{x_{1},x_{2},x_{3}\}$ est $\varepsilon$-colinéaire.
\end{enumerate}
\end{definition}
Modulo le théorème de Kalton (\cite{kal}) qui affirme que tout espace métrique a la propriété de colinéarité si, et seulement si, toute suite
admet une sous-suite alignée, le résultat suivant apparait également dans \cite{mel}.
\begin{corollary}(Corollaire $3$ dans \cite{mel})
Un espace polonais $X$ a la propriété de colinéarité si, et seulement si, $E(X)$ est séparable.
\end{corollary}
Ce qui suit est un résultat de Kat\v etov démontré dans \cite{mel}.
\begin{proposition}Soit $X$ un espace métrique. $E(X,\omega)$ est séparable si $X$ est séparable. De plus, l'application de Kuratowski plonge isométriquement $X$ dans $E(X,\omega)$, de telle façon que toute isométrie de $X$ s'etend (uniquement) en une isométrie de $E(X,\omega)$, et le morphisme d'extension est continu.
\end{proposition}
En partant d'un espace métrique polonais $X$; on construit par récurrence une
suite $(X_{i})_{i\in\mathbb{N}}$ d'espaces métriques définie par $X_0=X$ et  $X_{i+1} = E(X_{i},\omega)$ (Ceci a un sens,
puisque l'application de Kuratowski permet d'identifier $X_{i}$ à un sous-espace de $X_{i+1}$).\\
Les résultats précédents permettent de voir que toute isométrie de $X_{i}$ s'étend de manière unique en une isométrie de $X_{i+1}$, et que ceci définit un morphisme continu de $Iso(X_{i})$ dans $Iso(X_{i+1})$.\\
Par conséquent, si on pose $X_{\infty}=\underset{i=1}{\overset{\infty}{\bigcup}}X_i$, on a défini un morphisme continu $\Psi$ de $Iso(X)$ dans $Iso(X_{\infty})$, qui a de plus la propriété suivante :
Pour toute isométrie $\varphi \in Iso(X),\,\Psi(\varphi)$ est un prolongement de $\varphi$ (On dit encore comme dans \cite{us6}
que $X$ est $g$-plongé dans $X_{\infty}$).

\begin{theorem}(voir \cite{vp1})
Si $X$ est un espace métrique séparable, alors le complété de l'espace $X_{\infty}$ est isométrique à l'espace universel d'Urysohn.
\end{theorem}

\subsection{Le groupe des isométries}
Si $X$ est un espace métrique, alors la topologie de la convergence simple sur le groupe $Iso(X)$ des isométries de $X$ sur lui-même coïncide avec la topologie compact-ouvert et ces deux topologies sont compatibles avec la structure de groupe de $Iso(X)$ (\cite{vp1}). Une base de voisinages de l'identité pour cette topologie est constitué des ensembles de la forme $$V[F,\varepsilon]=\{g\in Iso(X):\,\forall\,x\in F,\,d(g(x),x)< \varepsilon\},$$ où $F\subseteq X$ est fini et $\varepsilon> 0$. Si $X$ est séparable et donc vérifie le deuxième axiome de dénombrabilité, alors $Iso(X)$ est aussi séparable. Si $X$ est séparable et complet, alors $Iso(X)$ est métrisable et complet. En effet, si $(x_n)_{n\in \mathbb{N}}$ est une partie dense de $X$, on vérifie que l'application $$d(f,g)=\underset{i=1}{\overset{\infty}{\sum}}2^{-i}[d_X(f(x_i),g(x_i))+d_X(f^{-1}(x_i),g^{-1}(x_i))]$$ est une métrique complète sur $Iso(X)$ qui engendre la topologie de la convergence simple (Voir \cite{vp1}, Proposition $5.2.1$)

\begin{proposition}(\cite{vp1})
Soient $X$ un espace métrique, $g\in Iso(X)$ et $f\in E(X,\omega)$ une fonction de Kat\v etov contrôlée par un sous-ensemble fini $A\subset X$, alors l'application
\[\begin{array}{lll}
^{g}f:& X &\longrightarrow \mathbb{R}\\
&x &\longmapsto (^{g}f)(x)=f(g^{-1}x)\\
\end{array}\]
est une fonction de Kat\v etov contrôllée par le sous-ensemble fini $g(A)$.
\end{proposition}
\begin{remark}
D'après ce qui précède, l'application
\[\begin{array}{lll}
\tau:& Iso(X)\times E(X,\omega) &\longrightarrow E(X,\omega)\\
&(g,f) &\longmapsto g.f=\,^{g}f\\
\end{array}\]
définie une action par isométries de $Iso(X)$ sur l'espace métrique $E(X,\omega)$. Cette action sera appelée action par translations à gauche.
\end{remark}
\begin{lemma}\cite{vp1}\label{lemact}
L'action par translations à gauche du groupe $Iso(X)$ dans $E(X,\omega)$ est continue. De plus, l'application de Kuratowski
\[\begin{array}{lll}
\delta:& X &\longrightarrow E(X,\omega)\\
&x &\longmapsto \delta_{x}\\
\end{array}\]
 est équivariante.
\end{lemma}
\begin{remark}
Comme conséquence du lemme \ref{lemact}, le groupe $Iso(X)$ opère continûment par isométries sur chaque extension de Ketetov $X_{i+1} = E(X_{i},\omega)$. Ceci nous permet d'avoir une action continue par isométries de $G$ sur $X_{\infty}$ qui se prolonge uniquement en une action continue de $G$ sur le complété $\mathbb{U}$ de $X_{\infty}$.
\end{remark}
\begin{theorem}(Uspenskij \cite{us5})\label{theousp} Le groupe topologique $Iso(\mathbb{U})$ est universel pour la classe des groupes topologiques de poids dénombrables. Autrement dit, tout groupe topologique de poids dénombrable est isomorphe à un sous-groupe de $Iso(\mathbb{U})$.
\end{theorem}
\begin{prof}
Soit $G$ un groupe topologique de poids dénombrable. Par le théorème de Teleman(Théorème \ref{tel}), on peut choisir un espace métrique séparable $X$ tel que $G$ est isomorphe à un sous-groupe de $Iso(X)$. D'après ce qui précède et la construction de Kat\v etov de l'espace d'Urysohn $\mathbb{U}$, le groupe $Iso(X)$ est isomorphe à un sous-groupe de $Iso(\mathbb{U})$
\end{prof}
\begin{remark}
Le théorème d'Uspenkij précédent a été récemment raffiné par Melleray (\cite{mel}) comme suit:
\end{remark}
\begin{theorem}(Melleray\cite{mel})
Tout groupe polonais est isomorphe au sous-groupe de $Iso(\mathbb{U})$ constitué par les isométries qui laissent fixe un certain fermé $F\subseteq \mathbb{U}$.
\end{theorem}

\begin{remark}
Avec les deux solutions ($Iso(\mathbb{U})$ et $Homeo(I^{\aleph_{0}})$) fourni par Uspenskij du problème d'existence des groupes topologiques universels de poids dénombrables, la question naturelle qui apparait dans \cite{us5} est celle de savoir si les deux solutions ci-dessus sont isomorphes en tant que groupes topologiques. Une réponse à cette question est apporté par le résultat suivant:
\end{remark}

\begin{theorem}(Pestov \cite{vp4})
Le groupe $Iso(\mathbb{U})$ muni de la topologie de la convergence simple est extrêmement moyennable.
\end{theorem}
\begin{corollary}
Les groupes topologiques $Iso(\mathbb{U})$ et $Homeo(I^{\aleph_{0}})$ ne sont pas isomorphes.
\end{corollary}
\begin{prof}
Le groupe $Homeo(I^{\aleph_{0}})$ possède une action continue sur $I^{\aleph_{0}}$ sans point fixes car le cube de Hilbert $I^{\aleph_{0}}$ est homogène (\cite{keller}): Pour tous $x,y\in I^{\aleph_{0}}$ distints, il existe $f\in Homeo(I^{\aleph_{0}})$ tel que $f(x)=y$. Ce groupe n'est donc pas extrêmement moyennable.
\end{prof}
\begin{remark}
 Le théorème d'Uspenkij présenté dans le paragraphe précédent sur l'universalité du groupe $Iso(\mathbb{U})$ pour la classe des groupes topologiques de poids dénombrable suggère de conjecturer comme beaucoup de mathématiciens que le groupe $Iso(\mathbb{U}_{\mathfrak m})$ où $U_{\mathfrak m}$ est l'espace généralisé d'Urysohn $\mathfrak m$-homogène et $\mathfrak m$-universel pour un cardinal infini non dénombrable $\mathfrak m$ est universel pour la classe des groupes topologiques de poids $\mathfrak m$.\\
Dans cette section, nous montrerons que cette conjecture est fausse.
\end{remark}
 \section{Construction de l'espace d'Urysohn-Kat\v etov $\mathbb{U}_{\mathfrak m}$}\index{espace!d'Urysohn-Kat\v etov}
Soit $X$ un espace topologique. Rappelons que le poids\index{espace topologique!poids d'un} de $X$ est le nombre cardinal $\omega(X)=\min\{|\mathcal{B}|:\,\,\mathcal{B}\, \text{est une base de}\, $X$\}$ et la densité\index{espace topologique!densité d'un} de $X$ est le nombre cardinal $d(X)=\min\{|A|:\,\,A\, \text{est une partie dense de}\, $X$\}$. Où $|A|$ désigne la cardinalité de $A$.\\
Soit $\mathfrak m$ un cardinal infini non dénombrable vérifiant
\begin{equation}
\label{eq:mn}
\sup\left\{{\mathfrak m}^{\mathfrak n}\colon {\mathfrak n}<{\mathfrak m}\right\}={\mathfrak m}.
\end{equation}
Nous allons présenter dans ce paragraphe en suivant \cite{kat} la construction de Kat\v etov de l'espace universel d'Urysohn de poids $\mathfrak m$. La construction dans le cas particulier où $\mathfrak m=\aleph_{0}$ a déjà été présentée au début de ce chapitre. Nous supposerons dans la suite que $\mathfrak m$ est un cardinal infini non dénombrable vérifiant la condition (\ref{eq:mn}).
\begin{definition}
Soit $\tau$ un cardinal infini. Un espace métrique $(X,d)$ est dit:
\begin{enumerate}
  \item $\tau$-homogène si pour tout $A,B\subseteq X$ avec $card(A)< \tau$ et $card(B)< \tau$, et pour toute isométrie $f:A\longrightarrow B$, il existe une isométrie $g:X\longrightarrow X$ telle que $g|_{A}=f$.
  \item $\tau$-universel si pour tout espace métrique $Y$ tel que $card(Y)\leq \tau$, il existe un plongement isométrique $i:Y\longrightarrow X$.
  \item fortement $\tau$-universel si pour tout espace métrique $Y$ tel que $\omega(Y)\leq \tau$, il existe un plongement isométrique $i:Y\longrightarrow X$.
      \item $\tau$-Urysohn universel (ou Urysohn universel de poids $\tau$) si $(X,d)$ est de poids $\tau,\,\tau$-homogène et fortement $\tau$-universel.
\end{enumerate}
\end{definition}
Le théorème suivant est démontré dans \cite{kat}
\begin{theorem}(Kat\v etov \cite{kat})
Si $\mathfrak m$ est un cardinal infini vérifiant (\ref{eq:mn}), alors il existe à isométrie près un unique espace universel d'Urysohn de poids $\mathfrak m$.
\end{theorem}
Si $\mathfrak{m}$ est un nombre cardinal, notons $E_{\mathfrak m}(X)=\{f\in E(X):\,\,|\text{support de f}|< \mathfrak m\}$ le sous-espace métrique de $E(X)$ constitué des fonctions de Kat\v etov donc le support a une cardinalité $< \mathfrak m$ (ou de façon équivalente de densité $< \mathfrak{m}$). La densité de $E_{\mathfrak m}(X)$ ne dépasse pas $\sup \{d(X)^{\mathfrak n}:\,\,\mathfrak n< \mathfrak m\}$, où $d(X)$ désigne la densité de $X$.\\
 L'application de Kuratowski plonge isométriquement $X$ dans $E_{1}(X)$ et de façon canonique dans tout espace $E_{\mathfrak m}(X)$ (\cite{kat}).\label{kat2}\\
Partons d'un espace métrique $X$ tel que $\omega(X)=\mathfrak m$. Considérons la chaîne des extensions de Kat\v etov $E^{\tau+1}_{\mathfrak m}(X)=E_{\mathfrak m}(E^{\tau}_{\mathfrak m}(X))$ où $\tau < \mathfrak{m}$ et $E^{\tau}_{\mathfrak m}(X)=\underset{\lambda<\tau}{\overset{}{\bigcup}}E^{\lambda}_{\mathfrak m}(X)$ pour les cardinaux limites. Sous la condition ($\ref{eq:mn}$), la densité de chaque extension itérée de Kat\v etov est bornée par $\mathfrak m$. Ainsi, en itérant cette construction $\mathfrak m$ fois et en prenant le complété de la réunion \begin{equation}
\label{eq:chain}
\bigcup_{{\mathfrak n}<{\mathfrak m}} E_{\mathfrak m}^{\mathfrak n}(X),\end{equation}
on obtient un plongement de $X$ dans un espace métrique complet $\mathbb{U}_{\mathfrak{m}}$ de densité $\mathfrak{m}$ (\cite{kat}).\\

\begin{remark}
\begin{enumerate}
  \item Si un cardinal non-dénombrable $\mathfrak{m}$ vérifiant la condition (\ref{eq:mn})est tel que $\mathfrak{m}\geq \mathfrak{c}$ ($\mathfrak{c}$ est le continu), alors la réunion (\ref{eq:chain}) est automatiquement complète (\cite{kat}).
  \item Un espace métrique $\mathbb{U}_{\mathfrak{m}}$ avec les propriétés précédentes est unique (à isométrie près) et $ \mathfrak{m}$-homogène (\cite{kat}).
  \item Par construction, l'espace métrique $\mathbb{U}_{\mathfrak{m}}$ est universel pour les espaces métriques de poids $\leq \mathfrak{m}$.
  \item Si $m<\sup\left\{{\mathfrak m}^{\mathfrak n}\colon {\mathfrak n}<{\mathfrak m}\right\}$, alors il n'existe pas d'espace métrique $\mathfrak{m}$-homogène et $\mathfrak{m}$-universel de poids $\mathfrak{m}$.
\end{enumerate}
\end{remark}

L'espace métrique $\mathbb{U}_{\mathfrak{m}}$ vérifie en plus la propriété suivante:
\begin{lemma}\label{extkat} Si $Y$ est un sous-espace de $\mathbb{U}_{\mathfrak{m}}$ de cardinalité $< \mathfrak{m}$ et $f$ une fonction de Kat\v etov sur $Y$, alors il existe un point $z\in \mathbb{U}_{\mathfrak{m}}$ tel que $f(y)=d(z,y)$ pour tout $y\in Y$.
\end{lemma}
\begin{prof}
Soit $Y\subset \mathbb{U}_{\mathfrak{m}}$ un sous-espace de cardinalité $<\mathfrak m$ de $\mathbb{U}_{\mathfrak{m}}$ et soit $f\in E(Y)$. L'espace $\mathbb{U}_{\mathfrak{m}}$ étant universel de poids $\mathfrak m$, l'espace métrique $Y_{f}=Y\cup\{f\}$ se plonge isométriquement dans $\mathbb{U}_{\mathfrak{m}}$. Donc il existe une copie isométrique $Y^{\prime}_{f}=Y^{\prime}\cup\{z\}$ de $Y_f$ contenue dans $\mathbb{U}_{\mathfrak{m}}$. Par définition de $Y_f$, il existe une isométrie $\varphi:Y\longrightarrow Y^{\prime}$ telle que $d(z,\varphi(y))=f(y)$ pour tout $y\in Y$. L'espace $\mathbb{U}_{\mathfrak{m}}$ étant $\mathfrak m$-homogène, l'isométrie $\varphi$ se prolonge en une isométrie $\widetilde{\varphi}:V\longrightarrow \mathbb{U}_{\mathfrak{m}}$. Soit $y\in Y$, on a: $d(\widetilde{\varphi}^{-1}(z),y)=d(z,\widetilde{\varphi}(y))=d(z,\varphi(y))=f(y)$.
\end{prof}

\section{Voisinages de l'identité des Sous-groupes de $Iso(\mathbb{U}_{\mathfrak m})$}
Le groupe $Iso(\mathbb{U}_{\mathfrak m})$ de toutes les isométries de $\mathbb{U}_{\mathfrak m}$ sur lui-même est muni de la topologie de la convergence simple sur $\mathbb{U}_{\mathfrak m}$. Une base de voisinage de l'élément neutre est constitué des ensembles de la forme
$$V[x_{1},x_{2},...,x_{n};\varepsilon]=\{g\in Iso(\mathbb{U}_{\mathfrak m}):\,d(x_{i},gx_{i})<\varepsilon,\,i=1,2,...,n\},$$
où $\{x_{1},x_{2},...,x_{n}\}$ est un sous-ensemble fini de $\mathbb{U}_{\mathfrak m}$ et $\varepsilon>0$.\\
Le lemme suivant est fondamental pour la suite.
\begin{lemma}\label{lemmatopo}
Soit $\mathfrak m$ un cardinal infini vérifiant l'égalité (\ref{eq:mn}), et soit $G$ un sous-groupe de $Iso(\mathbb{U}_{\mathfrak m})$ de densité $<\mathfrak m$. Les ensembles $$V[x;\varepsilon]\cap G,\,\,\,x\in \mathbb{U}_{\mathfrak m},\,\,\,\varepsilon>0$$ forment une base de voisinage de l'identité pour la topologie de $G$.
\end{lemma}
\begin{prof}
Sans nuire à la généralité, nous pouvons remplacer $G$ par un sous-groupe dense de cardinalité $<\mathfrak m$.\\
Soit $X=\{x_{1},x_{2},...,x_{n}\}$ un sous-ensemble fini de $\mathbb{U}_{\mathfrak m}$ et $\varepsilon>0$ arbitraire.\\
 Nous cherchons $y\in \mathbb{U}_{\mathfrak m}$ et $\gamma>0$ tel que $$V[y;\gamma]\cap G\subseteq V[x_{1},x_{2},...,x_{n};\varepsilon].$$
 Notons $D$ le diamètre de $X$. Choisissons $\gamma>0$ tel que $\gamma\leq \varepsilon$ et les boules centrées en $x_{i},\,i=1,2,..,n$ et de rayons $n\gamma$ sont deux à deux disjointes. La fonction
 $$f(x_i)=D+i\gamma$$ est de Kat\v etov. En effet,  Si $x_{i},x_{j}\in X$, alors on a d'une part

  $$|f(x_i)-f(x_j)|=|i-j|\gamma \leq (n-1)\gamma < 2n \gamma < d(x_i,x_j),$$
 car $$ B_{n\gamma}(x_i)\cap B_{n\gamma}(x_j)=\emptyset  \Longrightarrow d(x_j,x_i)> 2n\gamma,$$
 et d'autre part $$f(x_i)+f(x_j)=2D+(i+j)\gamma > D \geq d(x_i,x_j).$$
 Notons encore par $f$ son extension de Kat\v etov $k_{X}(f)$ sur $\mathbb{U}_{\mathfrak m}$.\\
 Soit $i=1,2,...,n$. Supposons $x\in \mathbb{U}_{\mathfrak m}$ et pour tout $j=1,2,...,i,\,\,d(x_{j},x)\geq (i-j)\gamma$. Alors pour tout $j=1,2,...,n$, on a: $$d(x,x_j)+f(x_j)\geq (i-j)\gamma +D+j\gamma.$$ Donc $$d(x,x_j)+f(x_j)\geq D+i\gamma.$$  Ainsi, $f(x)\geq D+i\gamma$. Donc $f$ vérifie la propriété suivante:

 \begin{equation}
\label{eq:prop}
\forall x\in \mathbb{U}_{\mathfrak m},\,\,\,(f(x)< D+i\gamma)\Longrightarrow x\in \underset{j=1}{\overset{i-1}{\bigcup}}B_{(i-j)\gamma}(x_{j})
\end{equation}
Notons $A$ la réunion de toutes les $G$-orbites des points $x_{1},x_{2},...,x_{n}$. Par le lemme \ref{extkat}, il existe $y\in \mathbb{U}_{\mathfrak m}$ tel que $f(a)=d(y,a)$ pour tout $a\in A$. Montrons que $$V[y;\gamma]\subseteq V[x_{1},x_{2},...,x_{n};\varepsilon].$$
Soit $g\in V[y;\gamma]$ i.e $d(y,gy)<\gamma$. Pour tout $x\in A$, on a:

$$|f(x)-f(gx)|=|d(y,x)-d(y,gx)|=|d(gy,gx)-d(y,gx)|\leq d(y,gy)< \gamma$$

En particulier, pour $i=1,2,...,n$, nous obtenons: $$f(gx_{i})< D+(i+1)\gamma.$$
Par conséquent, nous obtenons via (\ref{eq:prop})
$$gx_{i}\in \underset{j=1}{\overset{i}{\bigcup}}B_{(i-j+1)\gamma}(x_{j}).$$

Ainsi, pour tout $i\in \{1,2,...,n\}$, il existe $j=\phi(i)\in \{1,2,...,i\}$ tel que $$gx_{i}\in B_{(i-j+1)\gamma}(x_{j}).$$

Montrons que $\phi$ est injective. \\

Soit $i\neq k$, alors il existe d'une part $j=\phi(i)\in \{1,2,...,i\}$ tel que $$gx_{i}\in B_{(i-j+1)\gamma}(x_{j}).$$ D'autre part il existe $l=\phi(k)\in \{1,2,...,k\}$ tel que $$gx_{k}\in B_{(i-l+1)\gamma}(x_{l}).$$

Supposons $j=l$, alors on a: $$d(x_{i},x_{k})=d(gx_{i},gx_{k})<d(gx_{i},x_{j})+d(gx_{k},x_{j})<2(i-j+1)\gamma.$$

Ainsi, $d(x_{i},x_{k})< 2(i-j+1)\gamma$. Donc $B_{(i-j+1)\gamma}(x_i)\cap B_{(i-j+1)\gamma}(x_k)\neq\emptyset$. Ce qui est absurde par le choix de $\gamma$. Ainsi $\phi$ est une application injective de $\{1,2,...,n\}$ dans lui-même qui vérifie $\phi(i)\leq i$. Donc $\phi$ est l'application identique. Ceci implique $gx_{i}\in B_{\varepsilon}(x_i)$ pour tout $i$ comme souhaité.

\end{prof}

\section{Groupes SIN et Groupes FSIN}

\begin{definition}
Soit $G$ un groupe topologique. Il existe sur $G$ deux structures uniformes bien connues: la structure uniforme droite que nous noterons $\mathcal{U}_{R}(G)$ et la structure uniforme gauche que nous noterons $\mathcal{U}_{L}(G)$.
\begin{enumerate}
  \item Soit $V$ un voisinage de $e$ dans $G$. Posons $V_{d}=\{(x,y):\,\,\,xy^{-1}\in V\}$.\\
  La famille des ensembles $V_{d}$ lorsque $V$ parcourt une base de voisinages de $e$ forme une base d'entourages pour la structure uniforme droite $\mathcal{U}_{R}(G)$ sur $G$.
  \item Soit $V$ un voisinage de $e$ dans $G$. Posons $V_{g}=\{(x,y):\,\,\,x^{-1}y\in V\}$.\\
  La famille des ensembles $V_{g}$ lorsque $V$ parcourt une base de voisinages de $e$ forme une base d'entourages pour la structure uniforme gauche $\mathcal{U}_{L}(G)$ sur $G$.
\end{enumerate}
\end{definition}

\begin{theorem}\label{theosin}
Soit $G$ un groupe toplogique. Les propositions suivantes sont équivalentes:
\begin{enumerate}
  \item Les structures uniformes $\mathcal{U}_{R}(G)$ et $\mathcal{U}_{L}(G)$ coïncident.
  \item Pour tout voisinage $U$ de $e$, il existe un voisinage $V$ de $e$ tel que $x^{-1}Vx\subseteq U$ pour tout $x\in G$.
  \item L'inversion sur $G$ est uniformément continue pour la paire $(\mathcal{U}_{R}(G),\mathcal{U}_{R}(G))$ ou pour la paire $(\mathcal{U}_{L}(G),\mathcal{U}_{L}(G))$.
\end{enumerate}
\end{theorem}
\begin{prof}
\begin{enumerate}
  \item [$1)\Longrightarrow 2)$] Supposons que les structures uniformes $\mathcal{U}_{R}(G)$ et $\mathcal{U}_{L}(G)$ coïncident. Soit
  $U$ un voisinage de l'identité dans $G$. Alors il existe un voisinage symétrique $V$ de l'identité dans $G$ tel que $V_{d}\subseteq U_{g}$. Ceci implique que si $xy^{-1}\in V$ alors $x^{-1}y\in U$. Ainsi, si $y\in Vx$, alors $y\in xU$. Donc pour tout $x\in G$, on a: $Vx\subseteq xU$. Ainsi $x^{-1}Vx\subseteq U$.
  \item [$2)\Longrightarrow 3)$]Soit $A\in\mathcal{U}_{R}(G)$. Alors, il existe un voisinage de l'identité $U$ tel que $U_{d}\subseteq A$. Puisque $U^{-1}$ est également un voisinage de l'identité, alors par $2$, il existe un voisinage $V$ de l'identité tel que $x^{-1}Vx\subseteq U^{-1}$ pour tout $x\in G$. Ainsi, si $xy^{-1}\in V$ alors $x^{-1}(xy^{-1})x=y^{-1}x\in U^{-1}$. Donc $(y^{-1}x)^{-1}=x^{-1}(y^{-1})^{-1}\in U$ et ainsi $(x^{-1},y^{-1})\in U_{d}$. Donc $V_{d}\subseteq \{(x,y)\in X\times X:\,(x^{-1},y^{-1})\in U_{d}\}\subseteq A$. Ainsi, l'inversion est uniformément continue pour la paire $(\mathcal{U}_{R}(G),\mathcal{U}_{R}(G))$.
  \item [$3)\Longrightarrow 1)$]Soit $A\in \mathcal{U}_{L}(G))$. Alors il existe un voisinage $V$ de l'identité dans $G$ tel que $V_{g}\subseteq A$. Notons que $(x,y)\in V_{g}$ si et seulement si $(x^{-1},y^{-1})\in V_{d}$. Donc $V_{g}=\{(x,y):\,(x^{-1},y^{-1})\in V_{d}\}$. Par $3$; il existe un voisinage de l'identité $U$ tel que $U_{d}\subseteq V_{g}$. Ainsi $A\in \mathcal{U}_{R}(G)$. De la même façon, si $A\in \mathcal{U}_{R}(G))$, alors il existe un voisinage de l'identité $U$ tel que $U_{g}\subseteq A$.
\end{enumerate}
\end{prof}
\begin{definition}
Un groupe topologique $G$ est dit SIN(Small Invariant Neighbourhood)\index{groupe!SIN} s'il vérifie l'une des conditions équivalentes du théorème \ref{theosin}.
\end{definition}
\begin{remark}
Si $G$ est un groupe SIN, alors pour tout voisinage $U$ de l'élément neutre, il existe un voisinage $V$ de l'élément neutre tel que $x^{-1}Vx\subseteq U$ pour tout $x\in G$. Posons $W=\underset{x\in G}{\overset{}{\bigcup}}x^{-1}Vx$. Alors $W$ est un voisinage de l'identité dans $G$ contenu dans $U$. De plus, pour tout $g\in G$, on a:
$$g^{-1}Wg=g^{-1}(\underset{x\in G}{\overset{}{\bigcup}}x^{-1}Vx)g=\underset{x\in G}{\overset{}{\bigcup}}(xg)^{-1}Vxg=W.$$
Donc $W$ est invariant.
\end{remark}
Ainsi nous avons le théorème suivant qui justifie la terminologie SIN(Small Invariant Neibourhoods)
\begin{theorem}
Un groupe topologique $G$ est SIN si et seulement s'il possède une base de voisinages l'identité constitué des voisinages invariants.
\end{theorem}
\begin{example}
\begin{enumerate}
\item Tout groupe  précompact est SIN.\\
En effet, soit $V$ un voisinage symétrique de l'élément neutre dans $G$. Soit $U$ un voisinage symétrique de l'identité dans $G$ tel que $U^{3}\subseteq V$. Puisque $G$ est précompact, il existe $A\subseteq G$ fini tel que $G=AU$. Notons que $W=\underset{a\in A}{\overset{}{\bigcap}}a^{-1}Ua$ est un voisinage de l'identité dans $G$ et $aWa^{-1}\subseteq U$.
Soit $g\in G$, alors il existe $u\in U$ tel que $g=au$. Nous avons:
$$g^{-1}Wg=(au)^{-1}(aWa^{-1})(au)=u^{-1}a^{-1}(aWa^{-1})au\subseteq u^{-1}Uu\subseteq U^{3}\subseteq V.$$
\item Tout groupe abélien est SIN. En effet, tout sous-ensemble d'un groupe abélien est invariant.
 \item Le groupe des permutations d'un ensemble infini n'est pas SIN.
\end{enumerate}
\end{example}
En général, on a la proposition suivante:
\begin{proposition}\label{propsinfini}
Soit $X$ un ensemble non vide. Le groupe $S_{X}$ de toutes les bijections de $X$ sur lui-même muni de la topologie de la convergence simple sur $X$ où $X$ est muni de la topologie discrète est SIN si et seulement si $|X|<\infty$.
\end{proposition}
Rappelons la description de la topologie de $S_{X}$. Pour tout sous-ensemble fini $M=\{m_{1},...,m_{n}\}$ de $X$, posons $St_{M}=\{g\in S_{X}:\,\,\forall j=1,...,n,\,g(m_j)=m_j\}$. La famille des ensembles de la forme $St_{M}$ où $M$ est un sous-ensemble fini de $X$ forme une base de voisinages de l'identité pour la topologie de la convergence simple sur $S_{X}$.
\begin{lemma}\label{lemsinfini}
Le groupe $S_{X}$ est non-discret si et seulement si $X$ est infini.
\end{lemma}
\begin{prof}
Soit $\iota$ l'application identique de $X$ sur lui-même i.e l'identité de $S_X$. Supposons que $|X|<\infty$, alors le sous-groupe $St_{X}$ est un élément de la base. Comme $St_{X}=\{\iota\}$, alors $\{\iota\}$ est un ouvert. Ainsi, $S_{X}$ est discret.\\
 Réciproquement, si $X$ est infini, alors il n'existe pas de sous-ensemble fini $M\subset X$ tel que $St_{M}=\{\iota\}$. Ainsi $\{\iota\}$ n'est pas ouvert. Donc $S_{X}$ est non-discret.
\end{prof}\\
\begin{prof}(\textbf{Proposition \ref{propsinfini}})\\
Supposons $|X|<\infty$. Par le lemme \ref{lemsinfini}, $St_{X}=\{\iota\}$ est un ouvert. Il s'en suit que pour tout sous-ensemble fini $M$ de $X$, on a: $g^{-1}St_{X}g\subseteq St_{M}$ pour tout $g\in S_X$. Donc $S_X$ est SIN.\\
Réciproquement, soit $M$ un sous-ensemble fini de $X$. Par hypothèse, il existe un sous-ensemble fini $N$ de $X$ tel que $g^{-1}St_{N}g\subseteq St_{M}$ pour tout $g\in S_X$. Nous allons montrer que $M\cup N=X$. Supposons le contraire i.e $X\setminus(M\cup N)\neq \emptyset$. Soit $m$ un point fixé de $M$ et soit $g\in S_{X}$ tel que $g(m)\in X\setminus(M\cup N)$. Maintenant, soit $f\in St_{N}$ tel que $f(g(m))\neq g(m)$. Alors $g^{-1}fg(m)\neq m$. Ainsi, $g^{-1}St_{N}g\nsubseteq St_{M}$. Ce qui est une contradiction. Donc $M\cup N=X$ et $X$ est fini.
\end{prof}
\begin{definition}
Un groupe topologique $G$ est FSIN\index{groupe!FSIN} si l'ensemble des fonctions sur $G$ à valeurs réelles uniformément continues à gauche coïncide avec l'ensemble des fonctions à valeurs réelles sur $G$ uniformément continue à droite.
\end{definition}
 \begin{remark}
Il est clair que tout groupe topologique SIN est FSIN. La question qui consiste à étudier la réciproque de cette assertion est appelée problème de Itzkowitz\index{problème d'Itzkowitz} qui fut le premier à s'intéresser à cette question dans \cite{itz}. Cette question a été répondu par l'affirmative entre autre pour la classe des groupes localement compacts (\cite{itz3}), celle des groupes métrisables (\cite{prota2}) et la classe des groupes localement connexes (\cite{megrepest}). Cependant, cette question reste ouverte dans le cas général. Les dévellopements récents sur cette question, peuvent être consulter dans (\cite{tro}).
\end{remark}
Le théorème suivant est une caractérisation très utile des groupes FSIN que nous utiliserons dans la suite.
\begin{theorem}(Protasov et Saryev \cite{pro})\label{prota}\index{théorème!de Protasov et Saryev}
Soit $G$ un groupe topologique et soit $\mathcal{U}$ une base de voisinages symétriques de $e$. Les propositions sont équivalentes:
\begin{enumerate}
  \item $G$ est FSIN
  \item Pour tout $A\subset G$ et pour tout $U\in \mathcal{U}$, il existe $V\in \mathcal{U}$ tel que $VA\subset AU$.
\end{enumerate}
\end{theorem}

Avant d'énoncer un important corollaire de ce théorème, introduisons une terminologie empruntée de \cite{megrepest}:
\begin{definition}
Un sous ensemble $A$ d'un groupe topologique $G$ est dit neutre à gauche si pour tout voisinage $V$ de l'identité dans $G$, il existe un voisinage $U$ de l'identité tel que $UA\subseteq AV$. De façon similaire, on définit un sous-ensemble neutre à droite. Un sous-ensemble à la fois neutre à gauche et à droite est dit neutre.
\end{definition}
\begin{example}
\begin{enumerate}
  \item Tout sous-espace compact $A$ d'un groupe topologique $G$ est neutre.\\
En effet, soit $V$ un voisinage de l'identité de $G$. Alors il existe un voisinage $W$ de l'identité tel que $W^{2}\subseteq V$. Les ensembles $\{aW:\,\,a\in A\}$ et $\{Wa:\,\,a\in A\}$ sont tous les recouvrement ouverts de $A$. Puisque $A$ est compact, il existe un nombre fini de points $a_{1},...,a_{n}\in A$ et  $b_{1},...,b_{m}\in A$ tel que $A\subseteq \underset{i=1}{\overset{n}{\bigcup}}a_{i}W$ et $A\subseteq \underset{j=1}{\overset{m}{\bigcup}}Wb_{j}$. Posons $U=\underset{i=1}{\overset{n}{\bigcap}}a_{i}Wa^{-1}_{i}$ et $U^{\prime}=\underset{j=1}{\overset{m}{\bigcap}}b^{-1}_{j}Wb_{j}$. Clairement, $U$ et $U^{\prime}$ sont des voisinages de l'identité. De plus, $Ua_{i}\subseteq a_{i}W$ et $b_{j}U^{\prime}\subseteq Wb_{j}$, pour $1\leq i\leq n$ et $1\leq j \leq m$ respectivement. Maintenant
$$UA\subseteq U(\underset{i=1}{\overset{n}{\bigcup}}a_{i}W)=\underset{i=1}{\overset{n}{\bigcup}}(Ua_{i})W\subseteq \underset{i=1}{\overset{n}{\bigcup}}a_{i}W^{2}\subseteq \underset{i=1}{\overset{n}{\bigcup}}a_{i}V\subseteq \underset{a\in A}{\overset{}{\bigcup}}aV=AV$$ et
$$AU^{\prime}\subseteq (\underset{j=1}{\overset{m}{\bigcup}}Wb_{j})U^{\prime}=\underset{j=1}{\overset{m}{\bigcup}}W(b_{j}U^{\prime})\subseteq \underset{j=1}{\overset{m}{\bigcup}}W^{2}b_{j}\subseteq \underset{j=1}{\overset{m}{\bigcup}}Vb_{j}\subseteq \underset{a\in A}{\overset{}{\bigcup}}Va=VA.$$
Ce qui démontre que $A$ est neutre.
  \item Un sous-groupe normal $H$ d'un groupe topologique $G$ est neutre. En effet, soit $V$ un voisinage de l'élément neutre dans $G$. Puisque $H$ est normal, nous avons $vH=Hv$ pour tout $v\in V$. Puisque $VH=\underset{v\in V}{\overset{}{\bigcup}}vH$ et $HV=\underset{v\in V}{\overset{}{\bigcup}}Hv$, il s'en suit que $VH=HV$ et ainsi, $H$ est neutre.
\end{enumerate}

\end{example}

\begin{definition}
Un sous ensemble $A$ d'un groupe topologique $G$ est dit uniformément discret à gauche s'il est uniformément discret par rapport à la structure uniforme gauche. En d'autres termes, s'il existe un voisinage $V$ de l'identité tel que $(aV)\cap (bV)=\emptyset$ dès que $a,b\in A$ et $a\neq b$
\end{definition}

\begin{corollary}
Un groupe topologique $G$ est FSIN si et seulement si pour tout sous-ensemble uniformément discret à gauche $A\subseteq G$ et pour tout voisinage de l'identité $V$, il existe un voisinage $U$ de l'identité tel que $UA\subseteq AV$.
\end{corollary}
\begin{prof}
Il suffit d'établir la suffisance. Nous allons le faire en utilisant le théorème de Protasov-Saryev (théorème \ref{prota}). Soit $A$ un sous-ensemble arbitraire de $G$ et soit $V$ un voisinage de l'identité dans $G$. Choisissons un voisinage symétrique de l'identité $W$ tel que $W^{4}\subseteq V$. Soit $B$ un sous-ensemble maximal de $AW$ vérifiant la propriété: pour tout $a,b\in B,\,\,a\neq b$, les ensembles $aW$ et $bW$ sont disjoints. L'existence de B est garantie par une application convenable du lemme de Zorn. Alors, il est clair que $A\subseteq BW^{2}$. Puisque $B$ est uniformément discret à gauche, il existe un voisinage $U$ de l'identité avec la propriété que $UB\subseteq BW$. Nous avons:
$$UA\subseteq UBW^{2}\subseteq BW^{3}\subseteq AW^{4}\subseteq AV.$$
\end{prof}
\begin{remark}\label{rem2}
Si $G$ est un groupe topologique FSIN et dense dans un groupe $H$, alors $H$ est FSIN.\\
En effet, soit $f$ une fonction uniformément continue à gauche sur $H$. Alors la restriction de $f$ à $G$ est uniformément continue à gauche. Puisque $G$ est FSIN, $f$ est uniformément continue à droite. Ainsi $f$ est uniformément continue à droite sur $H$ grâce à l'observation générale suivante:
\begin{proposition}(\cite{bourb})
Soit $f$ une fonction continue sur un espace uniforme $X$ et $Y$ un sous-espace dense de $X$. Si la restriction de $f$ à $Y$ est uniformément continue ($Y$ muni de l'uniformité induite), alors $f$ est uniformément continue sur $X$.
\end{proposition}
\end{remark}
Pour terminer cette section rappelons les résultats suivants qui donnent des réponses positives au problème d'Itzkowitz pour la classe des groupes métrisables et celles des groupes localement connexe.
\begin{theorem}(Protasov \cite{prota2})\label{protasin}
Tout groupe métrisable FSIN est SIN.
\end{theorem}

\begin{theorem}(M. Megrelishvili, P. Nickolas et V. Pestov \cite{megrepest})\label{megrepest}
Un groupe topologique $G$ localement connexe est SIN si et seulement s'il est FSIN.
\end{theorem}

\section{Sous-groupes FSIN du groupe $Iso(\mathbb{U}_{\mathfrak m})$}
Commençons par rappeler la notion de propriété (OB) introduite dans \cite{ros}.\index{propriété (OB)}
\begin{definition}
 \begin{enumerate}
   \item Une pseudo-métrique sur un ensemble non vide $X$ est une fonction $\rho:X\times X\longrightarrow \mathbb{R}^{+}$ telle que pour tous $x,y,z\in X$, on a:
       \begin{enumerate}
         \item $\rho(x,x)=0$
         \item $\rho(x,y)=\rho(y,x)$
         \item $\rho(x,y)\leq \rho(x,z)+\rho(z,y)$
       \end{enumerate}
   \item Une pseudo-métrique $d$ sur un groupe topologique $G$ est dite invariante à gauche si $d(ax,ay)=d(x,y)$ pour tout $a\in G$
 \end{enumerate}

\end{definition}
\begin{definition}
Un groupe topologique $G$ possède la propriété (OB) si toute pseudo-métrique continue invariante à gauche sur $G$ est borné.
\end{definition}
L'ensemble des groupes possédant la propriété (OB) comprend:
\begin{enumerate}
  \item Le groupe symétrique infini $S_{\infty}$ muni de sa topologie polonaise (\cite{berg}).
  \item Le groupe $Iso(\mathbb{U})$ des isométries de l'espace polonais universel d'Urysohn $\mathbb{U}$ sur lui-même muni de la topologie de la convergence simple (\cite{ros}).
  \item Le groupe $Homeo([0,1]^{\aleph_{0}})$ des homéomorphismes du cube de Hilbert sur lui-même muni de la topologie de la convergence uniforme (\cite{ros}).
\end{enumerate}

Nous pouvons à présent établir le principal résultat de cette section:
\begin{theorem}
Soit $\mathfrak m$ un cardinal infini non dénombrable vérifiant la condition (\ref{eq:mn}). Si $G$ est un sous-groupe topologique de $Iso(\mathbb{U}_{\mathfrak m})$ de densité $<\mathfrak m$ et possédant la propriété (OB), alors $G$ est FSIN.
\end{theorem}
\begin{prof}
En vertu de la Remarque \ref{rem2}, nous pouvons supposons sans perte de généralité que la cardinalité de $G$ ne dépasse pas $\mathfrak m$. Soit $A\subseteq G$ un sous-ensemble uniformément discret à gauche et soit $V$ un voisinage symétrique de $e$. Nous allons chercher un voisinage $W$ de $e$ tel que $WA\subseteq AV$.\\
 Comme $A$ est uniformément discret à gauche, on a:
\begin{equation}
\label{eq:ab}
aV\cap bV=\emptyset
\end{equation}
dès que $a,b\in A,\,\,a\neq b$.\\
 En utilisant le Lemma \ref{lemmatopo}, choisissons $x\in \mathbb{U}_{\mathfrak m}$ et $\varepsilon>0$ tel que $$G\cap V[x;2\varepsilon]\subseteq V.$$
Soit $(Ax)_{\varepsilon}$ le $\varepsilon$-voisinage ouvert de $Ax$ dans $\mathbb{U}_{\mathfrak m}$: $$(Ax)_{\varepsilon} = \bigcup_{a\in A}B_{\varepsilon}(ax).$$
Si $a\neq b$, alors les $\varepsilon$-boules ouvertes centrées en $ax$ et $bx$ de $Ax$ sont disjointes.\\

En effet, en supposant $B_{\varepsilon}(ax)\cap B_{\varepsilon}(bx)\neq \emptyset$, il existe $u\in \mathbb{U}_{\mathfrak m}$ tel que $d(u,ax)<\varepsilon$ et $d(u,ax)<\varepsilon$. Par application de l'inégalité triangulaire, on a: $$d(a^{-1}bx,x)=d(bx,ax)< 2\varepsilon.$$
Ainsi $a^{-1}b\in G\cap V[x;2\varepsilon]\subseteq V$ et $b\in aV$. Ce qui implique $a=b$ en vertu de (\ref{eq:ab}).
Soit $D$ un nombre réel positif à la fois plus grand que le diamètre de $Ax$ et $4\varepsilon$. Notons $F$ le complémentaire de $(Ax)_{\varepsilon}$ dans $\mathbb{U}_{\mathfrak m}$. Considérons la fonction $1$-lipschitzienne $f$ à valeurs réelles définie sur $\mathbb{U}_{\mathfrak m}$ par:
\[f(y) = D - d(y,F).\]
Montrons que $f(y)>D-2\varepsilon$ pour tout $y\in \mathbb{U}_{\mathfrak m}$.\\
Si $y\in F$ alors $f(y)=D> D-2\varepsilon$.\\
Si $y\notin F$, alors $d(y,F)\leq d(y,ax)+d(ax,F)$ avec $a\in A$. Donc $d(y,F)\leq 2\varepsilon$. Ainsi, $f(y)\geq D-2\varepsilon$ pour tout $y\in \mathbb{U}_{\mathfrak m}$ et on a: $f(x)+f(y)\geq 2D-4\varepsilon > D$ car $D> 4\varepsilon$. Donc $f(x)+f(y)> D> d(x,y)$. On conclut que $f$ est de Kat\v etov. Par le Lemme \ref{extkat} et le fait que le $G$-orbite de $x$ a une cardinalité $< \mathfrak m$, il existe $z\in \mathbb{U}_{\mathfrak m}$ tel que \[\forall g\in G,~~f(gx)=d(gx,z).\]
Posons
\[W = V[z;\varepsilon]\cap G.\]
Soient $w\in W$ et $a\in A$. Alors $w\in V[z;\varepsilon]$ et $d(wz,z)< \varepsilon$. On a: $$f(wax)=d(wax,z)\leq d(wax,wz)+d(wz,z)\leq d(ax,z)+d(wz,z)=f(ax)+d(wz,z).$$ Donc $f(wax)-f(ax)< \varepsilon$. Ainsi $f(wax)\neq D$  car $f(ax)=D-\varepsilon$. Ce qui entraine
 $$wax\in \{s\in \mathbb{U}_{\mathfrak m}:\,\,f(s)< D\}.$$
Si $t\notin (Ax)_{\varepsilon}$, alors $t\in F$ et $d(t,F)=0$. Ainsi $f(t)=D$. Or pour tout $x\in \mathbb{U}_{\mathfrak m}$, on a: $f(x)\leq D$. Donc
$$\{x\in \mathbb{U}_{\mathfrak m}:\,\,f(x)< D\}= (Ax)_{\varepsilon}.$$
Comme $wax\in (Ax)_{\varepsilon}$, il existe $b\in A$ tel que $wax\in B_{\varepsilon}(bx)$ car $(Ax)_{\varepsilon}=\underset{b\in A}{\overset{}{\bigcup}}B_{\varepsilon}(bx)$.\\
Comme $d(wax,bx)<\varepsilon$, on a: $d(b^{-1}wax,x)< \varepsilon$, i.e $b^{-1}wa\in V[x,\varepsilon]$. Donc $$b^{-1}wa\in G\cap V[x,\varepsilon]\subseteq G\cap V[x,\varepsilon]\subseteq V.$$ Ainsi $wa\in AV$ comme souhaité.
\end{prof}\\
En utilisant le théorème \ref{protasin} et le théorème \ref{megrepest}, nous obtenons:
\begin{corollary}
Soit $\mathfrak m$ un cardinal infini non dénombrable vérifiant la condition (\ref{eq:mn}). Si $G$ est un sous-groupe de $Iso(\mathbb{U}_{\mathfrak m})$ de densitét $< \mathfrak{m}$, possédant la propriété (OB) et métrisable ou localement connexe, alors $G$ est SIN.
\end{corollary}
Voici à présent une liste de groupes topologiques qui ne se plonge pas dans le groupe $Iso(\mathbb{U}_{\mathfrak m})$ pour $\mathfrak{m}$ non-dénombrable vérifiant (\ref{eq:mn}). Ceci permettra de conclure que ce groupe n'est pas universel:
\begin{enumerate}
  \item Notons $\mathbb{U}_{1}$ l'espace métrique d'Urysohn de diamètre $1$. Le groupe $Iso(\mathbb{U}_{1})$ est à la fois métrisable de localement connexe ( est en effet homéomorphe à $\ell^{2} (\cite{mel2})$), possède la propriété (OB) (\cite{ros}, Théorème 1.5), mais n'est pas SIN(Ceci provient du fait que ce groupe est universel pour les groupes polonais (\cite{us5},\cite{us6})). Ainsi, le groupe $Iso(\mathbb{U}_{1})$ ne se plonge pas dans le groupe $Iso(\mathbb{U}_{\mathfrak m})$.
  \item Le groupe symétrique infini $S_{\infty}$ est polonais non-SIN et possède la propriété (OB).
  \item Le groupe $Homeo([0,1]^{\aleph_{0}})$ des homéomorphismes du cube de Hilbert sur lui-même muni de la topologie de la convergence compact est un groupe polonais non-SIN possédant la propriété (OB).
  \item Le groupe $U(\ell^{2})$ muni de la topologie forte posséde la propriété (OB)(\cite{atkin})
\end{enumerate}
En particulier, si $\mathfrak m$ est un cardinal infini non dénombrable vérifiant la condition (\ref{eq:mn}), alors le groupe $Iso(\mathbb{U}_{\mathfrak m})$ ne contient pas une copie du groupe topologique $Iso(\mathbb{U})$.

\begin{remark}
D'après ce qui précède, le groupe $Iso(\mathbb{U}_{\mathfrak m})$ n'est pas universel pour la classe des groupes topologiques de poids non dénombrable $\mathfrak m$. Néanmoins, on a les résultats suivant:
\end{remark}

\begin{theorem}\label{thmimb1}
Soit $\mathfrak m$ un cardinal infini non dénombrable vérifiant la condition (\ref{eq:mn}). Tout groupe métrisable SIN de poids $\leq \mathfrak m$ se plonge dans $Iso(\mathbb{U}_{\mathfrak m})$.
\end{theorem}
La preuve de ce théorème utilise la même technique utilisée par Uspenkij (\cite{us5},\cite{us6}) pour établir l'universalité du groupe $Iso(\mathbb{U})$ pour la classe des groupes topologiques vérifiant le deuxième axiome de dénombrabilité. Cette technique a été présentée en début de chapitre pour la démonstration du Théorème \ref{theousp}. La continuité de l'action initiale est assurée par l'observation suivante:

\begin{lemma}\label{lemimb1}
Soit $G$ un groupe opérant par isométries sur un espace métrique $X$ et soit $V$ un voisinage de l'identité dans $G$ et $\varepsilon> 0$. Supposons que la propriété suivante est satisfaite: $$(\star):\,\,\,\,\forall x\in X\,\text{et}\,\,\,\,v\in V,\,\text{on a}:\,d(x,vx)\leq\varepsilon.$$
Alors l'action canonique de $G$ sur $E(X)$ vérifie la propriété $(\star)$.
\end{lemma}
\begin{prof}
Soit $f$ une fonction $1$-lipschitzienne sur $X$.\\
 Pour tout $g\in V$ et $x\in X$, on a:

\begin{tabular}{ccc}
$d_{X}^{E}(^{g}f,f)$
& $=$ &
$\underset{x\in X}{\overset{}{\sup}}|f(g^{-1}x)-f(x)|$ \\
& $\leq$ & $ \underset{x\in X}{\overset{}{\sup}}\, d(g^{-1}x,x)$ \\
& $\leq $ & $ \varepsilon$
\end{tabular}
\end{prof}\\
\begin{prof}(\textbf{du théorème \ref{thmimb1}})
Soit $G$ un groupe SIN métrisable de poids $\leq \mathfrak m$. Fixons une métrique bi-invariante $d$ sur $G$. L'action par multiplication à gauche de $G$ sur lui-même vérifie la propriété $(\star)$. \\
En effet, soit $\varepsilon> 0$. Considérons $V=B_{\varepsilon}(e)$. On peut supposer sans nuire à la généralité comme $G$ est SIN que $V$ est invariant par conjuguaisons. Soit $g\in G$ et $v\in V$, on a: $d(g,vg)=d(e,g^{-1}vg)<\varepsilon$. Par le lemme \ref{lemimb1}, l'action canonique de $G$ sur $E(G)$ possède également la propriété $(\star)$. De proche en proche, l'action itérée de $G$ sur $\mathbb{U}_{\mathfrak m}=\underset{\mathfrak n < \mathfrak m }{\overset{}{\bigcup}}E^{(\mathfrak n)}_{\mathfrak m}(G)$ possède la propriété $(\star)$ i.e est continue. Il s'en suit que cette action détermine un plongement de $G$ dans $Iso(\mathbb{U}_{\mathfrak m})$.
\end{prof}

\begin{remark}
\begin{enumerate}
  \item Le groupe $U(\ell^{2})$ muni de la topologie d'opérateurs uniformes est métrisable, SIN et poids le continu $\mathfrak{c}$ (\cite{dierolf}). Ainsi, ce groupe se plonge dans le groupe $Iso(\mathbb{U}_{\mathfrak m})$ pour $\mathfrak{m}$ non-dénombrable par le théorème \ref{thmimb1}.
  \item Nous ne savons pas si tous les groupes SIN de poids $\leq \mathfrak m$ se plongent dans $Iso(\mathbb{U}_{\mathfrak m})$. Au même temps, tous les sous-groupes de $Iso(\mathbb{U}_{\mathfrak m})$ ne sont pas SIN comme le montre l'observation suivante:
\end{enumerate}

\end{remark}

\begin{definition}Soit $\mathfrak m$ un cardinal infini non dénombrable vérifiant la condition (\ref{eq:mn}). Un groupe topologique $G$ est un $P_{\mathfrak m}$-groupe si toute intersection d'une famille d'ouverts de $G$ de cardinalité $< \mathfrak m$ est un ouvert.
\end{definition}
Rappelons qu'un sous-ensemble $Q$ d'un ensemble ordonné $(P,\preceq)$ est dit cofinal si pour tout $x\in P$, il existe $y\in Q$ tel que $x\preceq y$ et le plus petit cardinal d'un sous-ensemble cofinal est appelé cofinalité de $(P,\preceq)$.\\
Soit $\mathfrak m$ un cardinal infini non dénombrable vérifiant (\ref{eq:mn}). D'après \cite{hafner}, il existe un filtre ordonné $\mathbb{K}$ de cofinalité $\mathfrak m$.\\
Le groupe linéaire $GL(\mathbb{K},n)$ muni de sa topologie naturelle induite de celle de $\mathbb{K}^{n^{2}}$ est un $P_{\mathfrak m}$-groupe (\cite{vp5}). Ce groupe est par ailleurs non SIN(\cite{vp5}).\\
Notons $(\mathcal{C})$ la classe des $P_{\mathfrak m}$-groupes, on a le résultat:
\begin{theorem}
Soit $\mathfrak m$ un cardinal infini non dénombrable vérifiant (\ref{eq:mn}). Si $G\in(\mathcal{C})$ est de poids $\mathfrak m$, alors $G$ est isomorphe à un sous-groupe de $Iso(\mathbb{U}_{\mathfrak m})$.
\end{theorem}
La démonstration de cet important théorème utilise la même technique utilisée par Uspenkij dans \cite{us5} pour établir l'universalité  du groupe $Iso(\mathbb{U})$. La continuité de l'action initiale est garantie par le lemme suivant:
\begin{lemma}
Soit $G\in (\mathcal{C})$ un groupe topologique de poids $\mathfrak m$. Si $G$ opère continûment par isométries sur un espace métrique $X$, alors l'extension canonique de cette action sur $E_{\mathfrak m}(X)$(la representation régulière à gauche définie par  $^gf(x)=f(g^{-1}x)$) est continue.
\end{lemma}
\begin{prof}
Soit $f$ une fonction $1$-lipschitzienne sur $X$ controllée par un sous-ensemble $A$  de cardinalité $<{\mathfrak m}$. Pour $\varepsilon>0$, posons
\[V=\bigcap_{a\in A}V[a;\varepsilon].\]
$V$ est un voisinage de l'identité dans $G$ puisque $G\in \mathcal{C}$. Si $g\in V$, alors pour tout $a\in A$ on a: $|f(ga)-f(a)|<\varepsilon$. Cette condition implique que pour tout $x\in X$ et $g\in V$ on a $|f(gx)-f(x)|<\varepsilon$, et ainsi, l'action de $G$ sur $E_{\mathfrak m}(X)$ est continue.
\end{prof}
\begin{remark}
Le groupe $Iso(\mathbb{U}_{\mathfrak m})$ est extrêmement moyennable. Ceci découle des résultats de \cite{vp4} (Théorèmes $6.5$ et $6.6$)
\end{remark}
\section{Conclusion et perspectives}
Au terme de cette étude du groupe des isométries $Iso(\mathbb{U}_{\mathfrak m})$, des questions démeurent néanmoins à l'esprit.
\begin{enumerate}
  \item Le groupe $Iso(\mathbb{U}_{\mathfrak m})$ est-il connexe? connexe par arcs? comme le groupe $Iso(\mathbb{U})$ (\cite{mel2})ou alors localement connexe?
  \item Un espace métrique ultrahomogène $X$ (< $\omega$-homogène) est dit stable par oscillations (voir chapitre $8$ dans \cite{vp1}) si toute partition $\gamma$ de $X$ contient un élément $A\in \gamma$ tel que pour tout $\varepsilon> 0$, son voisinage d'ordre $\varepsilon,\,A_{\varepsilon}$ contient une copie isométrique de $X$. Un théorème célèbre de Odell et Schlumprecht (\cite{odell}) affirme que la sphère unité d'un espace de Hilbert n'est pas stable par oscillations (a des distortions). De même un résultat récent de Nguyen Van Thé et Sauer (\cite{nguyen}) (basé sur un travail préalable de Lopez-Abad et Nguyen Van Thé (\cite{nguyen2})) établit la stabilité par oscillations de la sphère unité de l'espace universel polonais d'Urysohn $\mathbb{U}$. La sphère unité de l'espace $\mathbb{U}_{\mathfrak m}$ est-elle stable par oscillations?
\end{enumerate}

\appendix

\addtocontents{toc}{
\protect\renewcommand*\protect\cftchappresnum{\appendixname~}}

\chapter{ \quad Notions sur les ordinaux et cardinaux}
\begin{definition}
Un ensemble d'ensembles $A$ est dit transitif si tout élément de $A$ est un élément de $A$. De façon précise, $A$ est transitif si on a:
$$x\in a\in A\Longrightarrow x\in A,$$
c'est-à-dire si $a\in A$ entraîne $a\subseteq A$.
\end{definition}

La notion d'ensembles transitif n'est pas familière, et il semble clair que la plupart des ensembles ne sont pas transitifs.

\begin{definition}
On dit qu'un ensemble $\alpha$ est un ordinal si $\alpha$ est un ensemble transitif et que la restriction de $\in$ à $\alpha$ est un bon ordre strict.
\end{definition}
\begin{remark}

\begin{enumerate}
  \item $\emptyset$ est un ordinal, et si $\alpha$ est un ordinal, il en est de même de $S(\alpha)=\alpha\cup \alpha$.
  \item  Soient $\alpha$ et $\beta$ deux ordinaux. On déclare $\alpha< \beta$ vrai si $\alpha$ est un élément de $\beta$. Alors $<$ est un bon ordre sur les ordinaux.
  \item $S(\alpha)$ est successeur immédiat de $\alpha$. Tout ensemble non vide d'ordinaux $A$ a un plus petit élément, à savoir $\bigcap_{\alpha\in A}\alpha$ et une borne supérieure, à savoir $\bigcup_{\alpha\in A}\alpha$
  \item Aucun ensemble ne contient tous les ordinaux
  \item Les ordinaux non nuls se partagent en ordinaux successeurs, du type $S(\alpha)$, et ordinaux limites, qui vérifient $\lambda=\bigcup_{\alpha\in \lambda}\alpha$. Le plus petit ordinal limite est $\omega$ borne supérieure des ordinaux finis.
\end{enumerate}
\end{remark}

\begin{definition}
Un cardinal est un ordinal qui n'est en bijection avec aucun de ses prédécésseurs.
\end{definition}Ainsi,
\begin{enumerate}
  \item $\aleph_{0}:=\omega$
  \item $\aleph_{1}:=$ le premier ordinal non-dénombrable.
  \item $\aleph_{2}:=$ le premier ordinal qui n'est pas en bijection avec $\aleph_{1}$.
\end{enumerate}
En général, pour tout ordinal $\alpha$, on définit:
\begin{enumerate}
  \item $\aleph_{\alpha+1}:=$ le premier ordinal qui n'est pas en bijection avec $\aleph_{\alpha}$.
  \item $\aleph_{\alpha}=\bigcup_{\delta< \alpha}\aleph_{\delta}$ si $\alpha$ est un ordinal limite.
\end{enumerate}

\begin{definition}
Soient $\alpha$ et $\beta$ deux ordinaux. Rappelons les notions suivantes:
\begin{enumerate}
  \item  Une application $f:\alpha\longrightarrow \beta$ est dite cofinale de $\alpha$ dans $\beta$ si $Im(f)$ est non bornée dans $\beta$. (i.e $\forall\,\gamma \in \beta,\,\exists\,x\in \alpha$ tel que $f(x)\geq \gamma$.)
  \item La cofinalité de $\beta$, notée $cof(\beta)$, est le plus petit ordinal $\alpha$ tel qu'il existe une application cofinale  $f:\alpha\longrightarrow \beta$.
  \item $\beta$ est dit limite si $\beta=\bigcup_{\tau<\beta}\tau$ ou de manière équivalente $\beta$ est un point d'accumulation lorsque l'ensemble des ordinaux est muni de la topologie de l'ordre.
   \item $\beta$ est dit régulier si $\beta$ est un ordinal limite et $cof(\beta)=\beta$.
   \item Un cardinal $\kappa$ est dit régulier s'il n'est pas la limite d'une suite avec pour ensemble d'indice un ordinal plus petit que $\kappa$.\\
Tous les cardinaux successeurs sont par exemple réguliers. Par contre, $\aleph_{\omega}=\bigcup_{i<\omega}\aleph_{i}$ n'est pas régulier (est singulier) puisqu'il est la limite de la suite dénombrable des $\aleph_{n}$.
 ($\aleph_{\omega}=\lim_{i\in \omega}\aleph_{i}$).
\item Un cardinal $\alpha$ est dit fortement inaccessible s'il est régulier et si pour tout cardinal $\beta < \alpha$, on a $2^{\beta}< \alpha$.
\end{enumerate}
\end{definition}
Soit $\mathfrak m$ un cardinal infini non dénombrable vérifiant la condition (\ref{eq:mn}).\\
Le cardinal $\aleph_{0}$ par exemple vérifie la condition (\ref{eq:mn}).\\
Montrons que tout cardinal fortement inaccessible vérifie la condition (\ref{eq:mn}).\\
 Soit $\mathfrak m$ un cardinal inaccessible. Posons $\eta=\sup\left\{{\mathfrak m}^{\mathfrak n}\colon {\mathfrak n}<{\mathfrak m}\right\}$.\\
Pour tout $\mathfrak n< \mathfrak m$, on a: $\mathfrak m^{\mathfrak n}\geq \mathfrak m$. Ainsi $\eta\geq \mathfrak m$. \\
Il nous suffit pour conclure de montrer que $\eta\leq \mathfrak m$. Pour cela, montrons que $\mathfrak m^{\mathfrak n}\leq\mathfrak m$ pour tout $\mathfrak n< \mathfrak m$.\\
    Fixons $\mathfrak n< \mathfrak m$. Soient $M$ un ensemble bien ordonné de cardinalité $\mathfrak m$ et $N$ un ensemble de cardinalité $\mathfrak n$. Notons $\mathfrak{F}$ l'ensemble des fonctions $f:N\longrightarrow M$. Rappelons que $\mathfrak m^{\mathfrak n}$ désigne le cardinal de $\mathfrak{F}$. Nous devons montrer que le cardinal de $\mathfrak{F}$ est $<\mathfrak m$.\\
    Puisque $\mathfrak m$ est un cardinal régulier et $\mathfrak n< \mathfrak m$, toute fonction $f$ de $N$ dans $M$ est telle que $Im(f)$ est bornée par un élément $b$ dans $M$.\\
     Pour tout $b$ particulier, le nombre de fonctions $f$ donc l'image est borné par $b$ est $\kappa^{\mathfrak n}$, où $\kappa$ est le nombre de prédécesseurs de $b$.\\
    Notons que pour tout $b\in M$, le $\kappa$ correspondant est $< \mathfrak m$. Puisque $\mathfrak m$ est fortement inaccessible, $\kappa^{\mathfrak n}< \mathfrak m$. Ainsi l'ensemble $S_{b}$ des fonctions $f:N\longrightarrow M$ donc l'image est bornée par $b$ a un cardinal $< \mathfrak m$. \\
    Notons enfin que $\mathfrak{F}=\bigcup_{b\in M}S_{b}$. Donc $\mathfrak{F}$ est la réunion de $\mathfrak m$ ensembles de cardinal au plus $\mathfrak m$ (en réalité chacun des ensembles a un cardinal strictement plus petit que $\mathfrak m$). Ainsi la réunion a un cardinal au plus $\mathfrak m$, comme souhaité.
\begin{remark}
\begin{enumerate}
  \item L'existence des cardinaux fortement inaccessible est indépendante des axiomes ZFC. Par contre, nous avons eu besoin des axiomes ZFC pour montrer que tout cardinal fortement inaccessible vérifie la condition (\ref{eq:mn}).
  \item L'hypothèse du continu généralisée (HCG) affirme que $2^{\kappa}=\kappa^{+}$ pour tout cardinal infini $\kappa$. Moyennant
HCG, tout cardinal infini successeur $k^{+}$ vérifie la condition (\ref{eq:mn}).\\
En effet, pour $n < k^{+}$, on a: $(k^+)^n = (2^k)^n = 2^{kn} = 2^k = k^{+}.$
\item L'existence des cardinaux vérifiant la condition (\ref{eq:mn}) est indépendante des axiomes ZFC.
\end{enumerate}
\end{remark}
\chapter{ \quad Quelques notions de topologie}
\section{Poids et densité d'un espace topologique}
\begin{definition}
Soit $X$ un espace topologique. Rappelons que le poids\index{espace topologique!poids d'un} de $X$ est le nombre cardinal $\omega(X)=\min\{|\mathcal{B}|:\,\,\mathcal{B}\, \text{est une base de}\, $X$\}$ et la densité\index{espace topologique!densité d'un} de $X$ est le nombre cardinal $d(X)=\min\{|A|:\,\,A\, \text{est une partie dense de}\, $X$\}$.
\end{definition}
\begin{remark}
\begin{enumerate}
  \item
  Pour tout espace topologique, nous avons $d(X)\leq \omega(X)$.\\ En effet, si $\mathcal{B}=\{U_{s}\}_{s\in S}$ est une base pour $X$ constituée des ensembles non-vides et tels que $|S|=m=\omega(X)$. Choisissons pour tout $s\in S$ un point $a_{s}\in U_{s}$ et montrons que l'ensemble $A=\{a_{s}:\,\,s\in S\}$ est dense dans $X$. En effet, tout ouvert non vide de $X$ contient un élément $U_{s}\in \mathcal{B}$, ainsi il contient le point $a_{s}\in A$. Puisque $|A|\leq |S|=m$, on a $d(X)\leq w(X)$. En particulier, tout espace vérifiant le deuxième axiome de dénombrabilité est séparable
  \item Si $X$ est un espace métrisable, alors $d(X)=\omega(X)$.\\
   En effet, si $A$ est un sous-espace dense de $X$ de cardinalité $|A|=d(X)=k$ avec $k$ infini, considérons $\mathcal{B}=\{B(x,r):\,\,x\in A;\,\,r\in \mathbb{Q}\}$. Il est clair que $|\mathcal{B}|\leq \aleph_{0}\times k=k$. Montrons pour terminer que $\mathcal{B}$ est une base de $X$.\\
Soit $U$ un ouvert non vide de $X$ et soit $x\in U$ alors il existe $\varepsilon> 0$ tel que $B(x,\varepsilon)\subset U$. Puisque $A$ est dense dans $X$, il existe $a\in A$ et $a\in B(x,\frac{\varepsilon}{3})$. Soit $r\in \mathbb{Q}$ tel que $\frac{\varepsilon}{3}< r <\frac{\varepsilon}{2}$, alors $x\in B(a,r)$ et $B(a,r)\in \mathcal{B}$. De plus, $B(a,r)\subset B(x,\varepsilon)$. En effet, soit $y\in B(a,r)$, on a $d(x,y)\leq d(x,a)+d(a,y)<r+r<\varepsilon$. Posons $B_{x}=B(a,r)$. On a donc montrer que:
 pour tout $x\in U$, il existe $B_{x}\in \mathcal{B}$ tel que $x\in B_{x}\subset U$.
\end{enumerate}
\end{remark}

\section{Topologie compact-ouvert}\index{topologie!compact-ouvert}
Soient $X$ et $Y$ deux ensembles arbitraires et soient $A\subset X$ et $B\subset Y$. Nous écrirons $Y^{X}$ pour désigner l'ensemble des applications de $X$ dans $Y$ et $F(A,B)$ le sous-ensemble de $Y^{X}$ formé des applications appliquant $A$ dans $B$:
$$F(A,B)=\{f\in Y^{X}:\,\,f(A)\subset B\}$$
 A présent, soient $X$ et $Y$ deux espaces topologiques et soit $\mathcal{A}$ l'ensemble des compacts de $X$ et $\mathcal{G}$ l'ensemble des ouverts de $Y$. La topologie sur $Y^{X}$ engendrée par $$\mathcal{S}=\{F(A,B):\,\,A\in \mathcal{A},\,\,\,B\in \mathcal{G}\}$$
 est appelée la topologie définie par un compact et un ouvert ou tout simplement topologie compacte-ouverte et $\mathcal{S}$ est la sous-base de définition de cette topologie.

 \section{Groupes polonais}\index{groupe!polonais}
 Un espace métrisable à base dénombrable est un espace polonais\index{espace topologique!polonais} si sa topologie peut être définie par une distance qui en fait un espace complet. Cette terminologie a été introduite par le groupe Nicolas Bourbaki, dans son volume sur la topologie générale.\\
 Rappelons qu'un sous-espace $A$ d'un espace topologique $X$ est dit $G_{\delta}$ s'il existe une famille dénombrable $(O_{i})_{i\in I}$ d'ouverts de $X$ tels que $A=\underset{i=1}{\overset{\infty}{\bigcap}}O_{i}$. Tout sous-espace $G_{\delta}$ d'un espace polonais est également polonais (\cite{inf}).\\
 $\mathbb{N}$ étant muni de la topologie discrète, l'espace de Baire $\mathbb{N}^{\mathbb{N}}$ est un espace polonais.
 On dit qu'un groupe topologique $G$ est un groupe polonais\index{groupe!polonais} si la topologie
de $G$ est polonaise. Le lecteur pourra se reporter à \cite{becker} ou \cite{kechris} pour une introduction complète à la
théorie des groupes polonais.\\
Nous allons à présent introduire deux exemples classiques de groupes polonais qui seront très utile dans la suite de la thèse.
 \subsection{Le groupe symétrique infini $S_{\infty}$}
 $S_{\infty}$ désigne le groupe de toutes les bijections de $\mathbb{N}$ dans $\mathbb{N}$. $S_{\infty}$ muni de la topologie induite de celle de l'espace de Baire $\mathbb{N}^{\mathbb{N}}$, est un $G_{\delta}$-sous-espace.\\ En effet, en posant $$\mathcal{A}=\{\sigma \in \mathbb{N}^{\mathbb{N}}:\,\,\sigma\,\, \text{est injective}\}$$ et $$\mathcal{B}=\{\sigma \in \mathbb{N}^{\mathbb{N}}:\,\,\sigma\,\, \text{est surjective}\}.$$ Nous avons donc $S_{\infty}=\mathcal{A}\cap \mathcal{B}$. En remarquant que
  $$\mathcal{A}=\underset{m\neq m}{\overset{}{\bigcap}}\{\sigma \in \mathbb{N}^{\mathbb{N}}:\,\,\sigma(n)\neq \sigma(m)\}$$  et $$\mathcal{B}=\underset{m\in \mathbb{N}}{\overset{}{\bigcap}}\underset{n\in \mathbb{N}}{\overset{}{\bigcup}}\{\sigma \in \mathbb{N}^{\mathbb{N}}:\,\,\sigma(n)=m\},$$ on conclut que $\mathcal{A}$ et $\mathcal{B}$ sont des $G_{\delta}$ sous-espace de $\mathbb{N}^{\mathbb{N}}$. Ainsi, $S_{\infty}$ est un $G_{\delta}$ sous-espace de $\mathbb{N}^{\mathbb{N}}$. C'est donc un groupe polonais.\\
  Soit $A\subseteq \mathbb{N}$ tel que $|A|< \infty$. Posons $St_{A}=\{g\in S_{\infty},\,\,\,\forall a\in A\,\,\,\,\,g(a)=a\}$.\\
  La famille $\Sigma=\{St_{A}:\,\, |A|< \infty\}$ est une base de voisinages de $e$ pour la topologie polonaise de $S_{\infty}$ formé des sous-groupes ouverts. On dit encore que $S_{\infty}$ est non archimédien.
\subsection{Le groupe $Homeo(D^{\aleph_{0}})$}
Commençons par des rappels sur l'ensemble de Cantor.
Notons $D=\{0,1\}$. L'espace produit $D^{\aleph_{0}}$ est l'ensemble de Cantor.
Posons $I=[0,1]$. Considérons l'opération $T$: à enlever le "tiers central" sur $I$.\\On pose $I_{1}=T([0,1])$; $I_{n+1}=T(I_{n})$; $C=\underset{n=1}{\overset{\infty}{\bigcap}}I_{n}$. $C$ est appelé ensemble triadique de Cantor.
$I_{n}$ est la réunion des $2^{n}$ segments disjoints de longueur $3^{-n}$, d'origine $\underset{j=1}{\overset{n}{\sum}}\alpha_{j}3^{-j}$ où $\alpha_{j}=0,2$.\\En effet, la propriété a lieu à l'étape 1: $C_{1}=[0,\frac{1}{3}]\cup[\frac{2}{3},1]$.\\
 Si elle a lieu aux étapes $1,...,n$, on a: $$C_{n}=\underset{\alpha}{\overset{}{\bigcup}}I_{n}(\alpha)$$ où $\alpha=(\alpha_{1},...\alpha_{n})$ parcourt $\{0,2\}^{n}$ et où $$I_{n}(\alpha)=[\underset{j=1}{\overset{n}{\sum}}\alpha_{j}3^{-j}, \underset{j=1}{\overset{n}{\sum}}\alpha_{j}3^{-j}+3^{-n}].$$ on voit que
 $$T[I_{n}(\alpha)]=I_{n+1}(\alpha^{0})\cup I_{n+1}(\alpha^{1})$$ où $\alpha^{0}=(\alpha_{1},...\alpha_{n},0)$ et $\alpha^{0}=(\alpha_{1},...\alpha_{n},2)$, ce qui donne
 $$C_{n+1}=\underset{\alpha}{\overset{}{\bigcup}}(I_{n+1}(\alpha^{0})\cup I_{n+1}(\alpha^{1}))=\underset{\beta}{\overset{}{\bigcup}}I_{n+1}(\beta)$$ où $\beta$ parcourt $\{0,2\}^{n+1}$. On a donc: $I_{n}=\underset{j=1}{\overset{2^{n}}{\bigcup}}I_{n}^{j}$. De plus, l'application

 \[\begin{array}{lll}
\varphi:& D^{\aleph_{0}} &\longrightarrow C\\
&x=(x_{j}) &\longmapsto \underset{j=1}{\overset{\infty}{\sum}}x_{j}3^{-j}\\
\end{array}\]
 est un homéomorphisme. L'ensemble de Cantor désigne selon le cas l'ensemble triadique $C$ ou $D^{\aleph_{0}}$.\\
 $Homeo(D^{\aleph_{0}})$ désigne le groupe de tous les homéomorphismes de $D^{\aleph_{0}}$ sur lui même.
Comme $D^{\aleph_{0}}$ est compact, la topologie de la convergence uniforme sur $Homeo(D^{\aleph_{0}})$ coïncide avec la topologie compact-ouvert.\\
On montre comme dans le cas de $S_{\infty}$ que $Homeo(D^{\aleph_{0}})$ est un $G_{\delta}$ sous-espace de $C(D^{\aleph_{0}},D^{\aleph_{0}})$. Ce qui permet de conclure que $Homeo(D^{\aleph_{0}})$ est polonais\\
 Pour tout $n\in \mathbb{N}$, posons $H_{n}=\{g\in G:\,\,g(I_{n}^{j})=I_{n}^{j}\,\,\,\forall\,\,1\leq j\leq 2^{n}\}$.\\
 La famille $\{H_{n}:\,n\in \mathbb{N}\}$ constitue une base de voisinages de $e$ dans $Homeo(D^{\aleph_{0}})$ consititué de sous-groupes ouverts.

\section{Compactifié de Stone-\v Cech}\index{compactifié!de Stone-\v Cech}
Soit $X$ un espace séparé complètement régulier. Notons $\mathbb{I}=[0,1]$ et $C(X)$ l'ensemble des applications continues de $X$ dans $\mathbb{I}$. Posons $\widetilde{X}=\mathbb{I}^{C(X)}$ et considérons l'application

\[ \xymatrix@C=2cm@R=0.3em{
 \varepsilon:  X  \ar[r] &  \mathbb{I}^{C(X)} &  \\
 x \ar@{|->}[r] & \,\,\,\,\varepsilon(x) : C(X) \ar[r] &  \mathbb{I} \\
& \,\,\,\,\,\,\,\,f \ar@{|->}[r] &  \varepsilon(x)(f)=f(x) }
\]

Notons $\beta X$ l'adhérence de $\varepsilon(X)$ dans $\widetilde{X}$. $\beta X$ est un compactifié de $X$ appelé Compactifié de Stone-\v Cech de $X$. De plus, $\beta X$ vérifie les propriétés suivantes:
\begin{theorem}
\begin{enumerate}
  \item Soit $X$ un espace complètement régulier séparé. Si $Y$ est un espace compact et $g:X\longrightarrow Y$ est une application continue, alors $g$ se prolonge uniquement en une application continue $\widetilde{g}:\beta X\longrightarrow Y$
  \item La topologie induite sur $X$ par celle de $\beta X$ coïncide avec la topologie initiale de $X$
  \item $X$ est dense dans $\beta X$.
\end{enumerate}
\end{theorem}

\section{Dimension de Lebesgue}
\begin{definition}
Soient $\mathcal{A}$ et $\mathcal{B}$ deux recouvrements d'un espace topologique $X$.
\begin{enumerate}
  \item $\mathcal{B}$ est un \index{raffinement} raffinement de $\mathcal{A}$ si pour tout $B\in\mathcal{B}$, il existe $A\in \mathcal{A}$ tel que $B\subset A$. Ce sera un raffinement ouvert si $\mathcal{B}$ est également constitué d'ensembles ouverts.
  \item L'ordre de $\mathcal{A}$ est défini comme étant le plus grand nombre entier $n$ tel qu'il existe $n+1$ éléments de $\mathcal{A}$ ayant une intersection non vide.
\end{enumerate}
\end{definition}
\begin{definition}
Soit $X$ un espace topologique, la \index{dimension!de Lebesgue}dimension de Lebesgue $dim(X)$ de $X$ est le plus petit entier $n$ tel que tout recouvrement ouvert de $X$ possède un raffinement d'ordre $\leq n$. Si un pareil entier $n$ n'existe pas, on dit que $dim(X)=\infty$.
\end{definition}
\begin{example}(\cite{munkres}, Exemples $1$ et $2$, page $305$)
\begin{enumerate}
  \item Tout sous-espace compact $X$ de $\mathbb{R}$ est de dimension de Lebesgue au plus $1$.
  \item L'espace compact $X=[0,1]$ est de dimension de Lebesgue $1$.
\end{enumerate}
\end{example}
\begin{theorem}(\cite{ark})
Un espace compact $X$ est totalement discontinu si et seulement si $dim(X)=0$.
\end{theorem}
\section{Produit diagonal}\index{produit diagonal}
\begin{definition}
Soient $X$ un espace topologique, $(X_i)_{i\in I}$ une famille d'espaces topologiques et $(f_i)_{i\in I}$ une famille de fonctions de $X$ dans $X_{i}$. Le produit diagonal des fonctions $f_i$ est la fonction $\Delta_{i}f_{i}:X\longrightarrow \Pi_{i}X_{i}$ définie par $\Delta_{i}f_{i}(x)=(f_{i}(x))$.
\end{definition}
\begin{lemma}
Si toutes les fonctions $f_{i}$ sont continues, alors le produit diagonal des fonctions $f_i$ est également continue.
\end{lemma}

\section{Systèmes projectifs}
Rappelons la notion de \index{système!projectif}système projectif de $G$-espaces en nous appuyant sur \cite{Bour}.\\
$G$ étant un groupe topologique, un système projectif de $G$-espaces est la donnée d'un triplet $(X_{i},\pi_{ij}, I)$ où $(I,\preceq )$ est un ensemble ordonné, $\{X_{i}\}_{i\in I}$ une famille  de $G$-espaces compacts et une famille d'applications équivariantes $\pi_{ij}:X_{j}\longrightarrow X_{i}$, pour $i\preceq j$ tel que $\pi_{ii}$ est l'identité de $X_{i}$ pour tout $i\in I$ et $\pi_{i k}=\pi_{i j}\circ \pi_{j k}$ pour tout $i \preceq j \preceq k$. La limite projective $X=\underset{\longleftarrow}{\overset{}{\lim}}\,(X_{i},\pi_{ij})$ du système projectif de $G$-espaces $(X_{i},\pi_{ij}, I)$ est le $G$-espace définie comme suit:
$\underset{i\in I}{\overset{}{\prod}}X_{i}$ est l'ensemble produit des ensembles $(X_{i})_{i\in I}$ et $pr_{i}:\underset{i\in I}{\overset{}{\prod}}X_{i}\longrightarrow X_{i}$ la projection sur le facteur $X_{i}$.\\
$$\underset{\longleftarrow}{\overset{}{\lim}}\,(X_{i},\pi_{ij})=\{x=(x_{i})\in \underset{i\in I}{\overset{}{\prod}}X_{i}/\,\,\, x_{i}=\pi_{ij}(x_{j})\,\,\,\,i\preceq j\}.$$
Ou encore  $$\underset{\longleftarrow}{\overset{}{\lim}}\,(X_{i},\pi_{ij})=\{x=(x_{i})\in \underset{i\in I}{\overset{}{\prod}}X_{i}/\,\,\, pr_{i}(x)=\pi_{ij}\circ pr_{j}(x)\,\,\,\,i\preceq j\}.$$ On note tout simplement $X=\underset{\longleftarrow}{\overset{}{\lim}}\,X_{i}$
 si aucune confusion n'est possible.
$\underset{\longleftarrow}{\overset{}{\lim}}\,X_{i}$ est un sous-espace fermé, donc compact de $\underset{i\in I}{\overset{}{\prod}}X_{i}$. $G$ opère continûment sur $\underset{i\in I}{\overset{}{\prod}}X_{i}$ par $g.x=(g.x_{i})$ où $x=(x_{i})$.
La restriction $\pi_{i}$ de la projection $pr_{i}$ à X est l'application canonique de $X$ dans $X_{i}$. On a la relation:
 $\pi_{i}=\pi_{ij}\circ \pi_{j}$ pour tout $i\preceq j$.

\section{Structure uniforme}
Commençons par introduire les notations suivantes:\\
Soit $X$ un ensemble non vide et soient $V$ et $W$ deux parties de $X\times X$.
\begin{enumerate}
  \item $V\circ W$ désigne l'ensemble des couples $(x,y)\in X\times X$ tel qu'il existe $z\in X$ tel que $(x,z)\in W$ et $(z,y)\in V$.
  \item $V^{-1}$ est l'ensemble des couples $(x,y)\in X\times X$ tel que $(y,x)\in V$.
\end{enumerate}
\begin{definition}
On appelle structure uniforme sur un ensemble $X$, une famille non vide $\mathcal{U}_{X}$ de parties de $X\times X$ vérifiant les conditions suivantes:
\begin{enumerate}
\item Tout élément $U\in \mathcal{U}_{X}$ contient la diagonale $\Delta_{X}$ de $X\times X$.
\item Si $U\in \mathcal{U}_{X}$, alors $U^{-1}\in \mathcal{U}_{X}$
\item Si $U\in \mathcal{U}_{X}$, alors il existe $V\in \mathcal{U}_{X}$ tel que $V\circ V\subseteq U$.
\item Si $U,V\in \mathcal{U}_{X}$, alors $U\cap V\in \mathcal{U}_{X}$.
\item Si $U\in \mathcal{U}_{X}$ et $U\subseteq V$ alors $V\in \mathcal{U}_{X}$
\end{enumerate}
\end{definition}
La paire $(X,\mathcal{U}_{X})$ constituée d'un ensemble et d'une structure uniforme est appelé un espace uniforme et les éléments de $\mathcal{U}_{X}$ sont appelés les entourages\index{entourage} de la diagonale.
\begin{definition}
Une sous-famille $\mathcal{B}\subseteq \mathcal{U}_{X}$ est une base de la structure uniforme $\mathcal{U}_{X}$ si pour tous $U,V\in \mathcal{B}$, il existe $C\in \mathcal{B}$ tel que $C\subseteq U\cap V$ et tout entourage $V\in \mathcal{U}_{X}$ contient un élément $U\in \mathcal{B}$.\\
Une famille $\mathcal{B}$ de parties de $X\times X$ est une base pour une structure uniforme si et seulement si elle vérifie les conditions suivantes:
\begin{enumerate}
  \item Pour tous $U,V\in \mathcal{B}$, il existe $C\in \mathcal{B}$ tel que $C\subseteq U\cap V$.
  \item Tout élément $U\in \mathcal{B}$ contient la diagonale $\Delta_{X}$ de $X\times X$
  \item Si $V\in \mathcal{B}$, alors il existe $U\in \mathcal{B}$ tel que $U\subseteq V^{-1}$
  \item Si $V\in \mathcal{B}$, alors il existe $U\in \mathcal{B}$ tel que  $U\circ U\subseteq V^{-1}$
\end{enumerate}
\begin{example}
\begin{enumerate}
  \item Toute métrique $d$ sur un ensemble $M$ engendre une structure uniforme $\mathcal{U}_{d}$ sur $M$ ayant pour base les ensmebles $$U_{d}^{\varepsilon}=\{(x,y)\in M\times M:\,\,d(x,y)< \varepsilon\}$$ où $\varepsilon> 0$.
  \item L'ensemble de tous les sous-ensembles de $X\times X$ contenant la diagonale forme une structure uniforme sur $X$ dite structure uniforme discrète.
  \item L'ensemble constitué uniquement de $X\times X$ est une structure uniforme sur $X$ dite structure uniforme triviale sur $X$.
\end{enumerate}
\end{example}
\end{definition}
\begin{definition}
Soit $G$ un groupe topologique. Il existe sur $G$ deux structures uniformes bien connues: La structure uniforme droite que nous noterons $\mathcal{U}_{R}(G)$ et la structure uniforme gauche que nous noterons $\mathcal{U}_{L}(G)$.
\begin{enumerate}
  \item Soit $V$ un voisinage de $e$ dans $G$. Posons $V_{d}=\{(x,y):\,\,\,xy^{-1}\in V\}$.\\
  La famille des ensembles $V_{d}$ lorsque $V$ parcourt une base de voisinages de $e$ forme une base d'entourages pour la structure uniforme droite $\mathcal{U}_{R}(G)$ sur $G$.
  \item Soit $V$ un voisinage de $e$ dans $G$. Posons $V_{g}=\{(x,y):\,\,\,x^{-1}y\in V\}$.\\
  La famille des ensembles $V_{g}$ lorsque $V$ parcourt une base de voisinages de $e$ forme une base d'entourages pour la structure uniforme gauche $\mathcal{U}_{L}(G)$ sur $G$.
\end{enumerate}
\end{definition}
\begin{definition}(Topologie d'un espace uniforme)\index{topologie!d'un espace uniforme}
Soit $(X,\mathcal{U}_{X})$ un espace uniforme. La structure uniforme $\mathcal{U}_{X}$ induie sur $X$ une topologie appelée toplogie uniforme.\\
Cette topologie est constituée des parties $\mathcal{O}\subseteq X$ telles que pour tout $x\in \mathcal{O}$, il existe $V\subseteq \mathcal{U}_{X}$
tel que $V[x]\subseteq \mathcal{O}$. Avec  $V[x]=\{y\in X:\,\,\,(x,y)\in V\}$.
\end{definition}
\begin{remark}
D'après ce qui précède, tout espace uniforme induit un espace topologique. La réciproque est donnée par le résultat suivant:
\end{remark}
\begin{theorem}\cite{will}
Un espace topologique $X$ est uniformisable si et seulement s'il est complètement régulier.
\end{theorem}
\begin{proposition}
Soit $(G,\tau)$ un groupe topologique. Les deux topologies uniformes induites respectivement par la structure uniforme droite et gauche sur $G$ coïncide avec la topologie initiale $\tau$
\end{proposition}
\begin{definition}
Une application $f$ d'un espace uniforme $(X,\mathcal{U}_{X})$ dans un espace uniforme $(Y,\mathcal{U}_{Y})$ est dite uniformément continue si,
 pour tout entourage $V$ de $\mathcal{U}_{Y}$, il existe un entourage $U$ de $\mathcal{U}_{X}$ tel que la relation $(x,y)\in U$ entraîne $(f(x),f(y))\in V$.
\end{definition}
\begin{remark} De manière équivalente, Une application $f$ d'un espace uniforme $(X,\mathcal{U}_{X})$ dans un espace uniforme $(Y,\mathcal{U}_{Y})$ est dite uniformément continue si,
 pour tout entourage $V$ de $\mathcal{U}_{Y}$, on a: $\{(x,y)\in X\times X:\,\,\,(f(x),f(y))\in V\}\subseteq \mathcal{U}_{X}$.
\end{remark}
\begin{remark}
\begin{enumerate}
  \item Si dans la définition précédente, $X$ est un groupe topologique, alors l'application $f$ sera dite uniformément continue à gauche\index{fonction!uniformément continue à droite} si elle est uniformément continue par rapport à la structure uniforme gauche. De manière similaire, l'application $f$ sera dite uniformément continue à droite si elle est uniformément
      \index{fonction!uniformément continue à gauche} continue par rapport à la structure uniforme droite.
  \item En particulier, une fonction sur un groupe topologique $G$ à valeures réelles sera dite uniformément continue à gauche si:
  pour tout $\varepsilon >0$, il existe un voisinage $V$ de $e$ tel que $x^{-1}y\in V\Longrightarrow |f(x)-f(y)|< \varepsilon$.
  \item De même, une fonction définie sur un groupe topologique $G$ à valeures réelles sera dite uniformément continue à droite si:
  pour tout $\varepsilon >0$, il existe un voisinage $V$ de $e$ tel que $xy^{-1}\in V\Longrightarrow |f(x)-f(y)|< \varepsilon$.
\end{enumerate}
\end{remark}
\begin{example}
\begin{enumerate}
  \item L'inversion sur $G$ est uniformément continue de la structure uniforme droite $\mathcal{U}_{R}(G)$ sur la structure uniforme gauche $\mathcal{U}_{L}(G)$.
  \item La translation à gauche et la translation à droite sont uniformément continues pour les paires de structures uniformes: $(\mathcal{U}_{R}(G),\mathcal{U}_{R}(G))$ et $(\mathcal{U}_{L}(G),\mathcal{U}_{L}(G))$.
\end{enumerate}
\end{example}

\chapter{\quad Quelques théorèmes classiques}
\section{Théorème de Krein-Milman}\index{théorème!de Krein-Milman}
\begin{definition}
Soit $A$ une partie d'un espace localement convexe.
\begin{enumerate}
  \item Un sous-ensemble non-vide $M$ est un sous-ensemble extrême de $A$ si pour tout $x\in M$ tel que $x=\alpha y+ (1-\alpha)z$ avec $0<\alpha< 1$ et $y,z\in A$ alors $y,z\in M$
  \item Un point $x_{0}\in M$ est un point extrême de $M$ si $\{x_{0}\}$ est un sous-ensemble extrême de $M$. Nous noterons $Extr(M)$ l'ensemble des points extémaux de $M$.
  \item On notera $conv(A)=\{\underset{i=1}{\overset{n}{\sum}}\alpha_{i}x_{i}:\,\,\,\,x_{i}\in A,\,\,\alpha_{i}\geq 0,\,\,\,et\,\,\underset{i=1}{\overset{n}{\sum}}\alpha_{i}=1\}$ l'envellope convexe de $A$.
\end{enumerate}
\end{definition}
\begin{definition}(Topologie faible)
Soit $L$ un espace vectoriel sur $\mathbb{R}$ et soit $L^{\sharp}$ l'espace vectoriel de toutes les formes linéaires sur $L$. Considérons un sous-espace $F\subseteq L^{\sharp}$ qui sépare les points de $L$: Pour tout $x_{1}\neq x_{2}$, il existe $f\in F$ tel que $f(x_{1})\neq f(x_{2})$. La topologie faible sur $F$, notée $w(F)$ est la topologie dont une base de voisinages de $x_{0}$ est constitué des ensembles
$$N(x_{0};A;\varepsilon)=\{x\in L:\,\,|f(x)-f(x_{0})|<\varepsilon,\,\,\,f\in A\}$$ Où $x_{0}\in L,\,\,\varepsilon> 0$ et $A$ est une partie finie de $F$.\\En particulier:
\begin{enumerate}
  \item Si $L$ est un espace de Banach et $F=L^{\prime}$ est son dual topologiue, il s'agit de la topologie faible
  \item Si $L=X^{\prime}$ est le dual topologique d'un espace de Banach, il s'agit de la topologie vague.
\end{enumerate}
\end{definition}
\begin{theorem}(Krein-Milman \cite{milbook})\label{krein}
Soit $L$ un espace vectoriel et soit $K$ un ensemble compact et convexe de $(L,w(F))$. Alors $Extr(K)\neq\emptyset$ et $\overline{conv(Extr(K))}=K$.
\end{theorem}
\section{Théorème de Gelfand}\index{théorème!de Gelfand}
\begin{definition}
On appelle algèbre normée un couple $(A,\,\|.\|)$ où $A$ est une algèbre et $\|.\|$ est une norme sur $A$ qui vérifie en plus:
 $\|ab\|\leq \|a\|\,\|b\|$ pour tous $a,b\in A$
\end{definition}
\begin{definition}
Une algèbre normée\index{algèbre!normée} est dite de Banach si la norme $\|.\|$ est une norme complète.
\end{definition}
\begin{definition}
Soit $A$ une algèbre. Une involution sur $A$ est une application $\star:A\longrightarrow A$ vérifiant
pour tous $a,b\in A$ et $\alpha,\beta \in \mathbb{C} $
\begin{enumerate}
  \item $(\alpha a + \beta b)^{\star}=\overline{\alpha}a^{\star}+\overline{\beta}b^{\star}$
  \item $(a^{\star})^{\star}=a$
  \item $(ab)^{\star}=b^{\star}a^{\star}$
\end{enumerate}
\end{definition}
\begin{definition}
$A$ est une $C^{\star}$-algèbre si $A$ est une algèbre de Banach\index{algèbre!de Banach} muni d'une involution vérifiant:
$\|a^{\star}a\|=\|a\|^{2}$ pour tout $a \in A$.
\end{definition}

\begin{definition}

Soit $A$ une algèbre commutative sur un corps(commutatif) $\mathbb{K}$. On appelle caractère\index{caractère} ou fonctionnelle multiplicative de $A$ tout morphisme d'algèbres $\chi: A\longrightarrow \mathbb{K}$ non identiquement nul.

\end{definition}
\begin{lemma}
Soit $A$ une algèbre de Banach commutative et unifère. Alors le spectre de $A$ noté $X_A$ muni de la topologie de Gelfand est un espace compact.
\end{lemma}

\begin{theorem} (Gelfand)\label{gel}
Si $A$ est une $C^{\star}$-algèbre de Banach commutative et unifére, alors l'application
 \[ \xymatrix@C=2cm@R=0.3em{
 \mathcal{G}:A \ar[r] &   \mathcal{C}(X_{A}) &  \\
 a \ar@{|->}[r] & \mathcal{G}(a): X_{A} \ar[r] &  \mathbb{C} \\
& \,\,\,\,\,\,\,\,\,\,\,\,\,\,\,\,\,\,\,\,\,\,\,\,\,\,\,\,\,\chi \ar@{|->}[r] &  \mathcal{G}(a)(\chi)=\chi(a) }
\]
appelée transformation de Gelfand est un isomophisme isométrique de $A$ sur l'algèbre des fonctions continues sur l'espace compact $X_A$.
\end{theorem}
\begin{lemma}(Hahn-Banach\cite{inf})\label{hahnban}
Soit $F\subseteq E$ un sous-espace d'un espace localement convexe $E$ tel que $\overline{F}\neq E$, alors il existe une forme linéaire continue
$f\in E^{\prime},\,\,\,f\neq 0$ tel que $f(x)=0$ pour tout $x\in F$.
\end{lemma}

\bibliographystyle{plain}

\begin{thebibliography}{100}
\addcontentsline{toc}{1}{\bigskip \noindent  {\bf References
\hfill}}

\bibitem{youssef} Y. Al-Gadid: {\it  A new characterization of topologically amenable groups,} Thesis
of Master of Science in Mathematics of University of Ottawa, 2007
\bibitem{brice1} Y. Al-Gadid, B. Mbombo, and V. Pestov. {\it Sur les espaces test pour la moyennabilité,} C. R. Math. Rep. Acad. Sci. Canada, 33 :65--77, 2011.

\bibitem{inf} C. Aliprantis and K. Border. {\it Infinite dimensional Analysis}. Springer-Verlag.

\bibitem{anan} Anantharaman-Delaroche. {\it Amenability and exactness for dynamical systems and their C$^{\star}$-algebras,} Trans. Amer. Soc. {\bf 10} (2002), 4153--4178.


\bibitem{} R. Arens. {\it Topologies for homéomorphism groups,} In M. Hazewinkel, editor,
Handbook of Algebra, volume 2, pages 1--77. North-Holland, 2000.
\bibitem{atkin} C. J. Atkin. {\it Boundedness in uniform spaces, topological groups, and homogeneous spaces,} Acta Math. Hungar {\bf 57} (1991), 213--232.

\bibitem{aus} J. Auslander. {\it Minimal flows and their extensions,} North Holland, 1988.

 \bibitem{bantark} S. Banach and A. Tarski. {\it Sur la décomposition des ensembles de points en parties
respectivement congruents,} Fund. Math. {\bf 6} (1924), 244--277.

 \bibitem{becker} H. Becker and A. S. Kechris. {\it The Descriptive Set theory of Polish Group Actions,} London Math. Soc. Lecture Note Series 232, Cambridge Univ. Press, 1996.

\bibitem{bekka} B. Bekka and P. de la Harpe. {\it Kazhdan's Property (T),} Prepublication.
 \bibitem{BYBHU} I. Ben Yaacov, A. Berenstein, C.W. Henson, and A. Usvyatsov, {\em Model Theory for Metric Structures,} in: Model theory with applications to algebra and analysis. Vol. 2, 315--427,
London Math. Soc. Lecture Note Ser., \textbf{350}, Cambridge Univ. Press, Cambridge, 2008.

 \bibitem{berg} G. M. Bergman. {\it Generating infinite symmetric groups,} Bull. London Math. Soc. {\bf 23} (2006), 429--440.


\bibitem{bes} C. Bessaga and A. Pelczynski, {\it Selected Topics in Infinte-Dimensional Topology,} PWN, Warszawa (1975).

\bibitem{boga} S. Bogatyi and V.V. Fedorchuk, {\it Schauder's fixed point theorem and amenability of a group,} Topological Methods in Nonlinear Analysis, Journal of the Juliusz Schauder center {\bf 29} (2007), 383--401.

 \bibitem{bourb1} N. Bourbaki, {\it El\'ements de Math\'ematiques: Espaces vectoriels topologiques,} Hermann, Paris, 1958-1961.

 \bibitem{Bou} N. Bourbaki, {\it El\'ements de Math\'ematiques: Livre III, Topologie G\'en\'erale,} Hermann, Paris, 1958-1961.

\bibitem{bour} N. Bourbaki, {\it Int\'egration}, Hermann, Paris, 1963.

\bibitem{tro} A. Bouziad and J.-P. Troallic. {\it Problems about the uniform structures of topological
groups, in: Open problems in topology. II, Edited by Elliott Pearl, Elsevier B.V.,
Amsterdam,,} 2007, 359--366.

\bibitem{brook} R. B. Brook, {\it A construction of the greatest ambit,} Mathematical Systems Theory {\bf 4} (1970), 243--248.

\bibitem{BO} N.P. Brown and N. Ozawa, {\em ${C}^*$-Algebras and Finite-Dimensional Approximations,} Graduate Studies in Mathematics \textbf{88},  American Mathematical Society, Providence, R.I., 2008.

\bibitem{cho} G. Choquet. {\it Lectures on analysis. Vol. I,} New York-Amsterdam, 1969.
\bibitem{CR} W.W. Comfort and K.A. Ross, {\em Pseudocompactness and uniform continuity in topological groups,} Pacific J. Math. \textbf{16} (1966), 483--496.
\bibitem{const} C. Constantinescu. {\it $C^{\star}$-algebras, Vol. 2: Banach Algebras and Compact Operators,} North Holland, 2001.

\bibitem{dlH}  P. de la Harpe, {\it Moyennabilit\'e de quelques groupes
topologiques de dimension infinie,} C.R. Acad. Sci. Paris, S\'er.
{\bf A 277} (1973), 1037--1040.

 \bibitem{vries}  J. de Vries, {\it Elements of Topological Dynamics,} Kluwer, 1993.
\bibitem{dran} A. N. Dranishnikov, {\it On generalized amenability,} Preprint {\bf } (1999).
 \bibitem{dun} N. Dunford and J.T. Schwartz, {\it Linear operators. Part I. General theory, Pure and applied mathematics,} Wiley Classics Library, New York, 1988.

 \bibitem{milbook} Y. Eidelman, V. Milman and A. Tsolomitis, {\it Functional Analysis: An Introduction,} AMS (2004).

\bibitem{el} R. Ellis, {\it Universal Minimal Sets,} Proc. of the Amer. Math. Soc. {\bf 11} (1960), 540--543.

\bibitem{eng} R. Engelking, {\it General Topology,} Publishers Warsaw, 2006.

\bibitem{ey} P. Eymard, {\it Initiation à la théorie des groupes moyennables}.
\bibitem{FS} I. Farah and S. Shelah, {\em A dichotomy for the number of ultrapowers,} arXiv:0912.0406v1 [math.LO].
\bibitem{Fl1} J. Flood, {\it Free Topological Vector Spaces,} Ph.D. thesis, Australian National University, Canberra, 1975, 109 pp.

\bibitem{Fl2} J. Flood,
{\it Free locally convex spaces,} Dissertationes Math.
{\bf CCXXI} (1984), PWN, Warczawa.
\bibitem{frem} D.H. Fremlin, {\it Measure theory, vol. 4,} Torres-Fremlin, 2003.


\bibitem{kei} K. Funano, {\it Concentration of maps and group actions,} Geom. Dedicata. {\bf 149} (2010), 103--119.

\bibitem{gio} T. Giordano and P. de la Harpe, {\it Moyennabilit\'e des groupes d\'enombrables et actions sur les ensembles de Cantor,} C.R.Acad.Sci.Paris S\'er. {\bf 324} (1997), 1255--1258.

\bibitem{giopes1} T. Giordano and V.G. Pestov, {\it Some extremely amenable groups,} C.R. Acad. Sci. Paris. S\'er. I {\bf 4} (2002), 273--278.

\bibitem{giopes2} T. Giordano and V.G. Pestov, {\it Some extremely amenable groups related to operator algebras and ergodic theory,} J. Inst. Math. Jussieu {\bf 6} (2007), 279--315.


\bibitem{gla3} E. Glasner, {\it On minimal actions of Polish groups,} Top. Appl. {\bf 85} (1998), 119--125.

\bibitem{granirer} E. Granirer, {\it Extremely amenable semigroups I,} Math. Scand. {\bf 17} (1965), 177--179.
 \bibitem{gro} M. Gromov, {\it Metric Structures for Riemannian and Non-Riemannian Spaces,} Progress in Mathematics 152, Birkhauser Verlag, 1999.

 \bibitem{grommil} M. Gromov and V.D. Milman, {\it A topological application of the isoperimetric inequality,} Amer. J. Math. {\bf 105}(1983), 843--854.
 \bibitem{guran}
I.I. Guran, {\em Topological groups similar to Lindel\"of groups} (Russian), Dokl. Akad. Nauk SSSR \textbf{256} (1981), no. 6, 1305--1307.

\bibitem{hafner} P. Hafner and G. Mazzola, {\it The cofinal character of uniform spaces and ordered field,} Z. Math. Log. Grund. Math. {\bf 17} (1971), 377--384.


\bibitem{halmos} P. R. Halmos, {\it Naive set theory,} Undergraduate Texts in Mathematics. Springer-Verlag, New York, 1974. Reprint of the 1960 edition.

\bibitem{hewitt} E. Hewitt, {\it A remark on density characters,} Bull. Amer. Math. Soc. {\bf 52}
(1946), 641--643.

\bibitem{hewit} E. Hewitt and K. Ross, {\it Abstract harmonic analysis. Vol. 1,} Springer-Verlag, Berlin, Heidelberg, New York, 1979.

\bibitem{higson} N. Higson and J. Roe, {\it Amenable groups actions and the Novikov conjecture,} J. reine angew. Math {\bf 519} (2000), 143--153.

\bibitem{hu} G.E Huhunaisvili, {\it On a property of Uryson universal metric space(en russe),} Dokl. Akad. Nauk. USSR (N. S) {\bf 101} (1955), 332--333.
\bibitem{itz3}
G. Itzkowitz, {\em Continuous measures, Baire category, and uniform continuity in topological groups,} Pacific J. Math. \textbf{54} (1974), 115--125.

\bibitem{itz} G. Itzkowitz, {\it Uniformities and uniform continuity on topological groups,} General topology and its Applications, {\bf 134}
(1991), 155--178.

\bibitem{itz1} G. Itzkowitz, {\it Projective limits and balanced topological groups,} General topology and its Applications, {\bf 110}
(2001), 163--183.
\bibitem{kac}
G.I. Kac, {\em Isomorphic mapping of topological groups into a direct product of groups satisfying the first countability axiom} (Russian), Uspehi Matem. Nauk (N.S.) \textbf{8} (1953). no. 6(58), 107--113.
\bibitem{kal} N. Kalton, {\it Extending Lipschitz maps in C(K) spaces,} preprint, {\bf 2005}
(2001), 163--183.

\bibitem{kat} M. Kat\v etov, {\it On universal metric spaces,} in: Gen. Topology and its Relations to Modern Analysis and Algebra VI, Proc. Sixth Prague Topol. Symp. 1986, Z. Frolik, ed., Heldermann Verlag (1988), 323--330.
\bibitem{kechris} A. S. Kechris, {\it Classical descriptive set theory,} Springer-Verlag (1995).


\bibitem{kpt} A.S. Kechris, V.G. Pestov, S. Todorcevic, {\it Fraiss\'e limits, Ramsey theory, and topological dynamics of automorphism groups.,} GAFA, Geom. Funct. Anal. {\bf 15} (2005), 106--189.


\bibitem{KR} A.S. Kechris and C. Rosendal,
\textit{Turbulence, amalgamation and generic automorphisms of
homogeneous structures,} Proc. Lond. Math. Soc. (3) \textbf{94} (2007), 302--350.

\bibitem{keller} O. H. Keller, {\it Die Homeomorphie der kompakten konvexen Mengen in Hilbertschen Raum,} Math. Ann. {\bf 14} (1931), 748--758.


\bibitem{krieger} F. Krieger, {\it Sur les invariants topologiques des actions de Groupes Moyennables Discrets,} thèse de doctorat de l'université de Strasbourg.

\bibitem{ledoux} M. Ledoux, {\it The concentration of measure phenomenon,} Math. Surveys and Monographs, $89$, Amer. Math. Soc., $2001$.
\bibitem{lieberman}
A. Lieberman, {\em The structure of certain unitary representations of infinite symmetric groups,} Trans. Amer. Math. Soc. \textbf{164} (1972), 189--198.
\bibitem{LANVT}
J. Lopez-Abad and L. Nguyen Van Th\'e, {\em The oscillation stability problem for the Urysohn sphere: a combinatorial approach,} Topology Appl. \textbf{155} (2008), 1516--1530.
\bibitem{ulam} R. D. Mauldin, {\it The Scottish Book,} Boston-Basel-Stuttgart 1981.


\bibitem{mau} B. Maurey, {\it Constructions des suites symétriques,} C. R. Acad. Sci. Paris, Sér. A. {\bf 288} (1979), 679--681.

\bibitem{brice2} B. Mbombo and V. Pestov, {\it Subgroups of isometries of Urysohn-Kat\v etov
metric spaces of uncountable density,} Accepted to appear in Topology and its Applications

\bibitem{megrepest} M. Megreslishvili, P. Nickolas and V. Pestov, {\it Uniformities and uniformly continuous functions on locally connected groups,} Bull. Austral. Math. Soc. {\bf 56} (1997), 279--283.

\bibitem{megre} M. Megrelishvili and T. Scarr, {\it The Equivariant Universality and Couniversality of the Cantor Cube,} Fundamenta Mathematicae

 \bibitem{mel}
J. Melleray, {\em Géométrie de l'espace d'urysohn et théorie descriptives des ensembles, thèse de doctorat,} l'université Paris 6, Decembre 2005.


 \bibitem{melleray}
J. Melleray, {\em On the geometry of Urysohn's universal metric space,} Topology Appl. \textbf{154} (2007), 384--403.

\bibitem{mel2}
J. Melleray, {\em Topology of the isometry group of the Urysohn space,} Fund. Math. \textbf{207} (2010), 273--287.


  \bibitem{grommil1}
V. D. Milman and G. Schechtman, {\em Asymptotic Theory of Finite Dimensional Normed Spaces,} Lecture Notes in Math. $1200$. Springer, $1986$.

\bibitem{morris}
S.A. Morris, {\em Varieties of topological groups (a survey),} Colloq. Math. \textbf{XLVI} (1982), 147--165.

 \bibitem{munkres}
J. Munkres, {\em Topology, Second Edition} Prentice Hall, Inc. 1975.


\bibitem{von}
J. V. Neumann, {\em Zur allgemeinen Theorie des Masses,} Fund. Math. \textbf{13} (1929), 73--116.
\bibitem{NVT} L. Nguyen Van Th\'e, \textit{Structural Ramsey Theory of Metric Spaces and Topological Dynamics of Isometry Groups}, Mem. Amer. Math. Soc. \textbf{206} (2010), no. 968.

    \bibitem{NVTS}
L. Nguyen Van Th\'e and N.W. Sauer, {\em The Urysohn sphere is oscillation stable,} Geom. Funct. Anal. \textbf{19} (2009), 536--557.
    \bibitem{OS}
E. Odell and T. Schlumprecht, \textit{The distortion problem,}
Acta Math. {\bf 173} (1994), 259--281.
\bibitem{ozawa} N. Ozawa, {\it Amenable actions and exactness for discrete groups.,} C. R. Acad. Sci. Paris S\'er. I Math {\bf 08} (2000), 691--695.

\bibitem{paterson} A.T. Paterson, {\it Amenability,} University Math. Surveys and Monographs 29, Amer. Math. Soc., Providence, RI, 1988.

\bibitem{vp5} V.G. Pestov, {\em Embeddings and condensations of topological groups} (Russian), Mat. Zametki \textbf{31} (1982), no. 3, 443--446.

\bibitem{vp1} V. Pestov, {\em Dynamics of Infinite-Dimensional Groups: the Ramsey-Dvoretzky-Milman phenomenon,} Amer. Math. Soc. University Lecture Series \textbf{40}, 2006.

 \bibitem{vp2} V.G. Pestov, {\it On free actions, minimal flows, and a problem by Ellis,} Trans. of the American Mathematical Society {\bf 350} (1998), 4149--4165.

\bibitem{vp4} V.G. Pestov, {\em Ramsey-Milman phenomenon, Urysohn metric spaces, and extremely amenable groups,} Israel Journal of Mathematics \textbf{127} (2002), 317--358. Corrigendum, ibid., \textbf{145} (2005), 375--379.
\bibitem{vpbresil} V.G. Pestov, {\it Dynamics of Infinite-dimensional Groups and Ramsey-type phenomena,} Impa, Brazil 2005.

\bibitem{pier} J-P. Pier, {\it Amenable Locally Compact groups,} New York, 1984.

\bibitem{pondi} E.S. Pondiczery, {\it Power problems in abstract spaces,} Duke Math. J. {\bf 11} (1944), 835--837.

\bibitem{prota2}I.V. Protasov, {\em Functionally balanced groups,} Math. Notes \textbf{49} (1991), 614--616.

\bibitem{pro} I.V. Protasov and A. Saryev, {\em The semigroup of closed subsets of a topological group} (Russian),
Izv. Akad. Nauk Turkmen. SSR Ser. Fiz.-Tekhn. Khim. Geol. Nauk 1988, no. 3, 21--25.

\bibitem{dierolf} W. Roelcke and S. Rierolf, {\em Uniform Structure on Topolgical Group and Their Quotient}, McGraw-Hill, $1981$.

\bibitem{ros} C. Rosendal, {\it A topological version of the Bergman property,} Forum Math. {\bf 21} (2009), 299--332.

\bibitem{ru} W. Rudin, {\it Analyse réelle et complexe, troisième édition,} Dunod, Paris, 1998.

 \bibitem{samuel} P. Samuel, {\it Ultrafilters and compactification of uniform spaces,} Trans. Amer. Math. Soc. {\bf 64} (1948), 100--132.

\bibitem{schauder} J. Schauder, {\it Zur Theorie stetiger abbildungen in Funktionalräumen,} Math. Z. {\bf 26} (1927), 47--65.

\bibitem{si} W. Sierpinski, {\it Sur un espace métrique universel,} Fund. Math. {\bf 33} (1945), 115--122.

\bibitem{tel} S. Teleman, {\it Sur la representation linéaire des groupes topologiques,} Ann. Sci. Ecole Norm. Sup. {\bf 74} (1957), 319--339.

\bibitem{ury1} P. S. Urysohn, {\it Sur un espace métrique universel,} C. R. Acad. Sci. Paris {\bf 180} (1925), 803--806.

\bibitem{ury2} P. S. Urysohn, {\it Sur un espace métrique universel,} Bull. Sci. Math. {\bf 52} (1927), 43--64 et 74--90.

\bibitem{us4} V.V. Uspenskij, {\it A universal topological group with countable base,} Funct. Anal. Appl. {\bf 20} (1986), 160--161.

\bibitem{us5} V.V. Uspenskij, {\it On the group of isometries of the Urysohn universal metric space,} Comment. Math. Univ. Carolinae {\bf 31} (1990), 181--182.

\bibitem{us6} V.V. Uspenskij,
{\em On subgroups of minimal topological groups,}
Topology Appl. \textbf{155} (2008), 1580--1606.

\bibitem{veech} W. A. Veech, {\it Topological dynamics,} Bull. Amer. Math. Soc. {\bf 83} (1977), 775--830.

  \bibitem{will} S. Willard, {\it General topology,} Addison-Wesley Publishing Company,Inc., 1970.

\end{thebibliography}

\printindex
\end{document}